\documentclass[preprint,12pt]{elsarticle}

\usepackage[utf8]{inputenc}
\usepackage[T1]{fontenc}
\usepackage{lmodern}
\usepackage{amsmath,amssymb,amsthm,mathtools}
\usepackage{geometry}
\usepackage{microtype}
\AtBeginDocument{%
  \setlength{\abovedisplayskip}{4pt plus 1pt minus 2pt}%
  \setlength{\belowdisplayskip}{4pt plus 1pt minus 2pt}%
  \setlength{\abovedisplayshortskip}{2pt plus 1pt minus 1pt}%
  \setlength{\belowdisplayshortskip}{2pt plus 1pt minus 1pt}%
}
\usepackage{caption}
\allowdisplaybreaks
\usepackage{enumitem}
\setlist{itemsep=1pt,topsep=2pt,parsep=0pt,partopsep=0pt}
\usepackage{graphicx}
\graphicspath{{figures/}}
\usepackage{booktabs}
\usepackage{array}
\usepackage[dvipsnames]{xcolor}

\usepackage{aliascnt}

\theoremstyle{plain}
\newtheorem{theorem}{Theorem}[section]
\newaliascnt{lemma}{theorem}
\newtheorem{lemma}[lemma]{Lemma}
\aliascntresetthe{lemma}
\newaliascnt{proposition}{theorem}
\newtheorem{proposition}[proposition]{Proposition}
\aliascntresetthe{proposition}
\newaliascnt{corollary}{theorem}
\newtheorem{corollary}[corollary]{Corollary}
\aliascntresetthe{corollary}
\newaliascnt{fact}{theorem}
\newtheorem{fact}[fact]{Fact}
\aliascntresetthe{fact}
\theoremstyle{definition}
\newaliascnt{definition}{theorem}
\newtheorem{definition}[definition]{Definition}
\aliascntresetthe{definition}
\newaliascnt{assumption}{theorem}
\newtheorem{assumption}[assumption]{Assumption}
\aliascntresetthe{assumption}
\theoremstyle{remark}
\newaliascnt{remark}{theorem}
\newtheorem{remark}[remark]{Remark}
\aliascntresetthe{remark}

\usepackage{hyperref}
\hypersetup{colorlinks=true, linkcolor=NavyBlue,
  citecolor=NavyBlue, urlcolor=NavyBlue}
\makeatletter
\@for\@tmpcnt:=theorem,lemma,proposition,corollary,fact,definition,assumption,remark\do{%
  \expandafter\edef\csname theH\@tmpcnt\endcsname{\noexpand\theHsection.\noexpand\arabic{\@tmpcnt}}}
\makeatother
\usepackage[capitalise,nameinlink]{cleveref}

\crefname{theorem}{Theorem}{Theorems}
\crefname{lemma}{Lemma}{Lemmas}
\crefname{proposition}{Proposition}{Propositions}
\crefname{corollary}{Corollary}{Corollaries}
\crefname{definition}{Definition}{Definitions}
\crefname{remark}{Remark}{Remarks}
\crefname{fact}{Fact}{Facts}
\crefname{assumption}{Assumption}{Assumptions}

\DeclareMathOperator{\ran}{Ran}
\DeclareMathOperator{\ind}{ind}
\DeclareMathOperator{\Real}{Re}
\DeclareMathOperator{\Imag}{Im}
\newcommand{\R}{\mathbb{R}}
\newcommand{\norm}[1]{\left\lVert #1\right\rVert}
\newcommand{\opnorm}[1]{\left\lVert #1\right\rVert_{\mathrm{op}}}
\newcommand{\ip}[2]{\left\langle #1,#2\right\rangle}
\newcommand{\Fh}{F_{\#}}
\newcommand{\dd}{\,\mathrm d}
\newcommand{\LambertW}{\mathcal W_0}

\begin{document}

\begin{frontmatter}

\title{\textbf{Source Conditions and Sharp Transition Asymptotics for the
Factorization Method with Far-Field Data}}
\author{George Pelekanos\textsuperscript{$\ast$} and Koung Hee Leem}
\cortext[cor1]{Corresponding author. \emph{Email address:}
  \href{mailto:gpeleka@siue.edu}{\texttt{gpeleka@siue.edu}} (G.~Pelekanos).}
\affiliation{organization={Southern Illinois University},
            addressline={Department of Mathematics and Statistics},
            city={Edwardsville},
            postcode={62026},
            state={IL},
            country={USA}}

\begin{abstract}
The factorization method characterizes an obstacle through a range
test, but practical implementations replace the associated Picard series by a
regularized approximation.  This paper studies the resulting regularization
bias for inverse acoustic scattering by impedance obstacles.  For a disc, we
show that every positive-order H\"older source condition fails at each
noncentral interior sampling point, whereas logarithmic source conditions hold
uniformly on compact interior subsets for every exponent, with explicit
constants.  At the center, the bias decays linearly with the regularization
parameter.  We then derive sharp spectral transition asymptotics through the
Lambert W function.  These asymptotics determine the stabilization rate inside the obstacle, the
iterated-logarithmic growth at the boundary, and the divergence outside.  They
also yield a threshold rule for which the reconstructed disc converges in
measure, together with quantitative boundary-location estimates.
For general impedance obstacles, we reduce the spectral decay of the
symmetrized far-field operator to that of the data-to-pattern operator and
obtain two-sided transition bounds in terms of the inradius and any radius
larger than the circumradius.  For a simply connected obstacle with analytic
boundary, positive analytic impedance, and conformal center inside the
obstacle, a Faber--Vekua normal form gives the sharp spectral law in terms of
logarithmic capacity, with the same leading constant as for the disc.  This
class includes every convex analytic obstacle.  Numerical experiments confirm
the predicted spectral transition and illustrate its use in choosing the
regularization level.  The analytic normal-form result is proved only for a
single boundary component.
\end{abstract}

\begin{keyword} factorization method, inverse obstacle scattering,
source conditions, regularization rates, Lambert $W$-function,
impedance boundary condition, severely ill-posed problems
\end{keyword}

\end{frontmatter}

\section{Introduction}

Qualitative methods in inverse scattering \cite{CK19} identify an obstacle
through a range test rather than through an optimization over shapes.  For the
factorization method, one forms, from the far-field operator $F$, a
positive compact operator $\Fh$ whose square root has the same range as the
data-to-pattern operator of the scatterer, so that a point $z$ lies in the
obstacle $D$ exactly when the sampling function
$\varphi_z(\widehat x)=e^{-ik\widehat x\cdot z}$ belongs to
$\ran(\Fh^{1/2})$ \cite{Kir98,KG08}.  In the two-dimensional impedance
setting considered here, $\Fh=|\Real(e^{-i\pi/4}F)|+\Imag(e^{-i\pi/4}F)$.
The phase rotation is forced by the two-dimensional Hankel asymptotics
and is already in Kirsch~\cite[Section~5]{Kir98}, where the range theorems are
shown to carry over to two dimensions with $e^{-i\pi/4}F$ in place of
$F$.  For an impedance boundary, $F$ is generally nonnormal, so $|F|$
must be replaced by the symmetrized $\Fh$ above.  The range
identity in that form is (F2) of \cref{fact:fm-inputs}.

Membership is decided in practice by the Picard series.  With
$(\gamma_j,\psi_j)_{j\ge1}$ an orthonormal eigensystem of $\Fh$,
\begin{equation}\label{eq:intro-picard}
 z\in D
 \quad\Longleftrightarrow\quad
 W(z):=\sum_{j\ge1}\frac{|\ip{\varphi_z}{\psi_j}|^2}{\gamma_j}<\infty,
\end{equation}
and since $W$ divides by arbitrarily small eigenvalues it is replaced by
the filtered field
\begin{equation}\label{eq:intro-Walpha}
 W_\alpha(z)=\ip{\varphi_z}{\Fh(\Fh+\alpha I)^{-2}\varphi_z}
 =\sum_{j\ge1}\frac{\gamma_j}{(\gamma_j+\alpha)^2}|\ip{\varphi_z}{\psi_j}|^2,
 \qquad \alpha>0 .
\end{equation}
The filter increases pointwise to $1/t$ as $\alpha\downarrow0$, so
$W_\alpha\uparrow W$.  The dichotomy \eqref{eq:intro-picard} is therefore
recovered in the limit.  That regularized indicators of this kind do
recover the range test in the limit was first established, for the linear
sampling method, through the same factorization structure
\cite{Are04,AL09}.  Quantitative conclusions about the method (how
small $\alpha$ must be, how the imaging error is distributed between
regularization and data error, whether an imaging experiment at fixed
dimension means anything) depend on the \emph{rate} at which
$W(z)-W_\alpha(z)$ vanishes at interior points.  That rate is the subject
of this paper.

\subsection*{Source conditions for the impedance disc}

Rates for a filtered spectral method are not available without a
hypothesis on how the right-hand side sits relative to the spectrum of the
operator.  The classical device is a H\"older source condition.  With
$A:=\Fh$, one assumes $\varphi_z=A^{1/2+s}q_z$ for some $s>0$ and a
bounded family $q_z$, and obtains a bias of order $\alpha^{2s}$
\cite[Theorem~4.3 and Corollary~4.4]{EHN96}.

Our first result is that for the impedance disc these conditions have no content.  The
eigenvalues of $\Fh$ for a smooth obstacle decay super-exponentially.  If
the coefficients $\ip{\varphi_z}{\psi_j}$ decay only geometrically
relative to $\gamma_j^{1/2}$, the quotient defining a H\"older source
condition diverges for every positive exponent.  For the impedance disc both
behaviors, and the failure itself, can be proved outright.  For every $s>0$ and every interior point
other than the center, $\varphi_z\notin\ran(A^{1/2+s})$
(\cref{prop:hoelder-fails}).  The center is the only exception.
Only one Fourier coefficient of $\varphi_0$ is nonzero, and its bias is
linear in $\alpha$ (Remark~\ref{rem:center}).

A common replacement for severely ill-posed problems is a logarithmic
source condition \cite{Hoh97,MP03,Tau98}, here
$\varphi_z=A^{1/2}\ell_\theta(A)q_z$ with
$\ell_\theta(t)=(1+\log(\norm A/t))^{-\theta}$.  We prove the
corresponding bias estimate with explicit constants
(\cref{thm:log-rate}) and show that the impedance disc satisfies the
condition on every compact interior set, for \emph{every} exponent
$\theta>0$ (\cref{prop:disc-log}).  The resulting rates are logarithmic in
$1/\alpha$ and therefore slow.  Since the disc satisfies the condition for
every exponent, its interior bias in fact decays faster than every fixed
logarithmic power, and \cref{sec:disc-sharp} identifies the exact
transition behind this.  Logarithmic conditions of
variational type have meanwhile been verified for inverse medium
scattering \cite{HW15}.  Remark~\ref{rem:vsc} relates that theory to the
spectral form required here.

\subsection*{The transition index}

Bias estimates of this kind are upper bounds, and they say little about
which part of the spectrum is responsible.  For the disc, the explicit
Fourier representation permits sharp estimates.  The filter in
\eqref{eq:intro-Walpha} passes modes with $\gamma_m\gtrsim\alpha$ and
suppresses the rest, so the imaging behavior is controlled by the largest
index still passed,
\[
 m_\alpha=\max\{m:\gamma_m\ge\alpha\}.
\]
Solving $\gamma_m\approx\alpha$ against the sharpened eigenvalue
asymptotics leads to an equation of the form $x\log x=\text{const}$, whose
inverse is the Lambert $W$-function \cite{CGHJK96}.  This yields
\cref{thm:disc-sharp},
\[
 m_\alpha=\frac{L_\alpha}{2\LambertW\bigl(L_\alpha/(eka)\bigr)}+O(1),
 \qquad L_\alpha=\log\frac{C_k}{\alpha},\quad C_k=\sqrt{2\pi/k}.
\]
The number of usable modes therefore grows like a logarithm of $1/\alpha$
divided by a log-logarithm because the spectrum decays super-exponentially.
The same theorem shows that this one index controls three
different behaviors.  Interior points stabilize with
$\log(W(r)-W_\alpha(r))=2m_\alpha\log(r/a)+O(\log m_\alpha)$, boundary
points grow like $\log m_\alpha$, and exterior points diverge with
$\log W_\alpha(r)=2m_\alpha\log(r/a)+O(\log m_\alpha)$.  The proof rests
on a two-sided bound for the ratio of adjacent eigenvalues, sharpened to a
relative $O(m^{-1})$ (\cref{lem:gamma-ratio}).  A bracket rather than an
asymptotic equivalence is what is available for $\gamma_{m_\alpha}$
itself, since consecutive eigenvalues differ by factors of order $m^{-2}$.

\subsection*{Reconstruction and numerical implications}

For set reconstruction, the obstacle approximation
used here is a \emph{sublevel} set of $W_\alpha$, and no fixed
threshold can classify every interior point.  The threshold must grow as
the regularization decreases.  We give a diagonal choice under which the
sublevel sets converge to $D$ in measure (\cref{thm:set-consistency}), and
show that a transversality margin at the reconstructed level curve
converts a uniform field error into a Hausdorff bound
(\cref{prop:hausdorff}).  For the disc, the transition asymptotics make
all of this quantitative.  The correct threshold scale is the boundary
value itself,
$W_\alpha(a)=\tfrac{2}{C_k}\log m_\alpha+O(1)$, with leading behavior
$\sqrt{2k/\pi}\,\log\log(1/\alpha)$ depending on the wavenumber alone
(\cref{cor:boundary-constant}).  The additive threshold rule
$\tau_\alpha=\tfrac{2}{C_k}\log m_\alpha+\beta$ then recovers $D$ with
misclassified area
$O(1/m_\alpha)=O\bigl(\log\log(1/\alpha)/\log(1/\alpha)\bigr)$
(\cref{thm:disc-recon}), and the margin hypothesis of the Hausdorff
estimate holds under the same threshold rule, locating the reconstructed level
circle within $O(a/m_\alpha)$ of the boundary
(\cref{prop:disc-margin}).

For numerical practice, at a fixed Fourier
truncation the regularized field has a finite limit as $\alpha\downarrow0$
even at exterior points, so large pre-asymptotic growth in a
fixed-dimensional experiment is not evidence of the divergence that the
range test predicts.  Only a coupled limit (regularization decreasing
while the retained spectrum grows) reproduces the continuum behavior, as
the computations of \cref{sec:disc-num} show.

\subsection*{The range identity}

The range identity underlying the method is
Kirsch and Grinberg~\cite[Theorem~2.15]{KG08}, applied to the impedance condition in
Kirsch and Grinberg~\cite[Section~2.5.2]{KG08} and carried to two dimensions by the phase
shift of Kirsch~\cite[Section~5]{Kir98}.  \Cref{fact:fm-inputs} states the form
used here, together with the normalization transfer that carries it
from the standard far-field convention to the physical one.  The finite-data perturbation
estimate used in several conditional statements is isolated in the next
subsection.  \Cref{sec:general} proves a two-sided transition law
valid for every impedance obstacle, bracketed
below by the inradius and above by any radius exceeding the
circumradius, and, for a single analytic boundary component whose
conformal center lies in the obstacle, that the disc law holds with the
logarithmic capacity in place of the radius
(\cref{thm:capacity-connected}).  The analytic normal-form argument is
restricted to a single boundary component.

\subsection*{The finite-data input}
A few results below combine the regularization bias with the error of
forming the field from finitely many noisy far-field samples.  That
error is not analysed here; it enters only through
\eqref{eq:data-bound}.

Write $\widehat W_{\alpha,M,N}$ for the empirical regularized field
formed from $M$ observation directions and $N$ averaged measurement
frames, $R_K:=\norm{\varphi_z}=\sqrt{2\pi}$ (the same for every $z$), and let
$\zeta_M$ be the test-vector approximation rate, $e_{M,N}$ the resulting
operator error, and $\rho_\Lambda$ the modulus of continuity of the
symmetrization $F\mapsto\Fh$ on the ball of radius $\Lambda$.  Its only
nonlinear ingredient is the absolute value of a self-adjoint operator,
which is continuous but \emph{not} operator Lipschitz.  What holds
instead is a bound with an iterated logarithm \cite{AP11}, improving
the single logarithm of Kato~\cite{Kato73}: with an
absolute constant $C$ and $t_\Lambda=\min(t,2\Lambda)$,
\begin{equation}\label{eq:rho-modulus}
 \rho_\Lambda(0):=0,
 \qquad
 \rho_\Lambda(t)\le C\,t_\Lambda
 \log\Bigl(2+\log\frac{2\Lambda}{t_\Lambda}\Bigr),
 \qquad t>0.
\end{equation}
Only the consequence
$\rho_\Lambda(t)=O\bigl(t\log\log(1/t)\bigr)$ as $t\downarrow0$ is used
below.

We assume throughout those statements that there is an event on which,
for every $\alpha>0$,
\begin{equation}\label{eq:data-bound}
 \sup_{z\in K_0}\bigl|\widehat W_{\alpha,M,N}(z)-W_\alpha(z)\bigr|
 \le\frac{3(R_K+\zeta_M)^2}{\alpha^{2}}\,\rho_\Lambda(e_{M,N})
 +\frac{(2R_K+\zeta_M)\zeta_M}{4\alpha}.
\end{equation}

The first term is the operator error $e_{M,N}$ carried through the
symmetrization and amplified by the filter, the second the test-vector
error $\zeta_M$.  At fixed $\alpha$ both vanish as the data improve.
Only this shape is used below, not the particular constants.  Every statement below that mentions
$\widehat W_{\alpha,M,N}$ is conditional on it.

\section{Source conditions and regularization bias}\label{sec:rates}

\subsection{The impedance disc}\label{sec:disc-prelim}
The impedance disc is the model case throughout, because there every object
above is explicit.  Let $D$ be the disc of radius $a$ centered at the origin,
with constant impedance $\lambda>0$ and wavenumber $k>0$.  The far-field
operator is diagonal in the Fourier basis
$\psi_m(\widehat x(t))=(2\pi)^{-1/2}e^{imt}$,
$\widehat x(t)=(\cos t,\sin t)$, with scattering coefficients
\begin{equation}\label{eq:disc-am}
 a_m=-\frac{J_m'(ka)+i\lambda J_m(ka)}
 {H_m^{(1)\prime}(ka)+i\lambda H_m^{(1)}(ka)}
\end{equation}
and far-field eigenvalues $\kappa_m$ given in \eqref{eq:disc-kappa} below.  Throughout, the far field is normalized by
$u^{s}(x)=|x|^{-1/2}e^{ik|x|}\bigl(u^{\infty}(\widehat x)+O(|x|^{-1})\bigr)$,
without the prefactor $\gamma_2$ of the convention used in
Chamaillard et al.~\cite{CCH13}.  The two differ by $F=\gamma_2F_{\rm st}$, which is the
source of the rotation in \eqref{eq:disc-kappa}.  In this normalization
\begin{equation}\label{eq:disc-kappa}
 \kappa_m=2\pi\sqrt{\frac{2}{\pi k}}\,e^{-i\pi/4}a_m,
 \qquad
 \widetilde\kappa_m=e^{-i\pi/4}\kappa_m,
 \qquad
 \gamma_m=|\Real\widetilde\kappa_m|+\Imag\widetilde\kappa_m,
\end{equation}
the $\gamma_m$ being the eigenvalues of $A=\Fh$.  Since $2\pi\sqrt{2/(\pi k)}=\sqrt{8\pi/k}$, \eqref{eq:disc-kappa}
gives $\widetilde\kappa_m=-i\sqrt{8\pi/k}\,a_m$ and hence
$\Imag\widetilde\kappa_m=\sqrt{8\pi/k}\,(-\Real a_m)$.  For an
absorbing obstacle the energy bound $|1+2a_m|\le1$ \cite{CK19} is
equivalent to $-\Real a_m\ge|a_m|^{2}\ge0$, so
$\Imag\widetilde\kappa_m\ge0$ for every $m$.  The imaginary part of
$\widetilde F$ enters $\Fh$ unmodified and no clipping is required.  Since
$\gamma_{-m}=\gamma_m$, indices are taken nonnegative wherever convenient.

\subsection{Source conditions and bias rates}
Let $A=\Fh$.  For an interior sampling point, $\varphi_z\in\ran(A^{1/2})$,
but this alone gives no rate for $W_\alpha(z)\uparrow W(z)$.  A
standard spectral source condition yields an algebraic rate.

\begin{assumption}[Local source condition]\label{ass:source}
For a set $K_0\Subset D$, there are $s\in(0,1/2]$ and $C_q<\infty$ such that
for every $z\in K_0$,
\begin{equation}\label{eq:source}
 \varphi_z=A^{1/2+s}q_z,
 \qquad
 \norm{q_z}\le C_q.
\end{equation}
\end{assumption}

\begin{theorem}[Regularization-bias rate]\label{thm:source-rate}
Under Assumption~\ref{ass:source},
\begin{equation}\label{eq:bias-rate}
 0\le W(z)-W_\alpha(z)
 \le C_s C_q^2\alpha^{2s}
 \qquad(z\in K_0),
\end{equation}
where
\[
 C_s=\sup_{x\ge0}
 \frac{x^{2s}(2x+1)}{(x+1)^2}<\infty.
\]
Consequently, on the event of \eqref{eq:data-bound},
\[
\begin{aligned}
 \sup_{z\in K_0}|\widehat W_{\alpha,M,N}(z)-W(z)|
 &\le C_sC_q^2\alpha^{2s}
 +\frac{3(R_K+\zeta_M)^2}{\alpha^2}
 \rho_\Lambda(e_{M,N})\\
 &\quad+
 \frac{(2R_K+\zeta_M)\zeta_M}{4\alpha}.
\end{aligned}
\]
If $\zeta_M=O(\alpha^{2s+1})$, the third term is $O(\alpha^{2s})$ and
balancing the first two terms gives
\begin{equation}\label{eq:optimal-alpha}
 \alpha\asymp\rho_\Lambda(e_{M,N})^{1/(2s+2)},
 \qquad
 \sup_{K_0}|\widehat W-W|
 =O\!\left(\rho_\Lambda(e_{M,N})^{s/(s+1)}\right);
\end{equation}
at that $\alpha$ the hypothesis reads
$\zeta_M=O\bigl(\rho_\Lambda(e_{M,N})^{(2s+1)/(2s+2)}\bigr)$.
\end{theorem}

\begin{proof}
Let $E_A$ be the spectral measure of $A$.  From
\eqref{eq:source},
\begin{align*}
 W(z)-W_\alpha(z)
 &=\int_{(0,\norm A]}
 \left(\frac1t-\frac{t}{(t+\alpha)^2}\right)
 t^{1+2s}\,\dd\ip{E_A(t)q_z}{q_z}\\
 &=\int_{(0,\norm A]}
 \frac{\alpha(2t+\alpha)t^{2s}}{(t+\alpha)^2}
 \,\dd\ip{E_A(t)q_z}{q_z}.
\end{align*}
With $x=t/\alpha$, the scalar factor is
$\alpha^{2s}x^{2s}(2x+1)/(x+1)^2$.  As $x\to\infty$ this behaves like
$2x^{2s-1}$, bounded precisely when $s\le1/2$.  The restriction is thus
the qualification saturation of Tikhonov regularization rather than an
artefact of the estimate, and \eqref{eq:bias-rate} follows.  Add \eqref{eq:data-bound} and minimize the leading two terms to obtain
\eqref{eq:optimal-alpha}.
\end{proof}

\begin{proposition}[The H\"older source condition is empty for the disc]\label{prop:hoelder-fails}
Let $D$ be the disc of radius $a$ centered at the origin, with constant
impedance $\lambda>0$ and wavenumber $k>0$, and let $A=\Fh$.  Then for every
$s>0$ and every $z\in D$ with $0<|z|<a$,
\[
 \varphi_z\notin\ran(A^{1/2+s}).
\]
Consequently Assumption~\ref{ass:source} fails, for every $s\in(0,1/2]$, on
every compact $K_0\Subset D$ that contains a point other than the center.
\end{proposition}

\begin{proof}
Write $r=|z|\in(0,a)$.  On the disc, $A$ is diagonal in the Fourier basis
$\psi_m(\widehat x(t))=(2\pi)^{-1/2}e^{imt}$, $\widehat x(t)=(\cos t,\sin t)$,
with eigenvalues
$\gamma_m=|\Real\widetilde\kappa_m|+\Imag\widetilde\kappa_m$,
$\widetilde\kappa_m=e^{-i\pi/4}\kappa_m$ as in \eqref{eq:disc-kappa}.
The Jacobi--Anger expansion
\[
 e^{-ik\widehat x\cdot z}
 =\sum_{m\in\mathbb Z}(-i)^{m}J_m(kr)\,e^{im(t-t_z)},
\]
with $t_z$ the polar angle of $z$, gives
\begin{equation}\label{eq:disc-coefficients}
 |\ip{\varphi_z}{\psi_m}|^2=2\pi J_m(kr)^2 .
\end{equation}
For a positive compact self-adjoint injective operator with eigensystem
$(\gamma_m,\psi_m)$, the Picard criterion applied to $A^{1/2+s}$ states that
$\varphi_z\in\ran(A^{1/2+s})$ if and only if
\begin{equation}\label{eq:power-picard}
 \sum_{m\in\mathbb Z}
 \frac{|\ip{\varphi_z}{\psi_m}|^2}{\gamma_m^{1+2s}}<\infty .
\end{equation}

We show that the terms of \eqref{eq:power-picard} are unbounded.  By the
large-order asymptotics \cite[Section~10.19]{DLMF}, uniformly for $x$ in
compact subsets of $(0,\infty)$,
\[
 J_m(x)=\frac{(ex/2m)^{m}}{\sqrt{2\pi m}}\bigl(1+o(1)\bigr),
 \qquad
 Y_m(x)=-\sqrt{\frac{2}{\pi m}}\Bigl(\frac{ex}{2m}\Bigr)^{-m}
 \bigl(1+o(1)\bigr),
\]
while the recurrences $J_m'=(m/x)J_m-J_{m+1}$, $Y_m'=(m/x)Y_m-Y_{m+1}$ give
$J_m'=(m/x)J_m(1+o(1))$, $H_m^{(1)\prime}=-(m/x)H_m^{(1)}(1+o(1))$, and
$H_m^{(1)}=iY_m(1+o(1))$.  In $a_m$ of \eqref{eq:disc-am} the impedance terms
are therefore of lower order, and
\[
 a_m=\frac{J_m(ka)}{H_m^{(1)}(ka)}\bigl(1+o(1)\bigr)
 =\frac{i}{2}\Bigl(\frac{eka}{2m}\Bigr)^{2m}\bigl(1+o(1)\bigr).
\]  Hence
$\widetilde\kappa_m=-i\sqrt{8\pi/k}\,a_m
=\tfrac12\sqrt{8\pi/k}\,(eka/2m)^{2m}(1+o(1))$ is positive real to leading
order, and there are constants $0<c_1\le c_2$ and $m_0$ with
\begin{equation}\label{eq:gamma-twosided}
 c_1\Bigl(\frac{eka}{2m}\Bigr)^{2m}
 \le\gamma_m\le
 c_2\Bigl(\frac{eka}{2m}\Bigr)^{2m},
 \qquad m\ge m_0 .
\end{equation}
Combining \eqref{eq:disc-coefficients}, \eqref{eq:gamma-twosided}, and the
$J_m$ asymptotics, the logarithm of the $m$th term of
\eqref{eq:power-picard} is
\[
 2m\log\frac{ekr}{2m}-(1+2s)\,2m\log\frac{eka}{2m}+O(\log m)
 =2m\log\frac ra+4sm\log\frac{2m}{eka}+O(\log m).
\]
The first summand is linear in $m$ with negative slope, but for every $s>0$
the second grows like $4sm\log m$.  The sum tends to $+\infty$, so the terms
of \eqref{eq:power-picard} are unbounded and the series diverges.  (For
$s=0$ the display reduces to $2m\log(r/a)+O(\log m)\to-\infty$, giving
geometric convergence, which is the range identity.)  At the center $z=0$
only the $m=0$ coefficient in \eqref{eq:disc-coefficients} is nonzero, so
$\varphi_0\in\ran(A^{t})$ for every $t>0$.  This is why the center is
excluded.
\end{proof}

\begin{remark}[The center]\label{rem:center}
At the exceptional point $z=0$ only the mode $m=0$ contributes, and the
bias is linear in the regularization:
\[
 W(0)-W_\alpha(0)
 =2\pi\,\frac{\alpha(2\gamma_0+\alpha)}{\gamma_0(\gamma_0+\alpha)^{2}}
 =\frac{4\pi}{\gamma_0^{2}}\,\alpha+O(\alpha^{2}),
\]
consistent with $\varphi_0\in\ran(A^{t})$ for every $t>0$.  The center
satisfies every power source condition.  Its bias is nonetheless only
$O(\alpha)$, as the display shows, because Tikhonov saturates at
qualification one.  The center therefore has the rate $O(\alpha^{p})$
for $0<p\le1$ and no better.  By \eqref{eq:sharp-interior} of
Theorem~\ref{thm:disc-sharp}, the bias at every other interior point is
$\exp(-O(m_\alpha))$ and so decays more slowly than every positive power
of $\alpha$; the center is the only interior point with a power rate.
\end{remark}

The disc calculation motivates a logarithmic replacement but does not
settle general obstacles, so the next result is formulated conditionally,
its hypothesis being verified for the disc in \cref{prop:disc-log}.

\begin{assumption}[Logarithmic source condition]\label{ass:log-source}
For a set $K_0\Subset D$ there are $\theta>0$ and $C_q<\infty$ such that for
every $z\in K_0$,
\begin{equation}\label{eq:log-source}
\begin{gathered}
 \varphi_z=A^{1/2}\,\ell_\theta(A)\,q_z,
 \qquad
 \norm{q_z}\le C_q,
 \\
 \text{where}\quad
 \ell_\theta(t)=\Bigl(1+\log\frac{\norm A}{t}\Bigr)^{-\theta},
 \quad t\in(0,\norm A],
 \qquad \ell_\theta(0):=0,
\end{gathered}
\end{equation}
so that $\ell_\theta(A)$ is defined by the spectral calculus on all of
$\sigma(A)$.  The value at zero matters because $0\in\sigma(A)$ for every
compact operator on an infinite-dimensional space.
\end{assumption}

\begin{theorem}[Logarithmic regularization-bias rate]\label{thm:log-rate}
Under Assumption~\ref{ass:log-source}, for every $0<\alpha\le\norm A$ and
every $z\in K_0$,
\begin{equation}\label{eq:log-bias}
 0\le W(z)-W_\alpha(z)
 \le 3\,c_\theta\,C_q^{2}\,
 \Bigl(1+\log\frac{\norm A}{\alpha}\Bigr)^{-2\theta},
 \qquad
 c_\theta:=\max\bigl\{1,\;e\,(2\theta/e)^{2\theta}\bigr\}.
\end{equation}
Consequently, on the event of \eqref{eq:data-bound}, every parameter choice
$\alpha=e_{M,N}^{\,\kappa}$ with a fixed $\kappa\in(0,1/2)$ and
test-vector rate
\[
 \zeta_M=O\Bigl(e_{M,N}^{\,\kappa}
 \bigl(\log(1/e_{M,N})\bigr)^{-2\theta}\Bigr)
\]
gives
\begin{equation}\label{eq:log-total}
 \sup_{z\in K_0}
 \bigl|\widehat W_{\alpha,M,N}(z)-W(z)\bigr|
 =O\!\left(\bigl(\log(1/e_{M,N})\bigr)^{-2\theta}\right).
\end{equation}
Here the operator-error contribution decays faster than every fixed
negative power of $\log(1/e_{M,N})$, while the test-vector contribution
is at most of the order of the bias.  The exponent $\kappa$ affects only
the constants.  Under the stronger hypothesis
$\zeta_M=O(e_{M,N}^{\,\kappa+\delta})$ with a fixed $\delta>0$ the
test-vector contribution too decays faster than every fixed negative
power of $\log(1/e_{M,N})$.
\end{theorem}

\begin{proof}
As in the proof of Theorem~\ref{thm:source-rate}, with $E_A$ the spectral
measure of $A$ and $\mu_z=\ip{E_A(\cdot)q_z}{q_z}$, the source representation
\eqref{eq:log-source} gives
$|\ip{\varphi_z}{\psi_j}|^2
=\gamma_j\,\ell_\theta(\gamma_j)^2\,|\ip{q_z}{\psi_j}|^2$ and hence
\[
 W(z)-W_\alpha(z)
 =\int_{(0,\norm A]}
 \frac{\alpha(2t+\alpha)}{(t+\alpha)^2}\,
 \ell_\theta(t)^2\,\dd\mu_z(t),
\]
using $\bigl(\tfrac1t-\tfrac{t}{(t+\alpha)^2}\bigr)t
=\tfrac{\alpha(2t+\alpha)}{(t+\alpha)^2}$.  We bound the integrand
$h(t)=\alpha(2t+\alpha)(t+\alpha)^{-2}\ell_\theta(t)^2$ on the two ranges
$t\le\alpha$ and $t>\alpha$.

The rational factor $\alpha(2t+\alpha)(t+\alpha)^{-2}$ is at most
$3\min(1,\alpha/t)$, so, $\ell_\theta$ being nondecreasing,
$h\le3\ell_\theta(\alpha)^2$ for $t\le\alpha$, while for $t>\alpha$ the
substitution $t=\norm A e^{-\ell}$, $L=\log(\norm A/\alpha)$, gives
$h\le3(\alpha/\norm A)e^{\ell}(1+\ell)^{-2\theta}$.  The latter function
has logarithmic derivative $1-2\theta/(1+\ell)$, so it decreases on
$[0,2\theta-1]$ and increases thereafter and its maximum over $[0,L]$
is attained at an endpoint, namely
$\max\{e^{-L},(1+L)^{-2\theta}\}$.  The constant $c_\theta$ controls
the first alternative: $L\mapsto
e^{-L}(1+L)^{2\theta}$ equals $1$ at $L=0$ and, for $\theta>1/2$, peaks
at $L=2\theta-1$ with value $e(2\theta/e)^{2\theta}$, so
$e^{-L}\le c_\theta(1+L)^{-2\theta}$ throughout.  Both ranges therefore
give $h\le3c_\theta\ell_\theta(\alpha)^2$, and integrating against
$\mu_z$, of mass $\norm{q_z}^2\le C_q^2$, yields \eqref{eq:log-bias}.  For \eqref{eq:log-total}, add \eqref{eq:data-bound} with $\alpha=e_{M,N}^{\kappa}$.  The
perturbation term is
$O\bigl(e_{M,N}^{1-2\kappa}\log\log(1/e_{M,N})\bigr)$ by
\eqref{eq:rho-modulus}, which vanishes faster
than any power of $1/\log(1/e_{M,N})$, while the bias term is
$\bigl(1+\kappa\log(1/e_{M,N})+\log\norm A\bigr)^{-2\theta}
\asymp\bigl(\log(1/e_{M,N})\bigr)^{-2\theta}$.
\end{proof}

\begin{proposition}[The disc satisfies the logarithmic condition for every exponent]\label{prop:disc-log}
For the impedance disc and every compact $K_0\Subset D$,
Assumption~\ref{ass:log-source} holds for every $\theta>0$, with a constant
$C_q$ depending on $\theta$, $K_0$, $k$, $a$, and $\lambda$.
\end{proposition}

\begin{proof}
By the spectral characterization, \eqref{eq:log-source} holds at $z$ if and
only if
\begin{equation}\label{eq:log-picard}
 \sum_{m\in\mathbb Z}
 \frac{|\ip{\varphi_z}{\psi_m}|^2}{\gamma_m\,\ell_\theta(\gamma_m)^2}
 =\sum_{m\in\mathbb Z}
 \frac{2\pi J_m(k|z|)^2}{\gamma_m}
 \Bigl(1+\log\frac{\norm A}{\gamma_m}\Bigr)^{2\theta}
 <\infty,
\end{equation}
and $C_q^2$ may be taken as the supremum of this sum over $z\in K_0$.  Let
$r_{\max}=\sup_{z\in K_0}|z|<a$.  The power series of $J_m$ gives, for
every $x\ge0$,
\[
 |J_m(x)|\le\frac{(x/2)^m}{m!}
 \exp\Bigl(\frac{x^2}{4(m+1)}\Bigr),
\]
so by \eqref{eq:disc-coefficients}, Stirling's bound
$m!\ge(m/e)^m\sqrt{2\pi m}$, and \eqref{eq:gamma-twosided} there are $C$
and $m_0$, uniform for $|z|\le r_{\max}$ including $z=0$, such that for
$m\ge m_0$
\[
 \frac{2\pi J_m(k|z|)^2}{\gamma_m}
 \le\frac{C}{m}\Bigl(\frac{r_{\max}}{a}\Bigr)^{2m},
 \qquad
 1+\log\frac{\norm A}{\gamma_m}\le C\,m\log m .
\]
The terms of \eqref{eq:log-picard} are therefore bounded by
$C^{1+2\theta}(r_{\max}/a)^{2m}\,m^{2\theta-1}(\log m)^{2\theta}$ for
$m\ge m_0$, so the geometric factor dominates the polynomial--logarithmic one, and
the series converges, uniformly in $z\in K_0$.  The finitely many remaining
terms are bounded because every $\gamma_m$ is strictly positive.
\end{proof}

\begin{remark}[Observed bias decay for the disc]\label{rem:disc-bias}
Proposition~\ref{prop:disc-log} predicts a bias decaying faster than any
fixed power of $1/\log(\norm A/\alpha)$.  Numerically ($a=1.3$,
$k=\lambda=1$, $r=1.0$) the product
$(W-W_\alpha)\cdot(1+\log(\norm A/\alpha))^{2}$ falls from $3.0$ at
$\alpha=10^{-2}$ to $1.2$ at $\alpha=10^{-12}$, so the bias tracks the
$\theta=1$ rate while slowly outperforming it, as expected when the condition
holds for every exponent.
\end{remark}

\begin{remark}[Relation to variational source conditions]\label{rem:vsc}
For nonlinear inverse problems, spectral source conditions have largely
been superseded by variational source conditions, which avoid the operator
calculus altogether and yield rates for Tikhonov regularization directly
\cite{HKPS07}.  Conditions of logarithmic type have been verified for the
acoustic inverse medium problem under Sobolev smoothness of the contrast
\cite{HW15}, and for linear problems the variational and spectral
formulations are connected by converse characterizations \cite{HW17}.
Those results concern the nonlinear map from the scatterer to the data.
The object regularized here is different.  The Picard series
\eqref{eq:intro-picard} is itself a spectral quantity of $A=\Fh$, and its
bias $W(z)-W_\alpha(z)$ is determined by how $\varphi_z$ distributes over
the eigensystem of $A$.  We therefore formulate the source condition
spectrally.  Propositions~\ref{prop:hoelder-fails} and~\ref{prop:disc-log}
show that the disc satisfies no positive-order H\"older source condition
away from its center, whereas it satisfies every logarithmic source
condition considered here.
\end{remark}

\section{Sharp transition asymptotics for the impedance disc}\label{sec:disc-sharp}
Let
\[
 C_k=\sqrt{\frac{2\pi}{k}},
 \qquad
 L_\alpha=\log\frac{C_k}{\alpha},
\]
and let $\LambertW$ denote the principal branch of the Lambert
$W$-function \cite{CGHJK96}.

\begin{lemma}[Sharpened eigenvalue asymptotics and adjacent ratio]
\label{lem:gamma-ratio}
For the impedance disc,
\begin{equation}\label{eq:gamma-sharp}
 \gamma_m
 =C_k\left(\frac{eka}{2m}\right)^{2m}\bigl(1+O(m^{-1})\bigr),
 \qquad
 \log\frac{C_k}{\gamma_m}
 =2m\log\frac{2m}{eka}+O(m^{-1}),
\end{equation}
and
\begin{equation}\label{eq:gamma-ratio}
 \frac{\gamma_{m+1}}{\gamma_m}
 =\frac{(ka)^{2}}{4m^{2}}\bigl(1+O(m^{-1})\bigr);
\end{equation}
in particular there are constants $0<c\le C<\infty$ and an index $m_1$ with
$c\,m^{-2}\le\gamma_{m+1}/\gamma_m\le C\,m^{-2}$ and
$\gamma_{m+1}<\gamma_m$ for all $m\ge m_1$.
\end{lemma}

The scale in \eqref{eq:gamma-sharp}, constant included, is already
implicit in the literature.  Ma and Hu~\cite[(2.11)--(2.12)]{MH22} compute the
exact eigensystem of the impedance disc and give
\[
 |\kappa_m|=\sqrt{\frac{2}{k\pi}}\,
 \frac{\pi^{2}(ka)^{2m}}{2^{2m-1}m!\,(m-1)!}\bigl(1+O(m^{-1})\bigr).
\]
Since
$m!(m-1)!=2\pi(m/e)^{2m}\bigl(1+O(m^{-1})\bigr)$, this is
$C_k(eka/2m)^{2m}\bigl(1+O(m^{-1})\bigr)$, the same constant.  Their
estimate is for the modulus of an eigenvalue of $F$ itself, and it is
all they require.  The spectrum of $\Fh$ is not computed there but
replaced by $|\kappa_m|$ through
$\tfrac1{\sqrt2}(|\Real\kappa_m|+|\Imag\kappa_m|)\le|\kappa_m|
\le|\Real\kappa_m|+|\Imag\kappa_m|$, which settles positivity of the
indicator but leaves the prefactor undetermined to within a factor
$\sqrt2$.  Equation~\eqref{eq:gamma-sharp} also determines the phase.  The leading term
of $\widetilde\kappa_m=e^{-i\pi/4}\kappa_m$ is positive real, so the
symmetrization $|\Real\widetilde\kappa_m|+\Imag\widetilde\kappa_m$
retains the constant rather than a bounded multiple of it, as required
for the capacity law in \cref{sec:general-capacity}.

\begin{proof}
For fixed $x>0$ the large-order expansions are asymptotic series in $1/m$
\cite[Section~10.19]{DLMF},
\[
 J_m(x)=\frac{(ex/2m)^{m}}{\sqrt{2\pi m}}\bigl(1+O(m^{-1})\bigr),
 \qquad
 Y_m(x)=-\sqrt{\frac{2}{\pi m}}\Bigl(\frac{ex}{2m}\Bigr)^{-m}
 \bigl(1+O(m^{-1})\bigr),
\]
and the recurrences give the same relative accuracy for the derivative
combinations entering $a_m$.  The impedance terms contribute relative errors
$i\lambda x/m+O(m^{-2})$.  Every step of the proof of
Proposition~\ref{prop:hoelder-fails} therefore carries a multiplicative
factor $1+O(m^{-1})$, so
$a_m=\tfrac i2(eka/2m)^{2m}(1+\delta_m)$ with a complex
$\delta_m=O(m^{-1})$, and, with $q_m:=(eka/2m)^{2m}$,
\[
 \gamma_m=|\Real\widetilde\kappa_m|+\Imag\widetilde\kappa_m
 =C_kq_m(1+\eta_m),
 \qquad \eta_m\in\R,\quad \eta_m=O(m^{-1}),
\]
which is \eqref{eq:gamma-sharp}.  Because the errors are multiplicative,
adjacent terms may be divided without loss, whence
\[
 \frac{q_{m+1}}{q_m}
 =\Bigl(\frac{eka}{2(m+1)}\Bigr)^{2}\Bigl(1+\frac1m\Bigr)^{-2m}
 =\frac{(ka)^{2}}{4m^{2}}\bigl(1+O(m^{-1})\bigr),
\]
using $(1+1/m)^{-2m}=e^{-2}(1+O(m^{-1}))$, and therefore
\[
 \frac{\gamma_{m+1}}{\gamma_m}
 =\frac{q_{m+1}}{q_m}\,\frac{1+\eta_{m+1}}{1+\eta_m}
 =\frac{(ka)^{2}}{4m^{2}}\bigl(1+O(m^{-1})\bigr),
\]
which is \eqref{eq:gamma-ratio}.  The right-hand side is eventually below
one, so the sequence decreases strictly from some $m_1$ on.
\end{proof}

By Lemma~\ref{lem:gamma-ratio} there is $m_0$ such that
$(\gamma_m)_{m\ge m_0}$ is strictly decreasing.  For sufficiently small
$\alpha$, define
\begin{equation}\label{eq:transition-index}
 m_\alpha=\max\{m\ge m_0:\gamma_m\ge\alpha\}.
\end{equation}

\begin{theorem}[Lambert-$W$ transition and two-sided disc asymptotics]
\label{thm:disc-sharp}
For the impedance disc,
\begin{equation}\label{eq:malpha-asymptotic}
 m_\alpha=
 \frac{L_\alpha}
 {2\LambertW\!\bigl(L_\alpha/(eka)\bigr)}
 +O(1),
 \qquad \alpha\downarrow0.
\end{equation}
For every fixed $r>0$, the following statements hold as $\alpha\downarrow0$.
\begin{align}
 \log\bigl(W(r)-W_\alpha(r)\bigr)
 &=2m_\alpha\log(r/a)+O(\log m_\alpha),
 &&0<r<a, \label{eq:sharp-interior}\\
 W_\alpha(a)&\asymp\log m_\alpha,
 &&r=a, \label{eq:sharp-boundary}\\
 \log W_\alpha(r)
 &=2m_\alpha\log(r/a)+O(\log m_\alpha),
 &&r>a. \label{eq:sharp-exterior}
\end{align}
The constants implicit in the $O(\cdot)$ and $\asymp$ notation may depend on
$r,a,k,\lambda$, but not on $\alpha$.
\end{theorem}

\begin{proof}
Throughout, \eqref{eq:gamma-sharp} and \eqref{eq:gamma-ratio} of
Lemma~\ref{lem:gamma-ratio} are used.  In particular
$m_\alpha\to\infty$ as $\alpha\downarrow0$.

\emph{The index.}  Definition \eqref{eq:transition-index} gives
$\gamma_{m_\alpha}\ge\alpha>\gamma_{m_\alpha+1}$, and
\eqref{eq:gamma-ratio} gives
$\gamma_{m_\alpha+1}\ge c\,m_\alpha^{-2}\gamma_{m_\alpha}$.  Hence
\begin{equation}\label{eq:malpha-bracket}
 \alpha\le\gamma_{m_\alpha}\le c^{-1}m_\alpha^{2}\,\alpha,
 \qquad\text{so}\qquad
 \log\frac{C_k}{\gamma_{m_\alpha}}=L_\alpha+O(\log m_\alpha).
\end{equation}
Because consecutive eigenvalues differ by factors $\asymp m^{-2}$,
one cannot strengthen this to $\gamma_{m_\alpha}\asymp\alpha$ with
$\alpha$-independent constants.
Combining \eqref{eq:malpha-bracket} with \eqref{eq:gamma-sharp} at
$m=m_\alpha$,
\[
 2m_\alpha\log\frac{2m_\alpha}{eka}
 =L_\alpha+E_\alpha,
 \qquad E_\alpha=O(\log m_\alpha),
\]
which in particular forces $L_\alpha\asymp m_\alpha\log m_\alpha$.  For
$L'>0$ let $m(L')$ solve $2m\log(2m/(eka))=L'$.  With $x=2m/(eka)$ and
$t=L'/(eka)$ the equation reads $x\log x=t$, whose solution is
$x=t/\LambertW(t)$, so that
\[
 m(L')=\frac{L'}{2\LambertW\bigl(L'/(eka)\bigr)},
 \qquad
 \frac{dm}{dL'}
 =\frac{1}{2\bigl(1+\LambertW(L'/(eka))\bigr)},
\]
using $\frac{d}{dt}\bigl(t/\LambertW(t)\bigr)=1/(1+\LambertW(t))$.  Since
$\LambertW(t)\ge\tfrac12\log t$ for large $t$, the derivative is
$O(1/\log L')$ on the interval with endpoints $L_\alpha$ and
$L_\alpha+E_\alpha$, and the mean value
theorem gives
\[
 m_\alpha=m(L_\alpha+E_\alpha)
 =m(L_\alpha)+O\Bigl(\frac{\log m_\alpha}{\log L_\alpha}\Bigr)
 =\frac{L_\alpha}{2\LambertW\bigl(L_\alpha/(eka)\bigr)}+O(1),
\]
by $L_\alpha\asymp m_\alpha\log m_\alpha$, which proves
\eqref{eq:malpha-asymptotic}.

\emph{Filter bounds.}  For fixed $r>0$, the Bessel and eigenvalue
asymptotics give constants $c_r,C_r>0$ such that, for all sufficiently
large $m$,
\begin{equation}\label{eq:disc-ratio}
 \frac{c_r}{m}\left(\frac ra\right)^{2m}
 \le
 R_m(r):=\frac{2\pi J_m(kr)^2}{\gamma_m}
 \le
 \frac{C_r}{m}\left(\frac ra\right)^{2m}.
\end{equation}
Write $q_{\alpha,m}=\gamma_m^2/(\gamma_m+\alpha)^2$, so the mode-$m$
contributions to $W_\alpha(r)$ and to $W(r)-W_\alpha(r)$ are
$R_m(r)\,q_{\alpha,m}$ and $R_m(r)\,(1-q_{\alpha,m})$.  From
$1-q_{\alpha,m}=\alpha(2\gamma_m+\alpha)/(\gamma_m+\alpha)^2$,
\begin{equation}\label{eq:filter-bracket}
 \frac{\alpha}{2\gamma_m}\le1-q_{\alpha,m}\le\frac{3\alpha}{\gamma_m}
 \quad(\gamma_m\ge\alpha),
 \qquad
 \frac14\Bigl(\frac{\gamma_m}{\alpha}\Bigr)^{2}
 \le q_{\alpha,m}\le\Bigl(\frac{\gamma_m}{\alpha}\Bigr)^{2}
 \quad(\gamma_m\le\alpha).
\end{equation}
Throughout, negative modes are symmetric ($\gamma_{-m}=\gamma_m$,
$J_{-m}^2=J_m^2$) and the finitely many modes below $m_0$ contribute at
most a constant multiple of $\alpha$ to the bias $W-W_\alpha$ (which is
$o\bigl((r/a)^{2m_\alpha}\bigr)$ because $\alpha\le\gamma_{m_\alpha}$
decays super-exponentially while $(r/a)^{2m}$ is geometric) and $O(1)$ to
$W_\alpha$ itself, absorbed in the boundary and exterior statements since
$W_\alpha(a)\to\infty$ and $\log W_\alpha$ carries an $O(\log m_\alpha)$
remainder.  Both effects only change constants.

\emph{Interior points $0<r<a$.}  For $m_0\le m\le m_\alpha$,
\eqref{eq:filter-bracket} bounds the terms by $3\alpha R_m/\gamma_m$.  By
\eqref{eq:disc-ratio} and \eqref{eq:gamma-ratio} the ratio of consecutive
values of $R_m/\gamma_m$ is at least $c'(r/a)^2m^{2}$ for a constant
$c'>0$, hence at least $2$ for large $m$, so the block is at most
twice its last term, namely
$\sum_{m\le m_\alpha}3\alpha R_m/\gamma_m
\le6\,(\alpha/\gamma_{m_\alpha})R_{m_\alpha}\le6R_{m_\alpha}$, using
$\alpha\le\gamma_{m_\alpha}$.  For $m>m_\alpha$, using $1-q_{\alpha,m}\le1$
and the geometric decay of \eqref{eq:disc-ratio}, the tail is at most
$C\,m_\alpha^{-1}(r/a)^{2m_\alpha}$.  For the lower bound, keep the single
mode $m_\alpha+1$, where $\gamma_m<\alpha$ gives
$q_{\alpha,m}=\gamma_m^2/(\gamma_m+\alpha)^2<\tfrac14$, so
$1-q_{\alpha,m}\ge\tfrac34$ and
$W(r)-W_\alpha(r)\ge\tfrac34R_{m_\alpha+1}(r)\ge
c\,m_\alpha^{-1}(r/a)^{2m_\alpha+2}$.  Combining,
\[
 c\,m_\alpha^{-1}\Bigl(\frac ra\Bigr)^{2m_\alpha+2}
 \le W(r)-W_\alpha(r)
 \le C\,m_\alpha^{-1}\Bigl(\frac ra\Bigr)^{2m_\alpha},
\]
and taking logarithms gives \eqref{eq:sharp-interior}.

\emph{Exterior points $r>a$.}  For $m\le m_\alpha$ use $q_{\alpha,m}\le1$.
Since $R_m$ grows geometrically, the block is at most
$C\,m_\alpha^{-1}(r/a)^{2m_\alpha}$.  For $m>m_\alpha$ use
$q_{\alpha,m}\le(\gamma_m/\alpha)^2\le1$ at $m=m_\alpha+1$ together with
the ratio bound, by which consecutive terms of $R_m(\gamma_m/\alpha)^2$ shrink by
factors at most $C'(r/a)^2m^{-4}\le\tfrac12$ for large $m$, with $C'$
determined by \eqref{eq:disc-ratio} and \eqref{eq:gamma-ratio}, so the tail is at
most $2R_{m_\alpha+1}\le C\,m_\alpha^{-1}(r/a)^{2m_\alpha+2}$.  For the
lower bound, the single mode $m_\alpha$ has
$q_{\alpha,m_\alpha}\ge\tfrac14$ by \eqref{eq:filter-bracket}, so
$W_\alpha(r)\ge\tfrac14R_{m_\alpha}(r)\ge
c\,m_\alpha^{-1}(r/a)^{2m_\alpha}$.  Taking logarithms gives
\eqref{eq:sharp-exterior}.

\emph{Boundary $r=a$.}  Here $R_m(a)\asymp m^{-1}$ by
\eqref{eq:disc-ratio}.  The lower bound follows from the modes
$m_0\le m\le m_\alpha$, where $q_{\alpha,m}\ge\tfrac14$, giving
$W_\alpha(a)\ge\tfrac14\sum_{m_0}^{m_\alpha}c\,m^{-1}\ge c'\log m_\alpha$.
For the upper bound,
$\sum_{m\le m_\alpha}C\,m^{-1}\le C\log m_\alpha+C$, while for
$m>m_\alpha$ the terms $C\,m^{-1}(\gamma_m/\alpha)^2$ start below
$C\,m_\alpha^{-1}$ and decay faster than geometrically, contributing $O(1)$.
Hence $W_\alpha(a)\asymp\log m_\alpha$.
\end{proof}

\begin{corollary}[Boundary growth law with explicit constant]\label{cor:boundary-constant}
For the impedance disc there is a constant $B=B(a,k,\lambda)$ such that,
for all sufficiently small $\alpha$,
\begin{equation}\label{eq:boundary-constant}
 \Bigl|\,W_\alpha(a)-\frac{2}{C_k}\log m_\alpha\Bigr|\le B .
\end{equation}
Consequently, by \eqref{eq:malpha-asymptotic},
\begin{equation}\label{eq:boundary-loglog}
 W_\alpha(a)=\sqrt{\frac{2k}{\pi}}\,
 \bigl(\log\log(1/\alpha)\bigr)\bigl(1+o(1)\bigr),
 \qquad\alpha\downarrow0,
\end{equation}
the radius and the impedance entering only through the bounded defect and
the $o(1)$.
\end{corollary}

\begin{proof}
Combining \eqref{eq:gamma-sharp} with the $J_m$ asymptotics of
Lemma~\ref{lem:gamma-ratio} at $x=ka$,
\begin{equation}\label{eq:Rm-boundary}
 R_m(a)=\frac{2\pi J_m(ka)^2}{\gamma_m}
 =\frac{1}{C_k m}\bigl(1+O(m^{-1})\bigr),
\end{equation}
the prefactors $(eka/2m)^{2m}$ cancelling exactly.  Fix $m_1\ge m_0$ so
large that the relative error in \eqref{eq:Rm-boundary} is at most
$\tfrac12$ and that $\gamma_{m+1}/\gamma_m\le\tfrac12$ for $m\ge m_1$
(Lemma~\ref{lem:gamma-ratio}).  With $T_m:=R_m(a)\,q_{\alpha,m}$ and
$q_{\alpha,m}=\gamma_m^2/(\gamma_m+\alpha)^2$ as in the proof of
Theorem~\ref{thm:disc-sharp}, symmetry gives
\[
 W_\alpha(a)=\sum_{|m|\le m_1}T_m
 +2\!\!\sum_{m_1<m\le m_\alpha}\!\!T_m
 +2\!\!\sum_{m>m_\alpha}\!\!T_m .
\]
The first block is at most $\sum_{|m|\le m_1}R_m(a)$, a constant.  In the
middle block,
\[
 2\!\!\sum_{m_1<m\le m_\alpha}\!\!\frac{1}{C_km}\bigl(1+O(m^{-1})\bigr)
 =\frac{2}{C_k}\log m_\alpha+O(1),
\]
so it remains to bound the deficit
$2\sum_{m_1<m\le m_\alpha}R_m(a)(1-q_{\alpha,m})$ and the tail.  For
$0\le j<m_\alpha-m_1$ the ratio bound gives
$\gamma_{m_\alpha-j}\ge2^{j}\gamma_{m_\alpha}\ge2^{j}\alpha$, whence
$1-q_{\alpha,m_\alpha-j}\le3\alpha/\gamma_{m_\alpha-j}\le3\cdot2^{-j}$ by
\eqref{eq:filter-bracket}, and
\[
 \sum_{m_1<m\le m_\alpha}\frac{1-q_{\alpha,m}}{m}
 \le\frac{6}{m_\alpha}\sum_{0\le j\le m_\alpha/2}2^{-j}
 +3\cdot2^{-m_\alpha/2}\!\!\sum_{m_1<m\le m_\alpha}\frac1m
 =O(m_\alpha^{-1}).
\]
In the tail, $\gamma_{m_\alpha+1+j}\le2^{-j}\gamma_{m_\alpha+1}
<2^{-j}\alpha$ gives $q_{\alpha,m_\alpha+1+j}\le4^{-j}$ by
\eqref{eq:filter-bracket}, so
$\sum_{m>m_\alpha}T_m\le\frac{2}{C_km_\alpha}\sum_{j\ge0}4^{-j}
=O(m_\alpha^{-1})$.  Collecting the estimates proves
\eqref{eq:boundary-constant}.  For \eqref{eq:boundary-loglog}, write
$P_\alpha=L_\alpha/\bigl(2\LambertW(L_\alpha/(eka))\bigr)$.  Then
\eqref{eq:malpha-asymptotic} gives
$\log m_\alpha=\log P_\alpha+O(P_\alpha^{-1})$, and since
$1\le\LambertW(t)\le\log t$ for $t\ge e$,
$\log P_\alpha=(1+o(1))\log L_\alpha=(1+o(1))\log\log(1/\alpha)$.
Finally $2/C_k=\sqrt{2k/\pi}$.
\end{proof}

\section{Set reconstruction}\label{sec:sets}
\subsection{Measure consistency of sublevel sets}
The obstacle is approximated by a \emph{sublevel} set
rather than by a superlevel set,
\begin{equation}
 D_{\alpha}^{\tau}=\{z\in K:W_\alpha(z)\le\tau\}.
\end{equation}

\begin{theorem}[Diagonal set consistency]\label{thm:set-consistency}
Let $K\subset\R^2$ be a bounded sampling window containing
$\overline D$.  There exist sequences
$\alpha_n\downarrow0$ and $\tau_n\uparrow\infty$ such that
\begin{equation}\label{eq:set-measure}
 \left|
 \{z\in K:W_{\alpha_n}(z)\le\tau_n\}\triangle D
 \right|\longrightarrow0.
\end{equation}
If empirical fields satisfy
$\sup_K|\widehat W_n-W_{\alpha_n}|\le\epsilon_n$, then
\begin{equation}\label{eq:set-sandwich}
 \{W_{\alpha_n}\le\tau_n-\epsilon_n\}
 \subseteq
 \{\widehat W_n\le\tau_n\}
 \subseteq
 \{W_{\alpha_n}\le\tau_n+\epsilon_n\}.
\end{equation}
Thus the same measure convergence follows whenever the measure of the ideal
transition band
$\{|W_{\alpha_n}-\tau_n|\le\epsilon_n\}$ tends to zero.
\end{theorem}

\begin{proof}
The argument is the standard diagonal one for monotone level sets.
With $D^\tau=\{z\in K:W(z)\le\tau\}$ one has $D^\tau\subset D$
increasing to $D$ by \eqref{eq:intro-picard}, so
$|D\setminus D^\tau|\to0$ by continuity from below, while for fixed
$\tau$ the sets $D_\alpha^\tau$ decrease to $D^\tau$ as
$\alpha\downarrow0$ because $W_\alpha\uparrow W$.  A diagonal choice of
$(\alpha_n,\tau_n)$ gives \eqref{eq:set-measure}, and
\eqref{eq:set-sandwich} is immediate from the uniform error bound.
\end{proof}

\subsection{Conditional Hausdorff boundary rate}
A quantitative boundary-location theorem requires a margin or transversality
assumption.  Uniform field convergence alone is insufficient.  For a chosen
pair $(\alpha,\tau)$, assume $\Gamma_{\alpha,\tau}$ is a connected
component of $\{W_\alpha=\tau\}$ forming a
simple closed $C^1$ curve with a tubular normal coordinate $x+s\nu(x)$, $|s|<r_0$,
and that there are $c>0$ and $p\ge1$ such that, after orienting $\nu$ from
the sublevel side to the superlevel side,
\begin{equation}\label{eq:transversality}
 W_\alpha(x+s\nu(x))-\tau
 \begin{cases}
 \le-c|s|^p,&-r_0<s<0,\\
 \ge c|s|^p,&0<s<r_0.
 \end{cases}
\end{equation}

\begin{proposition}[Hausdorff stability of the reconstructed boundary]\label{prop:hausdorff}
Assume \eqref{eq:transversality} and let $\widehat W$ be continuous with
$\norm{\widehat W-W_\alpha}_{L^\infty}\le\epsilon<cr_0^p$.  With
$\widehat\Gamma^{\rm tube}=\{\widehat W=\tau\}\cap
U_{r_0}(\Gamma_{\alpha,\tau})$ the empirical level set localized to the
tubular neighborhood, every normal fiber through the ideal curve contains a
point of $\widehat\Gamma^{\rm tube}$ and
\begin{equation}\label{eq:hausdorff}
 d_H(\widehat\Gamma^{\rm tube},\Gamma_{\alpha,\tau})
 \le(\epsilon/c)^{1/p};
\end{equation}
no conclusion is asserted for level-set components outside the tubular
neighborhood.  If in addition
$d_H(\Gamma_{\alpha_n,\tau_n},\partial D)\le h_n$, then
$d_H(\widehat\Gamma_n^{\rm tube},\partial D)\le
h_n+(\epsilon_n/c_n)^{1/p_n}$.
\end{proposition}

\begin{proof}
Both inclusions follow from the uniform perturbation against the
transverse margin.  At distance $r=(\epsilon/c)^{1/p}$ on either side of
an ideal boundary point, \eqref{eq:transversality} forces opposite signs
of $\widehat W-\tau$, so each normal fiber carries an empirical zero,
while a zero farther than $r$ inside the tube would contradict the
margin.  The final bound is the triangle inequality for $d_H$.
\end{proof}

\subsection{Quantitative rates for the disc}\label{sec:disc-recon}
Throughout this subsection $K\supset\overline D$ is a bounded
sampling window with $R_{\max}:=\sup_{z\in K}|z|>a$, and $W_\alpha(z)$, which
depends only on $r=|z|$, is written $W_\alpha(r)$.  The basic estimate,
obtained exactly as \eqref{eq:Rm-boundary} but keeping the radial factor,
is
\begin{equation}\label{eq:Rm-uniform}
 R_m(r)=\frac{2\pi J_m(kr)^2}{\gamma_m}
 =\frac{1}{C_k m}\Bigl(\frac ra\Bigr)^{2m}\bigl(1+O(m^{-1})\bigr),
\end{equation}
with the $O(m^{-1})$ uniform for $r$ in compact subsets of $(0,\infty)$.

\begin{lemma}[Boundary blow-up of the Picard series]\label{lem:picard-blowup}
For every $r_1\in(0,a)$ there is a constant $C=C(r_1,a,k,\lambda)$ with
\begin{equation}\label{eq:picard-blowup}
 \Bigl|\,W(r)-\frac{2}{C_k}\log\frac{1}{1-(r/a)^2}\Bigr|\le C,
 \qquad r_1\le r<a,
\end{equation}
and $W$ is bounded on $\{|z|\le r_1\}$.
\end{lemma}

\begin{proof}
Fix $m_1=m_1(r_1)$ so large that \eqref{eq:Rm-uniform} holds on $[r_1,a]$
for $m>m_1$.  Since $W(r)=\sum_{m\in\mathbb Z}R_m(r)$,
\[
 \sum_{|m|>m_1}R_m(r)
 =\frac{2}{C_k}\sum_{m>m_1}\frac1m\Bigl(\frac ra\Bigr)^{2m}
 +O\Bigl(\sum_{m>m_1}\frac1{m^2}\Bigr)
 =\frac{2}{C_k}\log\frac{1}{1-(r/a)^2}+O(1),
\]
uniformly for $r\in[r_1,a)$, by the Mercator series
$\sum_{m\ge1}x^m/m=\log\frac1{1-x}$ at $x=(r/a)^2$, the omitted initial
segment being bounded by $\sum_{m\le m_1}1/m$.  The modes $|m|\le m_1$
contribute at most $\sum_{|m|\le m_1}2\pi/\gamma_m<\infty$, since
$J_m(kr)^2\le1$.  This proves \eqref{eq:picard-blowup}.  For $|z|\le r_1$,
every $J_m(kr)$ with $m>kr_1$ is increasing on $[0,r_1]$, its first
critical point $j'_{m,1}$ exceeding $m$, so the tail of the series for $W$
is dominated by its value at $r_1$ and $W$ is bounded on $\{|z|\le r_1\}$.
\end{proof}

The sublevel sets $\{W\le\tau\}$ therefore exhaust the disc at speed
$1-(r/a)^2\asymp e^{-C_k\tau/2}$, and by
Corollary~\ref{cor:boundary-constant} the threshold scale is set by the
boundary value $W_\alpha(a)$.

\begin{lemma}[Near-boundary behavior of the regularized field]
There are constants $0<c\le C$ and $\alpha_0>0$, depending only on
$a,k,\lambda$, such that for $\alpha\le\alpha_0$, with $u=(r-a)/a$,
\begin{align}
 W_\alpha\bigl(a(1+u)\bigr)
 &\le \frac{2}{C_k}\log\frac1u
 +C\Bigl(\frac{e^{2m_\alpha u}}{m_\alpha u}+1\Bigr),
 &&\tfrac{2}{m_\alpha}\le u\le\tfrac12, \label{eq:Wa-upper}\\
 W_\alpha\bigl(a(1+u)\bigr)
 &\ge \frac{2}{C_k}\log\frac1u
 +c\,\frac{e^{2m_\alpha u}}{m_\alpha u}-C,
 &&\tfrac{2}{m_\alpha}\le u\le m_\alpha^{-1/2}, \label{eq:Wa-lower}\\
 W_\alpha(r)&\le \frac{2}{C_k}\log m_\alpha+C,
 &&0\le r\le a\bigl(1+\tfrac{2}{m_\alpha}\bigr). \label{eq:Wa-flat}
\end{align}
\end{lemma}

\begin{proof}
Fix $m_1\ge\max\{m_0,2ka,(2ka)^2\}$ so large that \eqref{eq:Rm-uniform}
holds on $[a,\tfrac32a]$ with relative error $\le\tfrac12$ and
$\gamma_{m+1}/\gamma_m\le\tfrac12$ for $m\ge m_1$.  Write
$T_m=R_m\,q_{\alpha,m}$ for the mode-$m$ term of $W_\alpha$.  The modes
$|m|\le m_1$ contribute between $0$ and a constant throughout, and by the
symmetry $\gamma_{-m}=\gamma_m$ the remaining sum is twice the sum over
$m>m_1$.

\emph{Upper bounds.}  For $r=a(1+u)$ with $0\le u\le\tfrac12$,
$T_m\le R_m\le\frac{1}{C_km}e^{2mu}\bigl(1+O(m^{-1})\bigr)$.  On
$m\le1/u$, $e^{2mu}-1\le2mu\,e^{2}$, so
$\sum_{m_1<m\le1/u}e^{2mu}/m\le\log\frac1u+2e^{2}+1$.  On
$1/u<m\le m_\alpha$, write $m=m_\alpha-j$.  The block with
$j\le m_\alpha/2$ is at most
$e^{2m_\alpha u}\tfrac{2}{m_\alpha}\sum_{j\ge0}e^{-2ju}
\le2e^{2m_\alpha u}/(m_\alpha u)$, using $1-e^{-2u}\ge u$, valid on
$0\le u\le\tfrac12$ because $u\mapsto1-e^{-2u}-u$ is concave and
vanishes at $u=0$ while being positive at $u=\tfrac12$; the block with
$j>m_\alpha/2$, that is $1/u<m\le m_\alpha/2$, obeys
\[
 \sum_{1/u<m\le m_\alpha/2}\frac{e^{2mu}}m
 \le u\!\!\sum_{m\le m_\alpha/2}\!\!e^{2mu}
 \le u\,\frac{e^{m_\alpha u}}{1-e^{-2u}}
 \le e^{m_\alpha u}
 \le\frac{e^{2m_\alpha u}}{e\,m_\alpha u},
\]
using $1/m\le u$ there, $1-e^{-2u}\ge u$, and $se^{-s}\le e^{-1}$ at
$s=m_\alpha u$.  The correction sums $\sum m^{-2}e^{2mu}$ are bounded by
splitting the same way: $O(1)$ on $m\le1/u$; at most
$u\sum_{m\le m_\alpha/2}e^{2mu}\le e^{2m_\alpha u}/(e\,m_\alpha u)$ on
$1/u<m\le m_\alpha/2$, by $m^{-2}\le u/m$ and the display above; and at
most $(4/m_\alpha^{2})\sum_{m\le m_\alpha}e^{2mu}
\le4e^{2m_\alpha u}/(m_\alpha^{2}u)$ on $m>m_\alpha/2$.  In the tail
$m=m_\alpha+1+j$, $q_{\alpha,m}\le4^{-j}$ as in
Corollary~\ref{cor:boundary-constant} and
$e^{2mu}\le e^{2m_\alpha u}e^{1+j}$ on $u\le\tfrac12$, so the tail is
dominated by the geometric series
$\sum_j(e/4)^je^{2m_\alpha u}/m_\alpha=O\bigl(e^{2m_\alpha u}/m_\alpha\bigr)$.
Collecting terms proves \eqref{eq:Wa-upper}.  For \eqref{eq:Wa-flat}:
on $a\le r\le a(1+2/m_\alpha)$ the same decomposition applies with
$e^{2mu}\le e^{4}$ for $m\le m_\alpha$, so
$\sum_{m_1<m\le m_\alpha}(e^{2mu}-1)/m\le2u\,e^{4}m_\alpha\le4e^{4}$ and
the sum is $\log m_\alpha+O(1)$ up to constants.  On $0\le r\le a$, every
$J_m(kr)$ with $m>m_1$ is increasing in $r$ (its first critical point
exceeds $m>2ka$), so $T_m(r)\le T_m(a)$ and
$W_\alpha(r)\le W_\alpha(a)+C'\le\frac{2}{C_k}\log m_\alpha+C$ by
Corollary~\ref{cor:boundary-constant}.

\emph{Lower bound.}  Let $\tfrac{2}{m_\alpha}\le u\le m_\alpha^{-1/2}$, so
that $1/u\le m_\alpha/2$.  On $m_1<m\le1/u$ use $(1+u)^{2m}\ge1$,
$R_m\ge\frac{1}{C_km}(1-Cm^{-1})$, and
$q_{\alpha,m}\ge1-3\cdot2^{-(m_\alpha-m)}$ from \eqref{eq:filter-bracket}
and $\gamma_{m_\alpha-j}\ge2^{j}\alpha$.  Since
$\sum_{m\le m_\alpha/2}2^{-(m_\alpha-m)}/m\le2^{-m_\alpha/2}\log m_\alpha=o(1)$,
this block is at least $\frac{2}{C_k}\log\frac1u-C$.  On
$m_\alpha-\lceil1/(2u)\rceil<m\le m_\alpha$ use $q_{\alpha,m}\ge\tfrac14$
and $R_m\ge\frac{1}{2C_km}(1+u)^{2m}$.  Here
$(1+u)^{2m}\ge e^{2mu(1-u/2)}\ge e^{2m_\alpha u}e^{-2ju}e^{-1}$, using
$m_\alpha u^{2}\le1$, so the block is at least
\[
 \frac{e^{2m_\alpha u}}{4eC_km_\alpha}\sum_{0\le j<1/(2u)}e^{-2ju}
 \ge\frac{e^{2m_\alpha u}}{4eC_km_\alpha}\cdot\frac{e^{-1}}{2u}
 =c\,\frac{e^{2m_\alpha u}}{m_\alpha u}.
\]
Adding the two blocks and dropping the remaining nonnegative modes proves
\eqref{eq:Wa-lower}.
\end{proof}

\begin{theorem}[Quantitative reconstruction of the disc]\label{thm:disc-recon}
There are $\beta_0$ and, for each fixed $\beta\ge\beta_0$, a constant
$C_\beta$, depending only on $a,k,\lambda,\beta$ and $R_{\max}$, such that the
additive threshold rule
\begin{equation}\label{eq:tau-rule}
 \tau_\alpha=\frac{2}{C_k}\log m_\alpha+\beta
\end{equation}
satisfies, for all sufficiently small $\alpha$,
\begin{equation}
 D\subseteq\{z\in K:W_\alpha(z)\le\tau_\alpha\}
 \quad\text{and}\quad
 \bigl|\{z\in K:W_\alpha(z)\le\tau_\alpha\}\,\triangle\,D\bigr|
 \le\frac{C_\beta}{m_\alpha}
 =O\Bigl(\frac{\log\log(1/\alpha)}{\log(1/\alpha)}\Bigr).
\end{equation}
If moreover $\sup_K|\widehat W-W_\alpha|\le\epsilon$ with
$\epsilon\le\beta-\beta_0$, the same bounds hold for
$\{\widehat W\le\tau_\alpha\}$ with $C_\beta$ replaced by
$C_{\beta+\epsilon}$.
\end{theorem}

\begin{proof}
\emph{Interior containment.}  By \eqref{eq:Wa-flat},
$W_\alpha\le\frac{2}{C_k}\log m_\alpha+C<\tau_\alpha$ on
$\{r\le a(1+2/m_\alpha)\}$ as soon as $\beta_0>C$.  This proves the
inclusion together with an exterior collar of width $2a/m_\alpha$.

\emph{Exterior.}  For $a<r\le R_{\max}$ and $\alpha$ so small that $m_\alpha$
exceeds the index $m_1(R_{\max})$ beyond which \eqref{eq:Rm-uniform} holds on
$[a,R_{\max}]$ with relative error $\le\tfrac12$, the passed mode $m_\alpha$
alone gives, by $q_{\alpha,m_\alpha}\ge\tfrac14$ from
\eqref{eq:filter-bracket},
\begin{equation}\label{eq:exterior-lower}
 W_\alpha(r)\ge\tfrac14R_{m_\alpha}(r)
 \ge\frac{c}{m_\alpha}\Bigl(\frac ra\Bigr)^{2m_\alpha},
 \qquad c>0\ \text{independent of $r$ and $\alpha$}.
\end{equation}
For $u\ge m_\alpha^{-1/2}$ this exceeds
$cm_\alpha^{-1}e^{\sqrt{m_\alpha}}>\tau_\alpha$ eventually, since
$(r/a)^{2m_\alpha}\ge e^{2m_\alpha u(1-u/2)}$ is increasing in $u$.  For
$C_\beta/m_\alpha\le u\le m_\alpha^{-1/2}$, write $u=t/m_\alpha$ with
$t\ge C_\beta$.  Then \eqref{eq:Wa-lower} gives
\[
 W_\alpha\bigl(a(1+u)\bigr)
 \ge\frac{2}{C_k}\bigl(\log m_\alpha-\log t\bigr)
 +c\,\frac{e^{2t}}{t}-C
 \ge\frac{2}{C_k}\log m_\alpha+\beta=\tau_\alpha,
\]
for all $t\ge C_\beta$ once $C_\beta$ is chosen with
$c\,e^{2t}/t\ge\beta+C+\frac{2}{C_k}\log t$ for $t\ge C_\beta$, which is
possible because the left side grows exponentially in $t$ and the
right side logarithmically.  Hence the over-inclusion
$\{W_\alpha\le\tau_\alpha\}\setminus D$ lies in the annulus
$\{a<r<a(1+C_\beta/m_\alpha)\}$, of measure at most
$3\pi a^{2}C_\beta/m_\alpha$ for small $\alpha$.  The rate in $\alpha$
follows from \eqref{eq:malpha-asymptotic}, which gives
$m_\alpha\ge cL_\alpha/\log L_\alpha$ with
$L_\alpha=(1+o(1))\log(1/\alpha)$.

\emph{Empirical sets.}  If $\sup_K|\widehat W-W_\alpha|\le\epsilon$, then
by \eqref{eq:set-sandwich},
$\{W_\alpha\le\tau_\alpha-\epsilon\}\subseteq\{\widehat W\le\tau_\alpha\}
\subseteq\{W_\alpha\le\tau_\alpha+\epsilon\}$, and the two displayed
estimates apply at the thresholds $\tau_\alpha\mp\epsilon$, that is, with
$\beta\mp\epsilon$ in place of $\beta$.  Interior containment needs
$\beta-\epsilon\ge\beta_0$, and the over-inclusion constant becomes
$C_{\beta+\epsilon}$.
\end{proof}

\begin{remark}[Multiplicative threshold rules]
For the multiplicative threshold rule
$\tau_\alpha=(1+\delta)\frac{2}{C_k}\log m_\alpha$ with fixed $\delta>0$,
the two-sided bounds
\eqref{eq:Wa-upper}--\eqref{eq:Wa-lower} localize the exterior crossing at
$2m_\alpha u^{*}=\log\log m_\alpha\,(1+o(1))$, since the crossing requires
$e^{2m_\alpha u}/(m_\alpha u)\asymp\delta\log m_\alpha$, while
Lemma~\ref{lem:picard-blowup} bounds the interior miss by
$O\bigl(m_\alpha^{-(1+\delta)}\bigr)$.  The symmetric difference is then
$O(\log\log m_\alpha/m_\alpha)$.  This is asymptotically larger than the
rate $O(1/m_\alpha)$ obtained by the additive threshold rule, which is why
\eqref{eq:tau-rule} is formulated additively.
\end{remark}

\begin{proposition}[Margin and Hausdorff rate for the disc]\label{prop:disc-margin}
Let $\tau_\alpha$ be the threshold defined by \eqref{eq:tau-rule}, with
$\beta\ge\beta_0$ as in Theorem~\ref{thm:disc-recon}, and assume in
addition that $K$ contains an exterior collar,
$\overline{B_{a+\rho}(0)}\subset K$ for some $\rho>0$.  There are positive
constants $c_1,C_1$, depending only on $a,k,\lambda,\beta,R_{\max}$, such that
for all sufficiently small $\alpha$ the following hold.
\begin{enumerate}
\item[\rm(i)] $\partial_rW_\alpha\ge c_1m_\alpha$ on $[a,R_{\max}]$.
In particular the level set $\{z\in K:W_\alpha(z)=\tau_\alpha\}$ is a single
circle $r=r^{*}_\alpha$, and
\[
 a\Bigl(1+\frac{2}{m_\alpha}\Bigr)\le r^{*}_\alpha
 \le a\Bigl(1+\frac{C_1}{m_\alpha}\Bigr).
\]
\item[\rm(ii)] With $r_0=a/m_\alpha$, the margin
\eqref{eq:transversality} holds on the tube
$\{|r-r^{*}_\alpha|\le r_0\}$ with $p=1$ and $c=c_1m_\alpha$.
\item[\rm(iii)] Consequently, if $\widehat W$ is continuous with
$\norm{\widehat W-W_\alpha}_{L^\infty(K)}\le\epsilon\le c_1a/2$, the
localized empirical level set $\widehat\Gamma^{\rm tube}$ of
Proposition~\ref{prop:hausdorff} satisfies
\[
 d_H\bigl(\widehat\Gamma^{\rm tube},\partial D\bigr)
 \le C_1\,\frac{a}{m_\alpha}
 +\frac{\epsilon}{c_1m_\alpha}
 =O\Bigl(\frac{a}{m_\alpha}\Bigr).
\]
\end{enumerate}
\end{proposition}

\begin{proof}
Fix $m_1\ge\max\{m_0,kR_{\max},(kR_{\max})^2\}$ so large that
\eqref{eq:Rm-uniform} holds on $[a,R_{\max}]$ with relative error at
most $\tfrac12$ for $m>m_1$.  For $m>m_1$ and $0<x\le kR_{\max}$ the
alternating series of $J_m$ has decreasing terms, so
$J_m(x)\ge\tfrac34(x/2)^m/m!$ and $0<J_{m+1}(x)\le(x/2)^{m+1}/(m+1)!$,
whence
\begin{equation}\label{eq:Jm-derivative}
 J_m'(x)=\frac mxJ_m(x)-J_{m+1}(x)\ge\frac{m}{2x}\,J_m(x)>0,
 \qquad 0<x\le kR_{\max} .
\end{equation}

\emph{(i)}  Write $T_m(r)$ for the mode-$m$ term of $W_\alpha$, so that
$T_m'(r)=4\pi kJ_m(kr)J_m'(kr)\,\gamma_m/(\gamma_m+\alpha)^2$.  Every
$T_m'$ with $|m|>m_1$ is nonnegative on $[a,R_{\max}]$ by
\eqref{eq:Jm-derivative}.  The finitely many modes $|m|\le m_1$ obey
$|T_m'|\le4\pi k/\gamma_m$, in total at most a constant $C_3$, since
$|J_m|\le1$ and $|J_m'|\le1$.  For $m_1<m\le m_\alpha$,
\eqref{eq:Jm-derivative}, $q_{\alpha,m}\ge\tfrac14$, and
\eqref{eq:Rm-uniform} give, on $[a,R_{\max}]$,
\[
 T_m'(r)\ \ge\ \frac{m}{r}\,T_m(r)
 \ \ge\ \frac{m}{r}\cdot\frac14\cdot\frac{1}{2C_km}
 \Bigl(\frac ra\Bigr)^{2m}
 \ \ge\ \frac{1}{8C_kR_{\max}},
\]
so that
$\partial_rW_\alpha\ge(m_\alpha-m_1)/(8C_kR_{\max})-C_3\ge c_1m_\alpha$ for
small $\alpha$.  By \eqref{eq:Wa-flat} and $\beta\ge\beta_0$,
$W_\alpha<\tau_\alpha$ on $\{r\le a(1+2/m_\alpha)\}$, while
$W_\alpha(R_{\max})>\tau_\alpha$ by \eqref{eq:exterior-lower}.
Strict monotonicity on $[a,R_{\max}]$ therefore yields exactly one crossing
$r^{*}_\alpha>a(1+2/m_\alpha)$, and no level point elsewhere in $K$.  Since
$r^{*}_\alpha-a=O(1/m_\alpha)$, the circle and the tube of (ii) lie in the
collar $\overline{B_{a+\rho}(0)}\subset K$ for small $\alpha$.  The
upper bound on $r^{*}_\alpha$ is the exterior estimate of
Theorem~\ref{thm:disc-recon}.  As a level circle of a radial smooth
function with nonvanishing radial derivative, $r=r^{*}_\alpha$ is a
simple closed $C^1$ curve.

\emph{(ii)}  The tube $\{|r-r^{*}_\alpha|\le a/m_\alpha\}$ is contained
in $[a(1+1/m_\alpha),R_{\max}]$ for small $\alpha$, where
$\partial_rW_\alpha\ge c_1m_\alpha$.  The mean value theorem gives
\eqref{eq:transversality} with $p=1$ and $c=c_1m_\alpha$.

\emph{(iii)}  By (i) and (ii), Proposition~\ref{prop:hausdorff} applies
to the circle $\Gamma_{\alpha,\tau_\alpha}=\{r=r^{*}_\alpha\}$ whenever
$\epsilon<c_1m_\alpha\cdot r_0=c_1a$, and gives
$d_H(\widehat\Gamma^{\rm tube},\Gamma_{\alpha,\tau_\alpha})
\le\epsilon/(c_1m_\alpha)$.  Since
$d_H(\Gamma_{\alpha,\tau_\alpha},\partial D)=r^{*}_\alpha-a
\le C_1a/m_\alpha$, the triangle inequality proves the claim.
\end{proof}

\begin{remark}[Resolution under the threshold rule]\label{rem:resolution}
Theorem~\ref{thm:disc-recon} and Proposition~\ref{prop:disc-margin} bound
both the misclassified area and the offset of the reconstructed level
circle from $\partial D$ by constants times $a/m_\alpha$, which by
\eqref{eq:malpha-asymptotic} is
$O\bigl(a\log\log(1/\alpha)/\log(1/\alpha)\bigr)$.  The computations in
\cref{sec:num-boundary} show that the constants are moderate in practice.
The guarantee improves only when $\log(1/\alpha)$ does.  Halving the
guaranteed error requires roughly squaring $1/\alpha$ and, by
\cref{sec:disc-coupled}, correspondingly many reliable Fourier modes.
\end{remark}

\section{Numerical study of the disc range test}\label{sec:disc-num}
\subsection{Coupled regularization and Fourier cutoff}\label{sec:disc-coupled}
\begin{figure}[tbhp]
\centering
\includegraphics{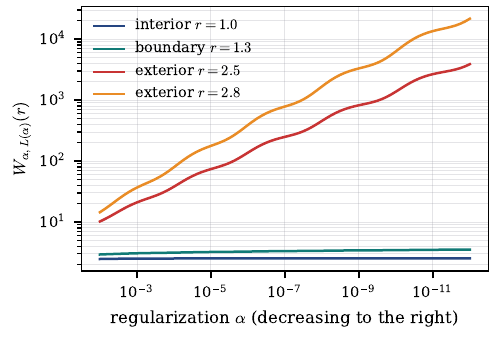}
\caption{Interior stabilization and exterior growth under a joint decrease of
$\alpha$ and increase of the Fourier cutoff.}
\label{fig:disc-convergence}
\end{figure}

With the disc spectral data \eqref{eq:disc-am}--\eqref{eq:disc-kappa} of
\cref{sec:disc-prelim}, the truncated regularized Picard field at radius $r$
is
\begin{equation}\label{eq:disc-W}
 W_{\alpha,L}(r)
 =2\pi\sum_{|m|\le L}
 J_m(kr)^2\frac{\gamma_m}{(\gamma_m+\alpha)^2}.
\end{equation}
At fixed $L$, \eqref{eq:disc-W} has a finite limit as $\alpha\downarrow0$
even at exterior points, so a fixed-dimensional experiment cannot
demonstrate the infinite-dimensional divergence theorem.
Figure~\ref{fig:disc-convergence} therefore couples decreasing
regularization to an increasing Fourier cutoff for $a=1.3$ and
$k=\lambda=1$.  The three regimes are distinct: as
$\alpha$ decreases from $10^{-2}$ to $10^{-12}$ and $L$ increases from
$8$ to $28$, the
interior value $W_{\alpha,L}(1.0)$ stabilizes at $2.535$, the boundary
value $W_{\alpha,L}(1.3)$ increases from $2.93$ to $3.52$, and the exterior
values grow to $3.9\times10^{3}$ at $r=2.5$ and
$2.2\times10^{4}$ at $r=2.8$.

By contrast, with $L=12$ fixed, the exterior values saturate at the
finite fixed-cutoff limits
$W_{0,12}(2.5)=4.319\times10^5$ and $W_{0,12}(2.8)=5.686\times10^6$, where
$W_{0,L}$ denotes the $\alpha\downarrow0$ limit of \eqref{eq:disc-W} at
fixed cutoff $L$.  The exact exterior field $W(r)$ itself is infinite.

Angular quadrature is separately verified at $\alpha=10^{-4}$ against the
analytic Fourier sum.  At the exterior point $r=2.8$, the relative errors for
$M=12,16,20,24$ are $3.45\times10^{-2}$, $1.90\times10^{-5}$,
$2.26\times10^{-9}$, and $2.94\times10^{-13}$ (the last value is at the
double-precision roundoff floor), confirming that increasing
angular resolution, rather than only decreasing $\alpha$, is the correct
continuum experiment.

\subsection{Sharp laws at extreme regularization}\label{sec:num-boundary}
\begin{figure}[tbhp]
\centering
\includegraphics{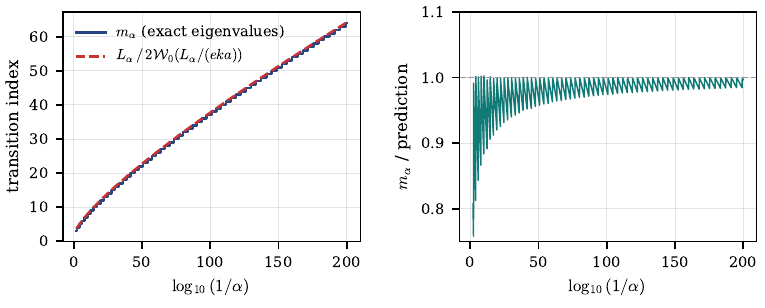}
\caption{Left: the transition index $m_\alpha$ computed from the exact
eigenvalues, against the Lambert-$W$ law \eqref{eq:malpha-asymptotic},
down to $\alpha=10^{-200}$.  Right: their ratio.  The sawtooth reflects
the integer jumps of $m_\alpha$.  By the additive form of
\eqref{eq:malpha-asymptotic} the ratio is $1+O(1/m_\alpha)$.}
\label{fig:lambertw}
\end{figure}

The sharp laws of \cref{cor:boundary-constant},
Lemma~\ref{lem:picard-blowup}, and Theorem~\ref{thm:disc-recon} involve
regularization levels far below double precision, so the spectral
quantities are checked in $50$-digit arithmetic against the exact
eigenvalues \eqref{eq:disc-kappa}.  Only the sublevel-set areas, which
involve no small differences, are computed in double precision through a
scaled evaluation of the filter.  Figure~\ref{fig:lambertw} summarizes
the results for $a=1.3$, $k=\lambda=1$, over $\alpha$ from $10^{-10}$
to $10^{-200}$, where $m_\alpha$ runs from $7$ to $64$.

The boundary defect
$W_\alpha(a)-\frac{2}{C_k}\log m_\alpha$ drifts slowly downward while
remaining bounded (from $1.900$ at $\alpha=10^{-10}$ to $1.647$ at
$\alpha=10^{-200}$), consistent with \eqref{eq:boundary-constant}, and shows
that the first-order term alone underestimates the boundary value at every
practical regularization.  Likewise the interior defect
$W(r)-\frac{2}{C_k}\log\bigl(1/(1-(r/a)^2)\bigr)$, evaluated with the
analytic tail of the Fourier series, remains bounded as $r\uparrow a$,
at $1.152$ for $r/a=1-10^{-5}$, which agrees with the $r\uparrow a$
limit, also $1.152$, to three decimals,
confirming the uniform bound \eqref{eq:picard-blowup}.  Finally, the
symmetric difference for the additive threshold rule with $\beta=2.2$ tracks the
$O(1/m_\alpha)$ bound of Theorem~\ref{thm:disc-recon} with the product $m_\alpha\,|D_\alpha^{\tau_\alpha}\triangle D|/|D|$
rising from $0.49$ and settling near $0.65$ across the range
$\alpha=10^{-10}$ to $10^{-200}$.  The symmetric
difference falls from $6.99\times10^{-2}$ to $1.02\times10^{-2}$ while
$m_\alpha$ grows from $7$ to $64$, the predicted decay law with a
moderate constant.  

\subsection{A reconstruction under data noise}\label{sec:num-recon}
The obstacle is
the impedance disc of the paper translated to the center $z_0=(0.4,-0.3)$,
which is not known to the imaging step.  By the translation formula
$u_\infty^{D+z_0}(\widehat x;d)
=e^{ikz_0\cdot(d-\widehat x)}\,u_\infty^{D}(\widehat x;d)$
the far-field matrix of the translated disc is exact, and the weighted
discrete operator is unitarily equivalent to the centered one, so the
spectral theory of the preceding sections applies verbatim.  Data are
$M=32$ far-field directions, and noise is complex Gaussian, scaled to a
relative Frobenius level $\delta$ of the kernel matrix, with twelve
independent draws per level.  On noisy data the imaginary part of the
discrete $\widetilde F$ is no longer positive semidefinite, so the
imaging step uses $|\Real\widetilde F|+|\Imag\widetilde F|$, positive
semidefinite by construction (round-off negative eigenvalues are set to
zero).  On exact data this is $\Fh$, since $\Imag\widetilde F\ge0$
there (\cref{sec:disc-prelim}).

The imaging procedure is calibrated once, on the exact-data phantom:
$\alpha=\opnorm{\widehat A}/10=0.312$, which retains
$m_\alpha=2$ Fourier mode pairs, and the threshold $\tau=2.09$ that
minimizes the misclassified area on exact data.  The calibrated threshold
is consistent with the boundary law of \cref{cor:boundary-constant}.  It exceeds
$\frac{2}{C_k}\log m_\alpha$ by $1.54$, the defect scale measured above,
so the calibration agrees with the additive threshold rule
\eqref{eq:tau-rule}.  Both parameters are then frozen across every noise level
and draw, so the experiment isolates stability.

\begin{table}[tbhp]
\centering
\caption{Stability of the calibrated reconstruction under relative data
noise $\delta$ (mean $\pm$ standard deviation over twelve draws;
misclassified area relative to $|D|$, Hausdorff distance of the
reconstructed boundary from $\partial D$).}
\label{tab:recon-noise}
\begin{tabular}{ccc}
\toprule
$\delta$ & $|\{\widehat W\le\tau\}\triangle D|/|D|$ &
$d_H(\widehat\Gamma,\partial D)$\\
\midrule
$0$    & $0.0001$          & $0.026$\\
$1\%$  & $0.0019\pm0.0004$ & $0.027\pm0.001$\\
$5\%$  & $0.0094\pm0.0022$ & $0.032\pm0.004$\\
$10\%$ & $0.0179\pm0.0077$ & $0.039\pm0.008$\\
$20\%$ & $0.0299\pm0.0091$ & $0.057\pm0.014$\\
\bottomrule
\end{tabular}
\end{table}

\begin{figure}[tbhp]
\centering
\includegraphics{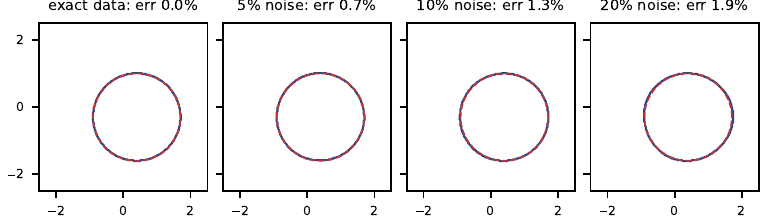}
\caption{Sublevel-set reconstruction of the off-center impedance disc
from $M=32$ far-field directions (first noise draw of each level;
solid blue: the boundary of the reconstructed set $\{\widehat
W\le\tau\}$; dashed red: $\partial D$).  The calibrated pair $(\alpha,\tau)$ is identical in every
panel.}
\label{fig:recon-noise}
\end{figure}

Table~\ref{tab:recon-noise} and Figure~\ref{fig:recon-noise} report the
results: the misclassified area grows from $0.01\%$ on exact data to
$3.0\%$ at $\delta=20\%$, roughly linearly in $\delta$, and the mean
Hausdorff distance of the reconstructed boundary from $\partial D$ stays
below $0.06$, under five percent of the radius.  The exact-data value
$0.026$ is the sampling-grid resolution.  The stability has a spectral
explanation consistent with the rest of the paper.  The retained
eigenvalues ($3.12$, $2.24$, $0.88$) sit an order of magnitude above the
noise floor of the weighted operator ($0.13$ at $\delta=10\%$), so the
few modes that survive the severe ill-posedness are themselves stable
under noise at this level.  Shrinking $\alpha$ buys accuracy only logarithmically, but the
coarse reconstruction that $m_\alpha=2$ mode pairs support is almost
insensitive to noise.  At $\delta=50\%$, where the noise floor reaches the smallest
retained eigenvalue, the misclassified area is still only about $9\%$ and
the obstacle remains clearly localized.

\subsection{Calibration protocol and solver validation}\label{sec:num-kite}
The sharp theory of this paper is disc-specific, but the imaging method
is not, and the experiments beyond the disc share the following setup.
Far-field matrices at $M=64$ directions are generated by a
combined-field integral equation, uniquely solvable at every
wavenumber, discretized by the Nystr\"om method with the standard
quadrature for the logarithmic kernel \cite{CK19}.  The implementation
is validated against the analytic disc at both wavenumbers used below
to a relative error under $4\times10^{-15}$.  For each
wavenumber and noise level the pair $(\alpha,\tau)$ is calibrated on a
\emph{single} noisy draw against the known phantom and then frozen, and
all reported errors are means over twelve fresh draws, so every
experiment measures variation across noise realizations.

\begin{figure}[tbhp]
\centering
\includegraphics{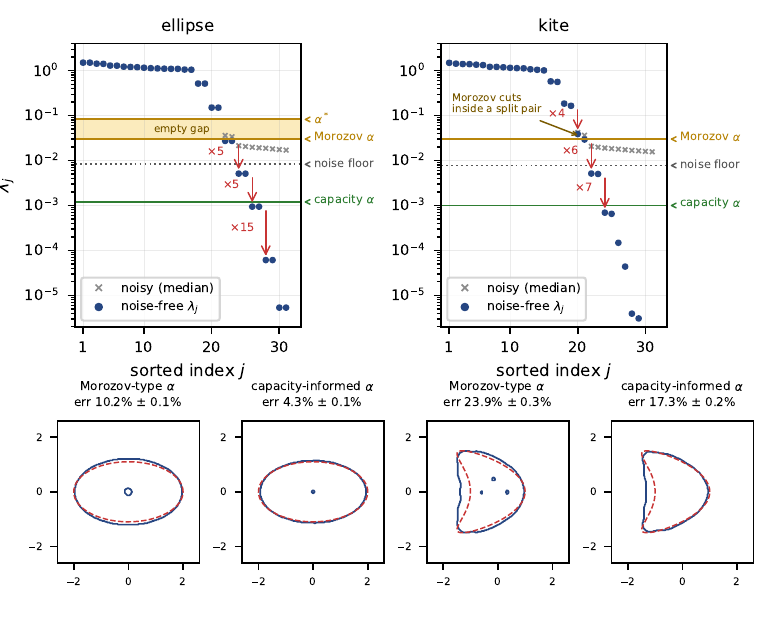}
\caption{\emph{Upper row:} measured eigenvalue ladders (dots) and
median noisy spectra (crosses) for the ellipse and the kite at $k=5$
with $2\%$ noise, with the candidate regularization levels on the
eigenvalue axis.  Arrows mark consecutive inter-pair gaps.  For the
ellipse the noise-proportional level and $\alpha^{*}$, ten times the noise floor,
bracket an \emph{empty} interval, so both retain $21$ modes.  For the kite the $\pm m$ pairs are
split and the noise-proportional level falls inside one.  \emph{Lower row:}
reconstructions from the same data at each level, thresholds
phantom-calibrated at the level shown so the panels isolate the choice
of $\alpha$; solid blue curves are the reconstructed boundaries and dashed red
curves the true boundaries, errors are means over twelve draws.  The noise-proportional panels carry spurious interior
components, three of them for the kite.  At the capacity-informed level
the kite's are gone and the ellipse's is six times smaller.}
\label{fig:worked-example}
\end{figure}

The protocol was also evaluated on the standard kite-shaped impedance
obstacle \cite{CK19}, translated to $(0.25,-0.2)$, with $\lambda=1$,
across noise levels from zero to ten percent at $k=1$ and $k=5$.  Its
summary is the practical reading of Remark~\ref{rem:resolution}.  At $k=1$ the
obstacle diameter is about half a wavelength, two to three mode pairs
carry all the information, and the method returns a stable convex
silhouette: correct location and size, no concavity.  At $k=5$ the
informative band widens to roughly a dozen pairs, the concavity
emerges, and the misclassified area moves only from $15\%$ at exact
data to $20\%$ at ten percent noise, with interior speckle appearing
before any loss of shape.  Resolvable detail is set by how many modes
the spectrum supports, which grows with $ka$ but only logarithmically
with $1/\alpha$.  Reconstructions of the kite appear in
Figure~\ref{fig:worked-example} of \cref{sec:worked-example}, where its
logarithmic capacity enters the choice of the regularization level.

\section{General obstacles}\label{sec:general}

The sharp laws of \cref{sec:disc-sharp,sec:sets} are proved for the disc.
This section extends the transition theory to \emph{every} impedance
obstacle and, for a single analytic boundary component whose conformal
center lies in the obstacle (in particular for every convex analytic
obstacle), identifies the \emph{effective radius}, the radius for which
the disc law reproduces the count, as the logarithmic capacity.
Its use is isolated in \cref{prop:trace-localization}, and
\cref{rem:noncentered-obstruction} identifies the obstruction without it.

\subsection{Reduction to the data-to-pattern operator}

The facts about the far-field operator used below are collected here as
a standing hypothesis, stated once and used throughout, so that the
inputs to \cref{sec:general-capacity} are explicit.

\begin{fact}[Inputs from the factorization theory]\label{fact:fm-inputs}
Let $D$ be an impedance obstacle with connected exterior and boundary
$\Gamma:=\partial D$ of class $C^2$.  Let the real-valued impedance
$\lambda\in L^\infty(\Gamma)$ satisfy $\lambda\ge\lambda_0>0$ almost
everywhere, and let
$G:X\to L^2(\mathbb S^1)$, $X:=H^{-1/2}(\Gamma)$, be the
data-to-pattern operator: $Gh$ is the far field of the radiating
solution $w_h$ of the exterior impedance problem
$\partial_\nu w_h+ik\lambda w_h=h$ on $\Gamma$.  This trace
convention is fixed throughout.
Write $Z_k:=ikM_\lambda:H^{1/2}(\Gamma)\to H^{-1/2}(\Gamma)$, where
$M_\lambda$ is multiplication by $\lambda$.  Then, with
$\gamma_2=e^{i\pi/4}/\sqrt{8\pi k}$ and
$T:H^{1/2}(\Gamma)\to H^{-1/2}(\Gamma)$ bounded and boundedly
invertible:
\begin{enumerate}
\item[\rm(F1)] \emph{Factorization}: $F=-\bar\gamma_2^{-1}GT^{*}G^{*}$,
 with $T$ the impedance middle operator of Chamaillard et al.~\cite{CCH13}.
\item[\rm(F2)] \emph{Range identity}: $\ran(\Fh^{1/2})=\ran(G)$, where
 $\Fh=|\Real\widetilde F|+\Imag\widetilde F$ and
 $\widetilde F=e^{-i\pi/4}F$.  This is the generalized factorization
 theorem of Chamaillard et al.~\cite{CCH13}, whose hypotheses hold here with boundary
 space $V(\Gamma)=L^2(\Gamma)$, the required compact embedding being
 $H^{1/2}(\Gamma)\Subset L^2(\Gamma)$, and $k^{2}$ no eigenvalue of the
 interior problem $\Delta u+k^{2}u=0$ in $D$, $\partial_\nu u\pm Z_ku=0$
 on $\Gamma$.  The last is immediate
 from $\lambda\ge\lambda_0>0$.  Green's first identity gives
 $\int_D(|\nabla u|^2-k^2|u|^2)=\int_\Gamma\partial_\nu u\,\bar u
 =\mp ik\int_\Gamma\lambda|u|^2$, whose left side is real, so
 $\int_\Gamma\lambda|u|^2=0$ and $u|_\Gamma=0$.  Then
 $\partial_\nu u|_\Gamma=0$ and $u\equiv0$ by unique continuation.
 The abstract antecedent is Kirsch and Grinberg~\cite[Theorem~2.15]{KG08}, stated for
 $\Fh=|\Real(e^{it}F)|+\Imag F$ with a free rotation parameter $t$.
 Kirsch and Grinberg~\cite[Section~2.5.2]{KG08} take $t=0$ in three dimensions, where no
 rotation is needed.  Here the rotation is forced by the normalization.
 Writing $F_{\rm st}$ and $G_{\rm st}$ for the operators formed in the
 convention of Chamaillard et al.~\cite{CCH13}, in which
 $u^{s}(x)=\gamma_2|x|^{-1/2}e^{ik|x|}\bigl(u^{\infty}_{\rm st}(\widehat
 x)+O(|x|^{-1})\bigr)$ against the physical
 $u^{s}(x)=|x|^{-1/2}e^{ik|x|}\bigl(u^{\infty}(\widehat
 x)+O(|x|^{-1})\bigr)$, comparison of the two gives
 $u^{\infty}=\gamma_2u^{\infty}_{\rm st}$ and hence
 $F=\gamma_2F_{\rm st}$, $G=\gamma_2G_{\rm st}$.  Therefore
 $\widetilde F=e^{-i\pi/4}\gamma_2F_{\rm st}=F_{\rm st}/\sqrt{8\pi k}$
 is a \emph{positive real} multiple of $F_{\rm st}$, so
 $\Fh=\Fh^{\rm st}/\sqrt{8\pi k}$ is a positive multiple of the
 standard-normalized symmetrization, while $\ran(G)=\ran(G_{\rm st})$
 because $\gamma_2\ne0$.  Multiplication by a positive scalar changes
 neither the range of a square root nor the identity itself.  Without
 the rotation the two parts of $\Fh$ would carry different phases.  The
 two-dimensional phase shift is that of Kirsch~\cite[Section~5]{Kir98}.
\item[\rm(F3)] \emph{Injectivity} \cite[Theorem~2.4(b)]{KG08}: $\ker G=\{0\}$,
 by uniqueness for the exterior impedance problem and Rellich's lemma.  The
 statement in Kirsch and Grinberg~\cite{KG08} is three-dimensional, and is obtained there by
 reduction to the Dirichlet data-to-pattern operator.  The direct argument
 just indicated uses only uniqueness for the exterior impedance problem and
 Rellich's lemma, and is independent of the dimension.
\item[\rm(F4)] \emph{Adjoint form}: with $B$ the map sending a Herglotz
 density to its impedance trace on $\Gamma$, the standard-normalized
 adjoint identity $B^{*}=G_{\rm st}T$ of Chamaillard et al.~\cite{CCH13} and
 $G=\gamma_2G_{\rm st}$ give
 \[
  B^{*}=\gamma_2^{-1}GT,
  \qquad\text{equivalently}\qquad
  G^{*}=\overline{\gamma_2}\,T^{-*}B .
 \]
 This is separate from (F1), which determines $G$ only through the
 composite $GT^{*}G^{*}$, whereas \eqref{eq:S-column} needs $G^{*}$.
\end{enumerate}
\end{fact}

 (F2) gives the Douglas factors $C_1,C_2$ in
\cref{prop:reduction}; (F3) closes the index-zero Fredholm argument in
\cref{prop:S-invertible}; (F1) fixes the constant $\gamma_2$ in
\eqref{eq:sandwich}; and (F4) gives the coordinate formula in
\cref{lem:S-coordinate}.  (F1), (F2) and (F4) rest on
Chamaillard et al.~\cite{CCH13} and (F3) on Kirsch and Grinberg~\cite{KG08}, with the two-dimensional phase
normalization of Kirsch~\cite{Kir98}.  The bookkeeping that carries them from
the standard convention to the physical one is the computation
displayed in (F2).

\begin{proposition}[Reduction to the data-to-pattern operator]\label{prop:reduction}
For every impedance obstacle, there are constants $0<c_1\le C$, depending on
$D$, $k$, and $\lambda$, such that
\begin{equation}\label{eq:reduction}
 c_1\,\sigma_j(G)\le\sigma_j\bigl(\Fh^{1/2}\bigr)\le C\,\sigma_j(G),
 \qquad j\ge1 .
\end{equation}
In particular $\log\lambda_j(\Fh)=2\log\sigma_j(G)+O(1)$.  Thus the
symmetrization changes the singular-value scale only by bounded factors;
the dependence on the impedance remains contained in $G$.
\end{proposition}

\begin{proof}
By (F2) of \cref{fact:fm-inputs}, $\ran(\Fh^{1/2})=\ran(G)$.
By Douglas' theorem on range inclusion and
factorization \cite{Dou66} there are bounded operators $C_1,C_2$ with
$\Fh^{1/2}=GC_1$ and $G=\Fh^{1/2}C_2$, whence
$\sigma_j(\Fh^{1/2})\le\opnorm{C_1}\,\sigma_j(G)$ and
$\sigma_j(G)\le\opnorm{C_2}\,\sigma_j(\Fh^{1/2})$ for every $j$
\cite{GK69}.  Since $\Fh^{1/2}$ is self-adjoint and nonnegative,
$\lambda_j(\Fh)=\sigma_j(\Fh^{1/2})^2$.
\end{proof}

\subsection{A two-sided transition law}
The next two lemmas bound $\sigma_j(G)$ from both sides.  The upper
bound cuts the radiated field on a circumscribed circle.  The lower bound
tests $G$ with multipoles placed at an interior point.  No monotonicity
or comparison principle is needed because a multipole \emph{is} the
exterior solution for its own impedance data.

\begin{lemma}[Radiating tail bound]\label{lem:tail}
Let $\overline D\subset\overline{B_R(z_1)}$ and $R'>R$.  There are $C$ and $m_2$,
depending on $D,k,\lambda,R'$, such that
\begin{equation}\label{eq:tail-bound}
 \sigma_{2m+2}(G)\le C\sqrt m\,\Bigl(\frac{ekR'}{2m}\Bigr)^{m},
 \qquad m\ge m_2 .
\end{equation}
\end{lemma}

\begin{proof}
After the unitary modulation used in \cref{sec:num-recon} we may take
$z_1=0$.  For $h\in X$ the radiating field $w_h$ has, for $r\ge R''$ with
$R''\in(R,R')$, an expansion
$w_h=\sum_{n\in\mathbb Z}b_nH_n^{(1)}(kr)e^{in\theta}$ converging in
$L^2$ on circles \cite{CK19}, and in the normalization of
\cref{sec:disc-prelim} its far field has Fourier coefficients of squared
modulus $(4/k)|b_n|^2$.  The map $h\mapsto w_h|_{\partial B_{R'}}$ is
bounded by the well-posedness of the exterior impedance problem
\cite{CK19}, which for $\lambda\in L^\infty(\Gamma)$ with
$\lambda\ge\lambda_0>0$ follows from the variational formulation on a
truncated annulus with a Dirichlet-to-Neumann condition, coercive up to
a compact perturbation and injective by the sign of $\lambda$:
$\norm{w_h}_{L^2(\partial B_{R'})}\le C_G\norm h$.  Let $\mathcal P_m$ denote the
Fourier projection onto $\operatorname{span}\{\psi_n:|n|\le m\}$ in
$L^2(\mathbb S^1)$.  Since $|H_n^{(1)}(kR')|$ is nondecreasing in $|n|$
for $|n|\ge m_2(kR')$ by the large-order asymptotics of
Lemma~\ref{lem:gamma-ratio},
\[
 \norm{(I-\mathcal P_m)Gh}^2
 =\frac4k\sum_{|n|>m}|b_n|^2
 \le\frac{2}{\pi kR'}\,
 \bigl|H_{m+1}^{(1)}(kR')\bigr|^{-2}
 \norm{w_h}^2_{L^2(\partial B_{R'})},
\]
and since $\mathcal P_mG$ has rank at most $2m+1$,
$\sigma_{2m+2}(G)\le\norm{(I-\mathcal P_m)G}
\le C|H_{m+1}^{(1)}(kR')|^{-1}$, which together with
$|H_{m}^{(1)}(x)|=\sqrt{2/(\pi m)}\,(2m/(ex))^{m}(1+o(1))$ gives
\eqref{eq:tail-bound}.
\end{proof}

\begin{lemma}[Multipole lower bound]\label{lem:multipole}
Let $B_d(z_0)\subset D$, with $\partial D$ of class $C^2$ as in
\cref{fact:fm-inputs}.  There are $c_1>0$ and $m_3$, depending on $D,k,\lambda$,
such that
\begin{equation}\label{eq:multipole-bound}
 \sigma_{2m+1}(G)\ge c_1\,m^{-1}\Bigl(\frac{ekd}{2m}\Bigr)^{m},
 \qquad m\ge m_3 .
\end{equation}
\end{lemma}

\begin{proof}
As before we may take $z_0=0$.  For
$n\in\mathbb Z$ let $w_n(x):=H_n^{(1)}(k|x|)e^{in\arg x}$, a radiating
Helmholtz solution on $\mathbb R^2\setminus\{0\}$, smooth in a
neighborhood of $\partial D$ because
$\operatorname{dist}(0,\partial D)\ge d$, and set
\[
 h_n:=\bigl(\partial_\nu w_n+ik\lambda w_n\bigr)\big|_{\partial D} .
\]
Since $\partial D$ is $C^2$ the normal field is continuous, so
$\partial_\nu w_n\in C(\partial D)$.  With $\lambda\in L^\infty(\Gamma)$
from \cref{fact:fm-inputs},
$h_n\in L^\infty(\partial D)\subset L^2(\partial D)\subset H^{-1/2}(\partial D)$,
with $\norm{h_n}_{X}\le C_0\norm{h_n}_{L^2(\partial D)}$.  By uniqueness
of the exterior impedance problem, the radiating solution with data $h_n$
\emph{is} $w_n$, so $Gh_n$ is the far field of $w_n$:
\[
 (Gh_n)(\widehat x)
 =\sqrt{\frac{2}{\pi k}}\,e^{-i\pi/4}(-i)^n e^{in\theta_{\widehat x}} .
\]
These functions are pairwise orthogonal in $L^2(\mathbb S^1)$ with the
$n$-independent squared norm $4/k$.  In particular the $h_n$ are linearly
independent and $V_m:=\operatorname{span}\{h_n:|n|\le m\}$ has dimension
$2m+1$.

We estimate $\norm{h_n}_{L^2(\partial D)}$.  On $\partial D$ the radial
variable ranges in $[d,R_1]$, $R_1:=\sup_{x\in\partial D}|x|$.  In polar
coordinates,
\[
 |\nabla w_n(x)|
 \le k\bigl|H_n^{(1)\prime}(k|x|)\bigr|
 +\frac{|n|}{d}\,\bigl|H_n^{(1)}(k|x|)\bigr|,
\]
and the recurrence
$H_n^{(1)\prime}=\tfrac12(H_{n-1}^{(1)}-H_{n+1}^{(1)})$ together with the
large-order asymptotics of Lemma~\ref{lem:gamma-ratio}, uniform on the
compact set $[kd,kR_1]$, gives
$|H_{n\pm1}^{(1)}(x)|\le C(n/x)^{\pm1}|H_n^{(1)}(x)|$ there, and
$\sup_{[kd,kR_1]}|H_n^{(1)}|\le2\,|H_n^{(1)}(kd)|$ for $n\ge m_3$, since
$|H_n^{(1)}(x)|=\sqrt{2/(\pi n)}\,\bigl(2n/(ex)\bigr)^{n}(1+o(1))$ is
decreasing in $x$ up to a factor $1+o(1)$.  Hence
\[
 \norm{h_n}_{L^2(\partial D)}
 \le|\partial D|^{1/2}
 \sup_{\partial D}\bigl(|\nabla w_n|+k\norm{\lambda}_{L^\infty(\partial D)}|w_n|\bigr)
 \le C\,\frac nd\,\bigl|H_n^{(1)}(kd)\bigr|
 \le C\,n^{1/2}\Bigl(\frac{2n}{ekd}\Bigr)^{n},
 \qquad n\ge m_3,
\]
and the same bound holds for $\norm{h_n}_X$.  For
$h=\sum_{|n|\le m}a_nh_n\in V_m$, exact orthogonality gives
$\norm{Gh}^2=(4/k)\sum|a_n|^2$, while by Cauchy--Schwarz
\[
 \norm h_X\le\Bigl(\sum|a_n|^2\Bigr)^{1/2}
 \Bigl(\sum_{|n|\le m}\norm{h_n}_X^2\Bigr)^{1/2}
 \le C\Bigl(\sum|a_n|^2\Bigr)^{1/2}\,m\,\Bigl(\frac{2m}{ekd}\Bigr)^{m},
\]
the last sum being dominated by its largest terms.  The Courant--Fischer
characterization
$\sigma_{2m+1}(G)\ge\min_{0\ne h\in V_m}\norm{Gh}/\norm h_X$ then proves
\eqref{eq:multipole-bound}.
\end{proof}

\begin{theorem}[Two-sided transition law for every obstacle]\label{thm:general-twosided}
Let $D$ be an impedance obstacle with $B_d(z_0)\subset D$ and
$\overline D\subset\overline{B_R(z_1)}$ (closure so that $R$ may be
the exact circumradius), and let $R'>R$.  Then, as $\mu\to\infty$,
\begin{equation}\label{eq:general-twosided}
 \mu\log\frac{\mu}{ekR'}-O(\log\mu)
 \;\le\;\log\frac1{\lambda_\mu(\Fh)}
 \;\le\;\mu\log\frac{\mu}{ekd}+O(\log\mu),
\end{equation}
and consequently, with $L_\alpha=\log(C_k/\alpha)$ as above (the
constant $\log C_k$ being absorbed in the $O(1)$ below) and
$N(\rho):=L_\alpha/\LambertW\bigl(L_\alpha/(ek\rho)\bigr)$, the number
$N_\alpha=\#\{j:\lambda_j(\Fh)\ge\alpha\}$ of passed modes satisfies,
as $\alpha\downarrow0$,
\begin{equation}\label{eq:general-count}
 N(d)-O(1)\;\le\;N_\alpha\;\le\;N(R')+O(1) .
\end{equation}
For the disc, $d=R=a$, and letting $R'\downarrow a$ recovers the
leading logarithmic form of \cref{thm:disc-sharp}.  It does not recover
its additive law.  The upper bound is proved for each fixed $R'>R$, and
$N(R')-N(a)\asymp L_\alpha(\log L_\alpha)^{-2}\log(R'/a)$ diverges as
$\alpha\downarrow0$ for every such $R'$, so the $O(1)$ remainder of
\cref{thm:disc-sharp} lies beyond this route.  Uniformity as
$R'\downarrow a$ is not established here.
\end{theorem}

\begin{proof}
By Proposition~\ref{prop:reduction} and Lemma~\ref{lem:tail}, at
$\mu=2m+2$,
$\lambda_\mu(\Fh)\le C\sigma_{2m+2}(G)^2\le Cm(ekR'/2m)^{2m}$, so
$\log(1/\lambda_\mu)\ge2m\log\bigl(2m/(ekR')\bigr)-O(\log m)
=\mu\log\bigl(\mu/(ekR')\bigr)-O(\log\mu)$.  By
Proposition~\ref{prop:reduction} and Lemma~\ref{lem:multipole}, at
$\mu=2m+1$,
$\lambda_\mu(\Fh)\ge c_1\,m^{-2}(ekd/2m)^{2m}$, which gives the other side.
Both extend from these parities to every $\mu$ by the monotonicity of
$\lambda_\mu$, at the cost of enlarging the $O(\log\mu)$.  For
\eqref{eq:general-count}: by \eqref{eq:general-twosided},
$\lambda_\mu<\alpha$ as soon as
$\mu\log\bigl(\mu/(ekR')\bigr)\ge L_\alpha+C\log\mu$, and solving this
threshold equation through $\LambertW$ as in the proof of
Theorem~\ref{thm:disc-sharp} (the derivative of the inverse map being
$1/(1+\LambertW)=O(1/\log L_\alpha)$ in the variable $\mu$, and half
that in $m=\mu/2$) shows that the
$C\log\mu$ correction moves the solution by $O(1)$.  Hence
$N_\alpha\le N(R')+O(1)$.  Likewise $\lambda_\mu\ge\alpha$ whenever
$\mu\log\bigl(\mu/(ekd)\bigr)\le L_\alpha-C\log\mu$, giving
$N_\alpha\ge N(d)-O(1)$.
\end{proof}

\begin{remark}[Scope of the two-sided law]
Arguments that factor through $F$ alone
introduce separate index losses for the incident and radiating
smoothing and reach only half of the exponent in
\eqref{eq:general-twosided}.  The full exponent is restored by the
reduction \eqref{eq:reduction}, which follows from the range identity.
The multipole construction uses
nothing about the boundary condition beyond uniqueness of the exterior
problem, so Lemma~\ref{lem:multipole} and the theorem hold verbatim for
Dirichlet and Neumann obstacles with their factorizations.  For those
boundary conditions, inclusion-based comparisons are also available
through the monotonicity relations of Albicker and Griesmaier~\cite{AG20}.
\end{remark}

\subsection{The capacity law}\label{sec:general-capacity}
Theorem~\ref{thm:general-twosided} brackets the mode count between the
law for the inradius and the law for any radius exceeding the
circumradius.  For a single analytic boundary component whose conformal
center lies in the obstacle, the effective radius is the
\emph{logarithmic capacity}
$\operatorname{cap}(\overline D)$ \cite{Ran95}.  The exterior conformal map
$\Psi$ of $\overline D$, normalized as
$\Psi(w)=cw+c_0+c_{-1}w^{-1}+\cdots$ on $|w|>1$, has
$\operatorname{cap}(\overline D)=c$.  From here on $c$ denotes the
capacity and generic constants are written $c_1,c_2,C$.  For the disc of radius $a$ the
capacity is $a$, and for an ellipse it is the mean of the two
semiaxes.  Since $B_d(z_0)\subset D\subset B_R(z_1)$ implies
$d\le\operatorname{cap}(\overline D)\le R$ by monotonicity of capacity
\cite{Ran95}, this value is admissible in
\eqref{eq:general-twosided}.

Assume in this subsection that $D$ is bounded and simply connected with
analytic Jordan boundary $\Gamma$ and that the impedance
$\lambda\in C^\omega(\Gamma)$ satisfies $\lambda\ge\lambda_0>0$.  Constant
impedance is included.  Translating $D$ (which conjugates $F$, hence
$\Fh$, by a unitary modulation as in \cref{sec:num-recon}), we may take
$c_0=0$ in the expansion of $\Psi$.  Analyticity of $\Gamma$ implies that
$\Psi$ and $\Psi^{-1}$ extend holomorphically to an annulus $|w|>\rho$
for some $\rho<1$.  Eliminating $G$ between (F1) and (F4) of \cref{fact:fm-inputs} gives the exact factorization
\begin{equation}\label{eq:sandwich}
 F=-\gamma_2 B^{*}T^{-1}B,
 \qquad B=Z_kH_0+H_1,
\end{equation}
where $H_0g$ and $H_1g$ are the Dirichlet and Neumann traces on $\Gamma$
of the Herglotz wave
$v_g(x)=\int_{\mathbb S^1}e^{ikx\cdot d}g(d)\,\dd s(d)$ and
$Z_k=ik\lambda$ acts by multiplication, so that $Bg$ is the impedance
trace of $v_g$ as in (F4) of \cref{fact:fm-inputs}, and
$T:H^{1/2}(\Gamma)\to H^{-1/2}(\Gamma)$ is the impedance middle
operator, an isomorphism for positive impedance \cite{CCH13}.  For
$m\ge1$ put
\begin{equation}
 d_m:=C_k\Bigl(\frac{ekc}{2m}\Bigr)^{2m},
 \qquad
 \mathcal H_m:=\operatorname{span}\{\psi_m,\psi_{-m}\},
\end{equation}
so that $d_{m+1}/d_m=O(m^{-2})$, and let $Q_M$ denote the orthogonal
projection onto $\bigoplus_{m\ge M}\mathcal H_m$.

\begin{definition}\label{def:centered}
Call $D$ \emph{centered} if its conformal center lies in $D$.  In the
normalization $c_0=0$ of this subsection, $0\in D$.
\end{definition}

Before the translation, the conformal center is the coefficient $c_0$
of $\Psi$, and it is the barycenter of the equilibrium measure of
$\overline D$.  That measure is the push-forward under $\Psi$ of
normalized arclength on $|w|=1$ \cite{Ran95}, whose mean the Laurent
expansion evaluates as $c_0$.  The measure has full support on
$\Gamma$, since $\Psi$ carries $|w|=1$ homeomorphically onto $\Gamma$,
and a barycenter lies in the closed convex hull of the support.  Were
the barycenter of a probability measure to lie on a supporting line of
that hull, its support would be contained in the line, which a Jordan
curve is not.  Since every boundary point of a closed convex set
carries a supporting line, $c_0$ lies in the \emph{interior} of the
convex hull of $\Gamma$.  For convex $D$ that hull is $\overline D$
and its interior is $D$, so every convex obstacle is centered.  Every
obstacle used below is centered: for the kite the equilibrium
barycenter, computed in the reproduction suite, is $(-0.4493,0)\in D$,
and the nonconvex four-lobed domain
$\psi(w)=1.3\,(w+0.3w^{-1}+0.2w^{-3})$ of \cref{tab:capacity} has
$\min_t|\psi(e^{it})|=0.98$, so $0\in D$ there too.

Centeredness enters only through
\cref{prop:trace-localization} below, which supplies the trace
estimate \eqref{eq:rmode-profile}, and there only through the single
estimate $|P_n(0)|\le C\vartheta_0^n$; it is assumed from here to the
end of this subsection, and it is \emph{not} needed for
Theorem~\ref{thm:general-twosided}, whose two-sided bracket holds for
every obstacle.

The orthogonalization of the mode family needs a preparatory
statement.  The family is nearly orthonormal, so a Gram correction is
available in many forms; but the correction is applied inside the
weight $D_c$, and conjugation by $D_c^{1/2}$ multiplies the entry
$(n,m)$ by $\sqrt{d_n/d_m}$, which is \emph{super-exponentially large}
when $n<m$.  A correction with entries above the diagonal is therefore
inadmissible, however small it is before conjugation.  In particular
the symmetric polar factor $\mathcal Gr^{-1/2}$, which mixes in both
directions, is not usable here.  What is usable is a block-lower
triangular factor, whose entries above the diagonal vanish identically.

\begin{lemma}[One-sided inverse Cholesky factor]
\label{lem:one-sided-cholesky}
Let $\mathcal Gr=I+H$ be a positive invertible block matrix on
$\bigoplus_{m\ge M}\mathbb C^2$ and suppose
\begin{equation}\label{eq:Gram-class}
 \opnorm{H_{nm}}
 \le C\frac{r^{|n-m|}}{\min(n,m)},
 \qquad n,m\ge M,
\end{equation}
with $\opnorm H<1/2$.  After enlarging the starting index $M$ once,
there is a boundedly invertible block-lower triangular matrix
$\mathcal Ch=I+Z$ such that
\begin{equation}\label{eq:inverse-Cholesky}
 \mathcal Gr^{-1}=\mathcal Ch\mathcal Ch^*,
 \qquad
 \mathcal Ch^*\mathcal Gr\mathcal Ch=I,
\end{equation}
and, for some $r_1\in(r,1)$,
\begin{equation}\label{eq:T-triangular-bound}
 Z_{nm}=0\quad(n<m),
 \qquad
 \opnorm{Z_{nm}}
 \le C_1\frac{r_1^{n-m}}{m}\quad(n\ge m).
\end{equation}
\end{lemma}

\begin{proof}
First convert the hypothesis \eqref{eq:Gram-class} into the class
$\mathcal B^{-1}$ required by \cref{lem:classes}, whose weight sits on
a single index.  Where $|n-m|\le\min(n,m)/2$ the two indices are
comparable, so $\min(n,m)^{-1}\asymp(1+n)^{-1}$.  Where
$|n-m|>\min(n,m)/2$ one has $\max(n,m)<3|n-m|$, so the polynomial
ratio is absorbed by passing to a strictly larger base
$r_0\in(r,1)$.  Hence $H\in\mathcal B^{-1}_{r_0}$.  The enlargement is
essential, since $(1+n)\,r^{|n-m|}$ is unbounded at the original base.
Now $\opnorm H<\tfrac12$, so
\cref{lem:classes}(ii) makes $\mathcal Gr^{-1}$ exponentially
localized, and \cref{lem:classes}(i) applied to
$\mathcal Gr^{-1}-I=-\mathcal Gr^{-1}H$ gives
$\mathcal Gr^{-1}=I+R$ with $R\in\mathcal B^{-1}_{r_1}$ for some
$r_1\in(r_0,1)$.  Since
$\mathcal Gr$ is self-adjoint, so is $R$.  The same bound therefore
holds with $m$ in place of $n$, and taking the smaller of the two
recovers the symmetric form
\begin{equation}\label{eq:G-inverse-entry}
 \opnorm{R_{nm}}
 \le C_2\frac{r_1^{|n-m|}}{\min(n,m)}.
\end{equation}
We construct the lower Cholesky factor of $I+R$.  For a finite leading
section, choose the diagonal blocks of $\mathcal Ch$ positive and use the
standard recursion
\begin{align}
 \mathcal Ch_{mm}\mathcal Ch_{mm}^*
 &=(I+R)_{mm}-\sum_{k<m}\mathcal Ch_{mk}\mathcal Ch_{mk}^*,
 \label{eq:Cholesky-diagonal}\\
 \mathcal Ch_{nm}\mathcal Ch_{mm}^*
 &=(I+R)_{nm}-\sum_{k<m}\mathcal Ch_{nk}\mathcal Ch_{mk}^*,
 \qquad n>m.\label{eq:Cholesky-offdiag}
\end{align}
Positivity of every leading section makes each diagonal square root
well defined.  We verify bounds independent of the section size.

From \eqref{eq:G-inverse-entry}, the diagonal right side in
\eqref{eq:Cholesky-diagonal} is $I+O(m^{-1})$.  Assuming the estimate
\eqref{eq:T-triangular-bound} for earlier columns, the sum in that
formula is
\[
 O\left(\sum_{k<m}\frac{r_1^{2(m-k)}}{k^2}\right)=O(m^{-2});
\]
the part $k<m/2$ is exponentially small and the rest is a geometric
sum.  Functional calculus for the positive square root therefore gives
$\mathcal Ch_{mm}=I+O(m^{-1})$ and a uniformly bounded inverse.

For $n>m$, insert the inductive estimates in
\eqref{eq:Cholesky-offdiag}.  The product sum is bounded by
\[
 C\sum_{k<m}\frac{r_1^{n-k}}{k}
                    \frac{r_1^{m-k}}{k}
 \le C'\frac{r_2^{n-m}}{m^2}
\]
after replacing $r_1$ by $r_2\in(r_1,1)$.  This is smaller than the
bound $C r_2^{n-m}/m$ for $(I+R)_{nm}$.  Multiplication by the uniformly
bounded $\mathcal Ch_{mm}^{-*}$ proves
\eqref{eq:T-triangular-bound}.  Induction establishes the estimate for
every finite leading section.

The row and column sums in \eqref{eq:T-triangular-bound} are $O(M^{-1})$.
After increasing the fixed starting index once, every finite factor has
$\opnorm{\mathcal Ch-I}\le1/2$.  Hence its inverse is uniformly bounded by
the Neumann series.  A diagonal subsequence converges entrywise and
strongly on finite sequences to a lower triangular $\mathcal Ch$ with
the same estimates.  Passing to the limit in the finite Cholesky
identity gives $\mathcal Gr^{-1}=\mathcal Ch\mathcal Ch^*$.  Both factors
are invertible because $\opnorm{\mathcal Ch-I}\le1/2$, and inversion gives the second identity in
\eqref{eq:inverse-Cholesky}.
\end{proof}

\begin{lemma}[Faber--Herglotz coordinate construction]
\label{lem:faber-coordinates}
Let $D$ be centered.  There are $M_0$, $\rho_1<1$, Herglotz
densities $g_{\pm m}$, $m\ge M_0$, a block-lower triangular
operator $\mathcal Ch$, and a unitary $U$ on
$L^2(\mathbb S^1)$ with the following properties:
\begin{enumerate}
\item
\[
 \left|\ip{g_{\pm m}}{\psi_{\pm j}}-\delta_{mj}\right|
 \le \frac{C}{m}\,2^{-(m-j)/2},
 \qquad 0\le j\le m .
\]
\item The vectors $g_{\pm m}\mathcal Ch$ form the orthonormal tail
$U^*\psi_{\pm m}$, and
\[
 \mathcal Ch_c:=D_c^{1/2}\mathcal ChD_c^{-1/2}=I+Z_c,
 \qquad
 \opnorm{(Z_c)_{nm}}\le \opnorm{Z_{nm}}\quad(n\ge m).
\]
\item The orthogonal complement of the tail is finite-dimensional and
has an orthonormal basis whose Fourier coefficients decay geometrically.
The constants in that decay can be chosen uniformly over the basis.
\end{enumerate}
\end{lemma}

\begin{proof}
\emph{1. Faber boundary modes.}
Let $P_m$ be the $m$-th Faber polynomial of $\overline D$, normalized by
$P_m(\Psi(w))=w^m+\sum_{n\ge1}b_{mn}w^{-n}$.  Moving the contour in the
generating formula for $P_m$ from $|w|=1$ to $|w|=\rho$ gives
$b_{mn}=O(\rho^{m+n})$, whence, for every fixed Sobolev index $s$,
\begin{equation}\label{eq:grunsky}
 \norm{P_m\!\circ\!\Psi-e^{imt}}_{H^s(\mathbb S^1)}
 +\norm{\partial_r(P_m\!\circ\!\Psi)-me^{imt}}_{H^s(\mathbb S^1)}
 \le C_sm^s\rho^m,
\end{equation}
the radial derivative taken at $|w|=1$ from the exterior.  This is the
Faber--Grunsky layer-potential calculation of
Jung and Lim~\cite[Sections~2--3]{JL21}.

\emph{2. Helmholtz continuation and its Herglotz density.}
Let $\mathcal V_k$ be the Vekua transform \cite{MHP11}, normalized so
that
$\mathcal V_k(z^m)=(2^mm!/k^m)J_m(kr)e^{im\theta}=z^m(1+O(m^{-1}))$
uniformly with all derivatives on compact sets.  With
$\mathcal A_m:=\sqrt{2\pi}\,i^m(kc/2)^m/m!$, set
$v_m:=\mathcal A_m\mathcal V_kP_m$, an entire Helmholtz solution, and
let $g_m$ be its Herglotz density.  Define the negative mode by
\[
 v_{-m}(z):=(-1)^m\overline{v_m(z)},
 \qquad
 g_{-m}(t):=(-1)^m\overline{g_m(t+\pi)}.
\]
Expanding $P_m(z)=\sum_{j\le m}p_{mj}z^j$ with $p_{mm}=c^{-m}$, the
Jacobi--Anger formula converts each monomial into a multiple of
$\psi_j$, so $g_{m}$ is supported on the modes $\psi_{j}$ with
$0\le j\le m$, and likewise for $g_{-m}$ on $\psi_{-j}$.  The degree-$m$
coefficient is exactly one by the choice of $\mathcal A_m$, and the
degree-$(m-1)$ coefficient vanishes because the translation removed
$c_0$.  The coefficients $p_{mj}$ themselves carry powers of $m$.  For
the ellipse $\Psi(w)=c(w+qw^{-1})$ one has
$P_m(\Psi(w))=w^m+q^mw^{-m}$ and $p_{m,m-2}=-mq\,c^{-(m-2)}$.  The
coefficient estimate is obtained from the inverse map.  Let
$\Phi:=\Psi^{-1}$ and let $\widetilde Q_m$ be the nonnegative Laurent
part of $\Phi^m$ at infinity.  Then $\widetilde Q_m(\Psi(w))-w^m$ is
holomorphic in $|w|>1$ and vanishes at $w=\infty$, so $\widetilde Q_m$
obeys the normalization above.  That normalization determines the
polynomial uniquely, a polynomial $R$ with $R(\Psi(w))=O(w^{-1})$ being
constant and hence zero.  So $P_m=\widetilde Q_m$ exactly.  Also
$\Phi(z)=(z/c)\bigl(1+\xi_2(c/z)^2+\xi_3(c/z)^3+\cdots\bigr)$, the
series converging for $|z|$ beyond the singularities of $\Phi$, so that
$|\xi_j|\le A_\Phi B_\Phi^{\,j}$ for some $A_\Phi,B_\Phi$.  No bound
$B_\Phi<1$ is claimed or needed, and none holds in general.  For the
ellipse $\Psi(w)=c(w+qw^{-1})$ one has $\xi_{2p}=-C_{p-1}q^{p}$ with
$C_{p-1}$ the Catalan numbers and $\xi_j=0$ for odd $j$, so
$\limsup_j|\xi_j|^{1/j}=2\sqrt q$, the focal distance in units of
$c$.  Also $\xi_1=0$, that is $\Phi$ has no
constant term, again because $c_0=0$.  Expanding $\Phi^m$, the coefficient of $z^{m-\ell}$ collects
partitions of $\ell$ into parts $\ge2$, at most $\lfloor\ell/2\rfloor$
of which can carry the multiplicity $m$:
\[
 |p_{m,m-\ell}|\,c^{\,m-\ell}
 \le C^{\ell}\,\frac{m^{\lfloor\ell/2\rfloor}}{\lfloor\ell/2\rfloor!}.
\]
For the ellipse, $\ell=2$ gives the factor $mq$ above.  In
the Herglotz coefficient this power of $m$ is offset by the factorial
ratio:
\[
 \Bigl(\frac{kc}2\Bigr)^{\ell}\,\frac{(m-\ell)!}{m!}\,
 C^{\ell}\,\frac{m^{\lfloor\ell/2\rfloor}}{\lfloor\ell/2\rfloor!}
 \le\frac{C'^{\ell}}{\lfloor\ell/2\rfloor!}\,
 m^{-\lceil\ell/2\rceil},
 \qquad \ell\le\frac m2,
\]
which is $O(m^{-1})$ at $\ell=2$ and decays at least geometrically in
$\ell$ for $m\ge M_0$.  On the range $m/2<\ell\le m$ enlarge $C$ so that
$kcC\ge1$.  In $m!/(m-\ell)!=m(m-1)\cdots(m-\ell+1)$ the first
$\lceil\ell/2\rceil$ factors are each at least
$m-\lceil\ell/2\rceil+1\ge m/2$ because $\ell\le m$, and the remaining
ones are each at least $m-\ell+1\ge1$, so
$m!/(m-\ell)!\ge(m/2)^{\lceil\ell/2\rceil}$.  With this,
$\lfloor\ell/2\rfloor\le\lceil\ell/2\rceil$,
$2^{\lceil\ell/2\rceil}\le\sqrt{2}\,2^{\ell/2}$,
$\ell\le2\lfloor\ell/2\rfloor+2$ and the bound
$m<4(\lfloor\ell/2\rfloor+1)$ forced by $\ell>m/2$,
\[
\begin{aligned}
&m2^{\ell/2}\Bigl(\frac{kc}2\Bigr)^{\ell}\frac{(m-\ell)!}{m!}
 C^{\ell}\frac{m^{\lfloor\ell/2\rfloor}}{\lfloor\ell/2\rfloor!}\\
&\qquad\le\frac{\sqrt{2}\,m\,(kcC)^{\ell}}{\lfloor\ell/2\rfloor!}\\
&\qquad\le
 \frac{4\sqrt{2}\,(\lfloor\ell/2\rfloor+1)
 (kcC)^{2\lfloor\ell/2\rfloor+2}}{\lfloor\ell/2\rfloor!},
 \qquad \frac m2<\ell\le m,
\end{aligned}
\]
whose right-hand side depends on $\lfloor\ell/2\rfloor$ alone and tends
to $0$, hence is bounded by some $C_1$.  On both ranges the bound is
therefore $C_1m^{-1}2^{-\ell/2}$.  This proves the \emph{two-index} estimate
\begin{equation}\label{eq:gmode-entry}
 \bigl|\ip{g_{\pm m}}{\psi_{\pm j}}-\delta_{mj}\bigr|
 \le\frac Cm\,2^{-(m-j)/2},
 \qquad 0\le j\le m,\quad m\ge M_0 .
\end{equation}
Summing \eqref{eq:gmode-entry} gives
$\norm{g_{\pm m}-\psi_{\pm m}}\le C/m$, and the row and column
$\ell^1$-norms of the matrix
$(\ip{g_{\pm m}}{\psi_{\pm j}}-\delta_{mj})_{m,j\ge M}$ are $O(1/M)$.
Hence the Gram matrix $\mathcal Gr$ of $(g_{\pm m})_{m\ge M}$ is
$I+O(M^{-1})$ in operator norm, with entries in the class
\eqref{eq:Gram-class}.  \Cref{lem:one-sided-cholesky} supplies a
block-lower triangular $\mathcal Ch=I+Z$ with
$\mathcal Ch^{*}\mathcal Gr\mathcal Ch=I$, so the vectors
$g_{\pm m}\mathcal Ch$ are orthonormal and define the unitary $U$.
Because $Z_{nm}=0$ for $n<m$ and $d_n\le d_m$ for $n\ge m$, the
weighted correction
$\mathcal Ch_c:=D_c^{1/2}\mathcal ChD_c^{-1/2}=I+Z_c$ obeys
\begin{equation}\label{eq:weighted-correction}
 (Z_c)_{nm}=\sqrt{d_n/d_m}\;Z_{nm},
 \qquad
 \opnorm{(Z_c)_{nm}}\le\opnorm{Z_{nm}}
 \qquad(n\ge m),
\end{equation}
so replacing $g_{\pm m}$ by $U^{*}\psi_{\pm m}$ perturbs all estimates
below by $O(M^{-1})$, entrywise with the same geometric off-diagonal
decay.

The orthonormal tail vectors define an \emph{isometry}, and extending
an isometry to a unitary is not by itself enough.  The complementary
columns must also be analytic, and an arbitrary extension need not be.
They can be chosen so.  Write $\mathcal L$ for the synthesis operator
$\mathcal Le_{\pm m}:=g_{\pm m}$, $m\ge M_0$.  Since $\mathcal L$ is a
compact perturbation of the inclusion of the Fourier tail, it has the
same Fredholm index, so $(\ran\mathcal L)^{\perp}$ and the fixed
Fourier complement have equal (finite) dimension, and
$U^{*}|_{\rm tail}=\mathcal L\mathcal Ch$ extends to a unitary $U^{*}$
by any isometry between them.  If $h\perp\ran\mathcal L$, then
orthogonality to each $g_m$ gives the triangular recursion
\[
 \widehat h(m)
 =-\sum_{j<m}\overline{\ip{g_m-\psi_m}{\psi_j}}\,\widehat h(j),
\]
and \eqref{eq:gmode-entry} with induction in $m$, followed by the
geometric convolution of \cref{lem:classes}, gives
$|\widehat h(m)|\le C_h\rho_1^{m}$.  The negative frequencies are
identical.  As $(\ran\mathcal L)^{\perp}$ is finite dimensional the
constants may be taken uniform on an orthonormal basis of it.  The
resulting analytic fixed complement is used later, in
\cref{prop:trace-localization} and \cref{lem:S-coordinate}; neither is
needed here.
\end{proof}

The classical Vekua and Faber references do not supply the trace-localization
estimate \eqref{eq:rmode-profile} needed below.  The degree-$m$ polynomial is composed with a
\emph{noncircular} analytic parametrization, so $z(t)^j$ is a
degree-dependent analytic multiplier whose strip norm may grow
exponentially in $j$.  Applying the monomial estimate term by term and
then invoking ``analytic multiplication preserves bandedness''
discards the cancellations that keep \eqref{eq:grunsky}
sharp.  We prove it instead by keeping $P_m$ intact, at the cost of
the geometric hypothesis of \cref{def:centered}.

\begin{lemma}[The exact Faber primitive]
\label{lem:Faber-primitive}
Write
\begin{equation}\label{eq:Psi-prime-coeff}
 \Psi'(w)=\sum_{\ell\ge0}\beta_\ell w^{-\ell},
 \qquad \beta_0=c,\qquad \beta_1=0.
\end{equation}
There are $C>0$ and $\vartheta_0<1$ such that
\begin{equation}\label{eq:three-analytic-estimates}
 |\beta_\ell|\le C\vartheta_0^\ell,
 \qquad
 |b_{ns}|\le C\vartheta_0^{n+s},
 \qquad
 |P_n(0)|\le C\vartheta_0^n.
\end{equation}
Let $\mathcal I Q(z)=\int_0^zQ(\zeta)\,d\zeta$ and let
$\mathcal JQ$ be the unique primitive of $Q$ whose $P_0$-coefficient
is zero.  Then the following identities hold:
\begin{equation}\label{eq:exact-Faber-primitive}
 \mathcal JP_m
 =\sum_{\ell=0}^{m}
   \frac{\beta_\ell}{m+1-\ell}P_{m+1-\ell},
 \qquad
 \mathcal IP_m=\mathcal JP_m-(\mathcal JP_m)(0)P_0,
\end{equation}
and
\begin{equation}\label{eq:primitive-constant-decay}
 |(\mathcal JP_m)(0)|\le C\vartheta_1^m
\end{equation}
for some $\vartheta_1<1$.
\end{lemma}

\begin{proof}
The first estimate in \eqref{eq:three-analytic-estimates} follows by
Cauchy's estimate from the continuation of $\Psi'$ to $|w|>\rho$.
The second is the two-variable Cauchy estimate obtained from the
Grunsky generating function used in \eqref{eq:grunsky}.  Polynomial
factors in $n$ and $s$ are absorbed by increasing $\vartheta_0$ while
keeping it below one.

For the third estimate use the Faber generating function
\begin{equation}\label{eq:Faber-generating-function}
 \frac{\Psi'(w)}{\Psi(w)-z}
 =\sum_{n=0}^{\infty}P_n(z)w^{-n-1}.
\end{equation}
Because $0\in D$, the equation $\Psi(w)=0$ has no solution with
$|w|\ge1$.  After increasing the inner radius of the continuation
annulus, if necessary, the left side of
\eqref{eq:Faber-generating-function} at $z=0$ is holomorphic for
$|w|>r_0$ with some $r_0<1$.  Cauchy's estimate on a circle
$r_0<|w|<1$ gives $|P_n(0)|\le C\vartheta_0^n$ after enlarging
$\vartheta_0$, if necessary.

It remains to prove the identity.  The nonnegative Laurent part of
$P_n'(\Psi(w))$ is
\[
 \left[\frac{n w^{n-1}}{\Psi'(w)}\right]_+ .
\]
Differentiating
$P_n(\Psi(w))=w^n+\sum_{s\ge1}b_{ns}w^{-s}$ produces a negative
Laurent tail, and multiplication by $1/\Psi'(w)$, which has only
nonpositive powers of $w$, cannot turn that tail into a nonnegative
power.  Differentiating the right side of
\eqref{eq:exact-Faber-primitive} and taking its nonnegative Laurent
part therefore gives
\[
 \left[
  \frac{\sum_{\ell=0}^{m}\beta_\ell w^{m-\ell}}
       {\Psi'(w)}
 \right]_+
 =\left[
  \frac{w^m\Psi'(w)-
        \sum_{\ell>m}\beta_\ell w^{m-\ell}}
       {\Psi'(w)}
 \right]_+=w^m.
\]
The omitted numerator and its product with $1/\Psi'$ contain only
negative powers.  Since Faber coefficients are precisely the
nonnegative Laurent coefficients after composition with $\Psi$, the
derivative is $P_m$.  This proves the first identity.  Subtracting its
value at zero gives the zero-based primitive.  Finally,
\eqref{eq:primitive-constant-decay} follows from the first and third
estimates in \eqref{eq:three-analytic-estimates}, the geometric
convolution estimate, and absorption of the denominators into the
constant.
\end{proof}

For $q\ge1$ set
\begin{equation}\label{eq:qfold-primitive}
 \mathcal I^qQ(z)
 :=\frac1{(q-1)!}\int_0^z(z-\zeta)^{q-1}Q(\zeta)\,d\zeta,
 \qquad \mathcal I^0Q:=Q,
\end{equation}
and introduce the reflected exterior map
\begin{equation}\label{eq:reflected-map}
 \Psi^\#(w):=\overline{\Psi(1/\overline w)}.
\end{equation}
On $|w|=1$, $\Psi^\#(w)=\overline{\Psi(w)}$, while
$\chi(w):=w\Psi^\#(w)$ is holomorphic in a disc $|w|<R_\chi$ for some
$R_\chi>1$.

\begin{lemma}[Uniform iterated-primitive Laurent estimate]
\label{lem:iterated-primitives}
For every fixed $k>0$ there are $C>0$ and $\vartheta<1$, independent
of $m$, such that, on $|w|=1$,
\begin{align}
 \mathcal R_m(w)
 &:={\sum_{q\ge1}\frac{(-k^2/4)^q}{q!}
  (\Psi^\#(w))^q(\mathcal I^qP_m)(\Psi(w))},
 \label{eq:Vekua-remainder-series}\\
 |\widehat{\mathcal R_m}(j)|
 &\le \frac C m\vartheta^{|m-j|}
       +C\vartheta^{m+|j|},
 \qquad j\in\mathbb Z.\label{eq:Vekua-function-profile}
\end{align}
If $\mathcal N$ denotes exterior normal differentiation in
the conformal parameter, namely
\begin{equation}\label{eq:conformal-normal-qterm}
 \begin{split}
 \mathcal N\!\left[(\Psi^\#)^q
       (\mathcal I^qP_m)\circ\Psi\right]
 ={}&(\Psi^\#)^q(w\Psi')
       (\mathcal I^{q-1}P_m)\circ\Psi\\
 &-q(\Psi^\#)^{q-1}(w(\Psi^\#)')
       (\mathcal I^qP_m)\circ\Psi,
 \end{split}
\end{equation}
then
\begin{equation}\label{eq:Vekua-normal-profile}
 \left|
  \widehat{\sum_{q\ge1}\frac{(-k^2/4)^q}{q!}
   \mathcal N\!\left[(\Psi^\#)^q
    (\mathcal I^qP_m)\circ\Psi\right]}(j)
 \right|
 \le C\vartheta^{|m-j|}+C\vartheta^{m+|j|}.
\end{equation}
The estimates remain valid after multiplication by any fixed function
holomorphic in a common annulus about $|w|=1$.
\end{lemma}

\begin{proof}
We give the coefficient argument, including the low-mode reset where
centeredness is used.  In the Faber basis write
\begin{equation}\label{eq:primitive-matrix}
 \mathcal IP_m=\sum_{n=0}^{m+1}A_{nm}P_n.
\end{equation}
Lemma~\ref{lem:Faber-primitive} gives, for $n\ge1$,
\begin{equation}\label{eq:primitive-matrix-entry}
 A_{nm}=\frac{\beta_{m+1-n}}{n}
 \quad(1\le n\le m+1),
 \qquad
 |A_{0m}|\le C\vartheta_1^m.
\end{equation}
Also $A_{10}=c$ and $A_{n0}=0$ for $n\ne1$, because
$\mathcal IP_0=z=cP_1$ when $c_0=0$.

An entry of $A^q$ is a sum over paths
$\nu_0=m,\nu_1,\ldots,\nu_q=n$.  As long as a path avoids zero, put
\begin{equation}\label{eq:path-loss}
 \ell_r:=\nu_{r-1}+1-\nu_r\ge0.
\end{equation}
Its weight is
\begin{equation}\label{eq:path-weight}
 \prod_{r=1}^q\frac{\beta_{\ell_r}}{\nu_r},
 \qquad
 \sum_{r=1}^q\ell_r=m+q-n.
\end{equation}
Choose numbers
$\vartheta_0<\eta_0<\eta_1<1$.  Since
$\sum_{\ell\ge0}|\beta_\ell|\eta_0^{-\ell}<\infty$, summing the
numerators of all paths with total loss $L$ gives at most
$C^q\eta_0^L$.  If every $\nu_r\ge m/2$, the denominator product in
\eqref{eq:path-weight} supplies a factor at most $2/m$.  If some
$\nu_r<m/2$, then the cumulative loss up to that step is at least
$m/2$, and
$\eta_0^L\le C m^{-1}\eta_1^L$.  Thus the part of $A^q_{nm}$ from
paths avoiding zero obeys the uniform estimate
\begin{equation}\label{eq:no-reset-path-bound}
 |(A^q_{nm})_{\mathrm{nr}}|
 \le\frac{C^q}{m}\eta_1^{m+q-n},
 \qquad 1\le n\le m+q,\quad q\ge1.
\end{equation}

Now consider a path which visits zero, and stop it at its first visit.
The no-reset part before that visit has accumulated loss
$m+r-1-\nu_{r-1}$, while the last entry into zero contributes
$|A_{0,\nu_{r-1}}|\le C\vartheta_1^{\nu_{r-1}}$.  Enlarging the
geometric base but keeping it below one shows that the initial
segment is bounded by $C^r\eta_1^m$.  After the visit to zero the
index can increase by at most one at each step.  The matrix in
\eqref{eq:primitive-matrix-entry} is bounded on every fixed
exponentially weighted finite-sequence norm.  Explicitly, for a fixed
$S>1$ sufficiently close to one,
\[
 \sum_{n\ge0}S^n|A_{np}|\le C S^p.
\]
In the nonconstant rows $n=p+1-\ell$, and
$\sum_{\ell\ge0}|\beta_\ell|S^{-\ell}<\infty$, while the zero row is
controlled by \eqref{eq:primitive-constant-decay}.  Hence the
remaining segment contributes at most $C^{q-r}$.  Consequently
\begin{equation}\label{eq:reset-path-bound}
 (A^q_{nm})_{\mathrm r}=0\quad(n>q),
 \qquad
 |(A^q_{nm})_{\mathrm r}|\le C^q\eta_1^m
 \quad(0\le n\le q).
\end{equation}
Equations \eqref{eq:no-reset-path-bound} and
\eqref{eq:reset-path-bound} are the required uniform
iterated-Faber-primitive estimate.

We now convert it into Laurent estimates without discarding the Faber
cancellations.  For $1<R<R_\chi$, write
\begin{equation}\label{eq:chi-power}
 \chi(w)^q=\sum_{r\ge0}h_{q r}w^r,
 \qquad |h_{q r}|\le M_R^qR^{-r}.
\end{equation}
Since $(\Psi^\#)^q=w^{-q}\chi^q$, the positive Faber part of a
no-reset term has frequency
\begin{equation}\label{eq:main-frequency-relation}
 j=n-q+r=m-L+r,
 \qquad L=m+q-n.
\end{equation}
The convolution of $\eta_1^L$ and $R^{-r}$ therefore gives
$C^qm^{-1}\eta_2^{|m-j|}$ for any
$\max(\eta_1,R^{-1})<\eta_2<1$.  For the negative Faber tail the two bases must be ordered before they
are combined.  Enlarging the Grunsky base once, write
$|b_{ns}|\le C\widetilde\vartheta^{\,n+s}$ with
$\widetilde\vartheta\in(\eta_1,1)$.  This is legitimate since
$\vartheta_0<\eta_1$ and a larger base is a weaker bound.  The
frequency is $j=-s-q+r$, and substituting $n=m+q-L$ gives
\[
 \eta_1^{L}\widetilde\vartheta^{\,n+s}
 =\widetilde\vartheta^{\,m+q+s}
  \Bigl(\frac{\eta_1}{\widetilde\vartheta}\Bigr)^{L}
 \le\widetilde\vartheta^{\,m+q+s},
\]
because $\eta_1<\widetilde\vartheta$.  If $j\ge0$, \eqref{eq:chi-power} supplies
$R^{-j}$.  If $j<0$, the relation $s+q-r=|j|$ supplies geometric
decay in $|j|$.  After one enlargement of the base this part is
$C^qm^{-1}\eta_2^{m+|j|}$.

For a reset term $n\le q$.  Its positive Faber frequency is
$j=n-q+r$.  Positive $j$ is controlled by $R^{-j}$; negative $j$
forces $q\ge|j|$.  The factor $1/q!$ in
\eqref{eq:Vekua-remainder-series} then gives decay faster than every
fixed geometric sequence in $|j|$.  For the negative Faber tail the frequency is $j=-s-q+r$.  If
$j\ge0$, then $r=j+s+q\ge j$, so \eqref{eq:chi-power} supplies
$R^{-r}\le R^{-j}$, while $\sum_{s\ge1}|b_{ns}|\le C$ uniformly in
$n$.  With the factor $\eta_1^m$ of \eqref{eq:reset-path-bound} this
part is $C^q\eta_2^{m+|j|}$ after one enlargement of the base.  If
$j<0$, then $s+q=|j|+r\ge|j|$, so either $s\ge|j|/2$, in which case
$b_{ns}$ gives the decay, or $q\ge|j|/2$, in which case $1/q!$ gives
it.  Summing
$((k^2/4)C)^q/q!$ proves \eqref{eq:Vekua-function-profile}.

Formula \eqref{eq:conformal-normal-qterm} follows from the chain rule.
Outward conformal differentiation acts as $w\partial_w$ on the
holomorphic factor and as $-w\partial_w$ on its reflected factor.
In its first line, the case $q=1$ contains $P_m$ and is bounded by the
same convolution argument without the factor $m^{-1}$.  For $q\ge2$,
\eqref{eq:no-reset-path-bound} applies to $\mathcal I^{q-1}P_m$.
The second line contains $q\mathcal I^qP_m$; the factor $q$ is absorbed
by $q!$.  The reset terms satisfy the same estimates.  This proves
\eqref{eq:Vekua-normal-profile}.  Finally, multiplication by a fixed
annular analytic function convolves with a two-sided geometric
sequence, so the last assertion follows from
the geometric convolution estimate in the proof of \cref{lem:classes}.
\end{proof}

\begin{lemma}[Pulled-back middle operator and amplitude scale]
\label{lem:middle-pullback}
Let $T_{\rm par}$ be the exact form pull-back of $T$.  Then
\begin{equation}\label{eq:T-principal}
 T_{\rm par}=T_{\rm circ}+R_0,
\end{equation}
where $T_{\rm circ}$ is the circular middle operator and
$R_0:H^{1/2}(\mathbb S^1)\to H^{-1/2}(\mathbb S^1)$ has order zero.  The
inverse factors as $T_{\rm par}^{-1}=(I+Y)T_{\rm circ}^{-1}$ with $Y$
banded of order $-1$, and for all densities on $\Gamma$ the form
transport
\begin{equation}\label{eq:assembly-form}
 \ip{T^{-1}\varphi}{\chi}_\Gamma
 =\ip{T_{\rm par}^{-1}\bigl(|\Psi'|\widetilde\varphi\bigr)}
      {|\Psi'|\widetilde\chi}_{L^2(\dd t)},
 \qquad
 T_{\rm par}^{-1}=(I+Y)\,T_{\rm circ}^{-1},
\end{equation}
holds.  Moreover,
\begin{equation}\label{eq:amp-stirling}
 |\mathcal A_m|^2=\frac{d_m}{C_km}\bigl(1+O(m^{-1})\bigr),
\end{equation}
and
\begin{equation}\label{eq:amp-ratio}
 (8\pi k)^{-1/2}|\mathcal A_m\mathcal A_n|
 \le C\frac{\sqrt{d_md_n}}{\sqrt{mn}} .
\end{equation}
These conclusions do not use trace localization.
\end{lemma}

\begin{proof}
Let $T_{\rm par}$ be the exact form pull-back,
characterized by
\[
 \ip{T\varphi}{\chi}_\Gamma
 =\ip{T_{\rm par}\widetilde\varphi}{\widetilde\chi}_{L^2(\dd t)}.
\]
The principal part of $T$ is the hypersingular operator, whose bilinear
form is conformally invariant to principal order, which gives
\eqref{eq:T-principal} in these coordinates.
This is proved in Appendix~\ref{app:pullback}, where the transported
kernels are decomposed on the analytic curve and $R_0$ is shown to have
a banded Fourier matrix with exponentially decaying off-diagonal
entries.  Equivalently, one may use the half-density conjugate
$|\Psi'|^{1/2}T|\Psi'|^{1/2}$, interpreted after pull-back.  The difference consists of order-zero
commutators and is absorbed into $R_0$.  The transfer to the inverse is
a parametrix argument.  For the circle,
$T_{\rm circ}$ is diagonal in the Fourier basis with symbol
$|t_n|\asymp1+|n|$, computable from Bessel functions and nonvanishing for
positive impedance \cite{CCH13}, so $X:=T_{\rm circ}^{-1}R_0$ has order
$-1$ and a banded matrix (Lemma~\ref{lem:classes}).  The operator $I+X$ is invertible on
$L^2(\mathbb S^1)$.  The identity
$T_{\rm par}=T_{\rm circ}(I+X)$ gives this on
$H^{1/2}(\mathbb S^1)$, and Appendix~\ref{app:pullback} carries it to
$L^2(\mathbb S^1)$ by a Fredholm argument of index zero.  By
Lemma~\ref{lem:classes}(ii) (see also \cite{Jaf90}),
$(I+X)^{-1}=I+Y$ with $Y$ banded and,
by the composition calculus, again of order $-1$.  Hence
$T_{\rm par}^{-1}=(I+Y)T_{\rm circ}^{-1}$
obeys the same weighted block estimates as $T_{\rm circ}^{-1}$, with
banded off-diagonal errors one order lower.
For densities on $\Gamma$, the form transport defining $T_{\rm par}$,
applied to $\psi:=T^{-1}\varphi$, gives
$T_{\rm par}\widetilde\psi=|\Psi'|\widetilde\varphi$, which is
\eqref{eq:assembly-form}.

The amplitude comparison follows from
\[
 \mathcal A_m=\sqrt{2\pi}\,i^m(kc/2)^m/m!,
 \qquad d_m=C_k(ekc/2m)^{2m},
\]
and Stirling's formula:
\[
 |\mathcal A_m|^2
 =\frac{d_m}{C_km}\bigl(1+O(m^{-1})\bigr),
\]
because
$(m/e)^m/m!=(2\pi m)^{-1/2}(1+O(m^{-1}))$.
The estimate \eqref{eq:amp-ratio} follows immediately by applying this
comparison to the two indices $m$ and $n$.
\end{proof}

\begin{proposition}[Faber--Vekua trace localization]
\label{prop:trace-localization}
For the modes $g_{\pm m}$ constructed in
Lemma~\ref{lem:faber-coordinates}, the impedance trace has, in the
conformal parameter, the form
\begin{equation}\label{eq:Bmode}
 Bg_m=\mathcal A_m\left(
   \frac{m}{|\Psi'(e^{it})|}e^{imt}+r_m
 \right).
\end{equation}
After setting $s_m=|\Psi'|r_m$, there are $C>0$ and $\rho_0<1$,
independent of $m$, such that
\begin{equation}\label{eq:rmode-profile}
 |\widehat s_m(j)|
 \le C\rho_0^{|m-j|}+C\rho_0^{m+|j|},
 \qquad j\in\mathbb Z,\quad m\ge M_0.
\end{equation}
The analogous estimate holds for $g_{-m}$ with the signs of the
principal frequencies reversed.  If $h$ belongs to the fixed
complement constructed in Lemma~\ref{lem:faber-coordinates}, then both
orderings of the middle pairing obey
\begin{equation}\label{eq:fixed-complement-form-bound}
 \left|\ip{T^{-1}Bh}{Bg_{\pm m}}_\Gamma\right|
 +\left|\ip{T^{-1}Bg_{\pm m}}{Bh}_\Gamma\right|
 \le C_h|\mathcal A_m|m\rho_0^m.
\end{equation}
In particular, after division by $\sqrt{d_m}$ the fixed-complement
coupling is $O(\rho_1^m)$ for some $\rho_1<1$ in either ordering.
\end{proposition}

\begin{proof}
The exact Vekua series is
\begin{equation}\label{eq:exact-Vekua-primitive-series}
 \mathcal V_kP_m(z)
 =\sum_{q=0}^{\infty}\frac{(-k^2/4)^q}{q!}
  \overline z^{\,q}(\mathcal I^qP_m)(z).
\end{equation}
For a monomial it reduces to
\[
 \mathcal V_k(z^j)=z^j
 \sum_{q=0}^{\infty}
 \frac{(-1)^qj!}{q!(j+q)!}
 \left(\frac{k^2|z|^2}{4}\right)^q,
\]
so \eqref{eq:exact-Vekua-primitive-series} has the
normalization fixed in Lemma~\ref{lem:faber-coordinates} and used in
\eqref{eq:amp-stirling}.  On the boundary its
$q\ge1$ part is \eqref{eq:Vekua-remainder-series}.  The $q=0$ part
satisfies
\[
 P_m(\Psi(w))=w^m+\sum_{s\ge1}b_{ms}w^{-s},
 \qquad |b_{ms}|\le C\vartheta_0^{m+s}.
\]
On the same boundary,
\[
 |\Psi'|\partial_\nu(P_m\circ\Psi)
 =m w^m-\sum_{s\ge1}s b_{ms}w^{-s}.
\]
Combining these identities with
\eqref{eq:Vekua-function-profile} and
\eqref{eq:Vekua-normal-profile} gives the required profile for the
normal derivative.  The impedance contribution is a fixed analytic
multiplier times the boundary value of $\mathcal V_kP_m$; it has order
zero in $m$ and is absorbed in $r_m$.  Its product
with the principal term $w^m$ is geometrically localized about $m$,
and all other terms are covered by the last assertion of
Lemma~\ref{lem:iterated-primitives}.  This proves
\eqref{eq:Bmode}--\eqref{eq:rmode-profile}.

The negative mode defined in Lemma~\ref{lem:faber-coordinates} is the
anti-Vekua continuation associated with $\overline{P_m(z)}$.
The shift by $\pi$ is the change of incident direction needed when a
Herglotz wave is conjugated.  The prefactor makes the coefficient of
$\psi_{-m}$ equal to one.
On the boundary,
\[
 \overline{P_m(\Psi(w))}
 =w^{-m}+\sum_{s\ge1}\overline{b_{ms}}w^s.
\]
Conjugating the path expansion \eqref{eq:path-weight} preserves every
absolute-value bound and sends the frequency relation
$j=m-L+r$ to $j=-m+L-r$.  Thus
\eqref{eq:Vekua-function-profile} and
\eqref{eq:Vekua-normal-profile} hold with $m,j$ replaced by $-m,-j$.
The real analytic multipliers $|\Psi'|$ and $\lambda$ do not change the
estimate.  The factor $i$ in the impedance term changes only its phase.
This proves the negative-frequency profile without any symmetry
assumption on $D$.

For the last assertion, the complementary density $h$ has
exponentially decaying Fourier coefficients by
Lemma~\ref{lem:faber-coordinates}.  Its Herglotz field and impedance trace
are therefore analytic in the conformal parameter.  The inverse
$T^{-1}$ preserves this analytic Fourier class by
\cref{lem:kernels} and the inverse-closed matrix lemma.
Pairing its geometric Fourier profile with the principal mode $m e^{\pm
imt}$ and with \eqref{eq:rmode-profile} gives the first term of
\eqref{eq:fixed-complement-form-bound}.  The second term is the
transposed coupling, and $T^{-1}$ is not self-adjoint, so it does not
follow from the first by conjugation.  Move the inverse onto the
analytic argument instead: by the form transport
\eqref{eq:assembly-form}, which rests on the transport identity of
Appendix~\ref{app:pullback} and on
Lemma~\ref{lem:middle-pullback}, neither of which uses the present proposition,
\[
 \ip{T^{-1}Bg_{\pm m}}{Bh}_\Gamma
 =\ip{|\Psi'|\widetilde{Bg_{\pm m}}}
     {T_{\rm par}^{-*}\bigl(|\Psi'|\widetilde{Bh}\bigr)}_{L^2(\dd t)},
 \qquad
 T_{\rm par}^{-*}=T_{\rm circ}^{-*}(I+Y^{*}),
\]
whose left argument is $\mathcal A_m(me_{\pm m}+s_{\pm m})$ by
\eqref{eq:Bmode}.  The banded classes are closed under adjoints by
\cref{lem:classes}(iii), so $Y^{*}$ again lies in $\mathcal B^{-1}$; and
$T_{\rm circ}^{-*}$ is diagonal with symbol $\overline{t_n}^{-1}$, of
the same modulus as $t_n^{-1}$.  Hence $T_{\rm par}^{-*}$ obeys the
weighted banded estimates of $T_{\rm par}^{-1}$ and produces the same
geometric Fourier profile from $|\Psi'|\widetilde{Bh}$.  Pairing that
profile with $me_{\pm m}$ and with
\eqref{eq:rmode-profile} as before gives the second term with
the same factor $|\mathcal A_m|m$ and the same geometric rate, only
the constant differing; $\rho_0$ and $C_h$ are enlarged once to serve
both orderings.  Finally
$|\mathcal A_m|m/\sqrt{d_m}=O(\sqrt m)$ by
\eqref{eq:amp-ratio}, which follows from the Stirling comparison
\eqref{eq:amp-stirling}.  Both \eqref{eq:amp-ratio} and
\eqref{eq:assembly-form} were established in
Lemma~\ref{lem:middle-pullback}, independently of the present proposition.  The polynomial
factor is absorbed by a slightly larger geometric base.
\end{proof}

\begin{lemma}[Faber--Vekua normal form]\label{lem:faber-vekua}
Let $D$ be centered.
There are a unitary $U$ on $L^2(\mathbb S^1)$, constants $C$,
$\rho_1<1$, and $M_0$, and a bounded block matrix
$E=(E_{mn})_{m,n\ge M_0}$ such that, with
$\widetilde F=e^{-i\pi/4}F$ and $D_c|_{\mathcal H_m}=d_mI_2$,
\begin{equation}\label{eq:normal-form}
 U\widetilde FU^{*}
 =D_c^{1/2}(I+E)D_c^{1/2}
 \quad\text{on }\bigoplus_{m\ge M_0}\mathcal H_m,
 \qquad
 \opnorm{E_{mn}}\le C\,\frac{\rho_1^{|m-n|}}{\min(m,n)} .
\end{equation}
By the Schur test, the entrywise bound implies that
$\opnorm{Q_MEQ_M}$ and $\opnorm{Q_ME(I-Q_M)}$ are at most $C'/M$ for
$M\ge M_0$, and the complementary finite block couples to the tail in
either order,
\[
 \opnorm{(I-Q_{M_0})U\widetilde FU^{*}Q_MD_c^{-1/2}}
 +\opnorm{D_c^{-1/2}Q_MU\widetilde FU^{*}(I-Q_{M_0})}
 \le C\rho_1^{M},
\]
where $I-Q_{M_0}$ projects onto the fixed finite complement.  The
second bound does not follow from the first, $\widetilde F$ being
non-self-adjoint, and the symmetrization used in
\cref{thm:capacity-connected} needs both.
\end{lemma}

\begin{proof}
Use the unitary and orthonormal tail constructed in
\cref{lem:faber-coordinates}. First consider the fixed complement.
By \eqref{eq:sandwich} and
\eqref{eq:fixed-complement-form-bound}, both matrix orderings between a
fixed complementary vector and a tail vector are bounded by
\[
 C_h|\mathcal A_m|m\rho_0^m .
\]
After division by $\sqrt{d_m}$, \eqref{eq:amp-stirling} turns the
polynomial factor into $O(\sqrt m)$, which is absorbed by enlarging
the geometric base. Square summation over $m\ge M$, followed by the
weighted triangular correction \eqref{eq:weighted-correction}, gives
\[
 \opnorm{(I-Q_{M_0})U\widetilde FU^*Q_MD_c^{-1/2}}
 +\opnorm{D_c^{-1/2}Q_MU\widetilde FU^*(I-Q_{M_0})}
 \le C\rho_1^M .
\]

The measure factors $|\Psi'|$ cancel the trace factors $|\Psi'|^{-1}$
of \eqref{eq:Bmode} on both sides: with $e_n:=e^{int}$ and
$s_n:=|\Psi'|\,r_n$ as in \eqref{eq:rmode-profile}, and with
$\rho_2\in[\rho_0,1)$ a common base for \eqref{eq:rmode-profile} and
the banded classes of \cref{lem:classes}, the
factorization \eqref{eq:sandwich}, inserted into
\eqref{eq:assembly-form}, becomes entry by entry
\begin{equation}\label{eq:assembly}
\begin{aligned}
 \ip{\widetilde Fg_n}{g_m}
 =-\frac{\overline{\mathcal A_m}\mathcal A_n}{\sqrt{8\pi k}}\Bigl[\,
 &mn\,\ip{(I+Y)T_{\rm circ}^{-1}e_n}{e_m}
 +n\,\ip{(I+Y)T_{\rm circ}^{-1}e_n}{s_m}\\[-1pt]
 &+m\,\ip{(I+Y)T_{\rm circ}^{-1}s_n}{e_m}
 +\ip{(I+Y)T_{\rm circ}^{-1}s_n}{s_m}\Bigr].
\end{aligned}
\end{equation}
Write
\[
 \theta_m:=-\frac{2\pi}{\sqrt{8\pi k}}\,\overline{\mathcal A_m}\mathcal A_m\,m^2t_m^{-1}
\]
for the diagonal of the first term.  The factor $2\pi$ is
$\ip{e_m}{e_m}_{L^2(\dd t)}$, the modes $e_n=e^{int}$ being
unnormalized.  It depends on $\Psi$ only through $c$, so it may be
identified by specializing to $\Psi(w)=cw$.  On the circle the
remainder $s_m$ does \emph{not} vanish.  It carries the impedance and
the lower Bessel-derivative terms, so the specialization gives not an
identity but the comparison
\begin{equation}\label{eq:theta-circle}
 \widetilde\kappa_m=\theta_m+O(d_m/m),
\end{equation}
with $\widetilde\kappa_m$ the circular eigenvalue of
\eqref{eq:disc-kappa}.  By \eqref{eq:amp-stirling} the circular remainder is $O(d_m/m)$, which is
\eqref{eq:theta-circle}.  Lemma~\ref{lem:gamma-ratio} then gives
$\theta_m=d_m\bigl(1+O(m^{-1})\bigr)$, including the constant, with the
$O(m^{-1})$ term absorbed into the diagonal of $E$.  The estimate
\eqref{eq:amp-ratio} is available from
Lemma~\ref{lem:middle-pullback}.
The four terms of \eqref{eq:assembly} are bounded as follows, uniformly
for $m,n\ge M_0$.
\emph{(a) Main--main.}  $T_{\rm circ}^{-1}e_n=t_n^{-1}e_n$, and $Y$
carries its weight on the row index by Lemma~\ref{lem:classes}(ii), $X$
doing so by Appendix~\ref{app:pullback}, so
$|Y_{mn}|\le C\rho_2^{|m-n|}(1+m)^{-1}$ and the off-diagonal part is at
most
$(8\pi k)^{-1/2}|\mathcal A_m\mathcal A_n|\,mn\,|t_n|^{-1}|Y_{mn}|
\le C\sqrt{d_md_n}\,\rho_2^{|m-n|}/\sqrt{mn}$ by
\eqref{eq:amp-ratio}.
\emph{(b) Main--remainder.}  Since $|n\,t_n^{-1}|\le C$, the two cross
terms are bounded by
$C\bigl(|\widehat{s_m}(n)|+|\widehat{s_n}(m)|\bigr)$ plus the
$Y$-smeared sums
\[
 \frac n{|t_n|}\sum_j|\widehat{s_m}(j)|\,|Y_{jn}|
 \qquad\text{and}\qquad
 m\sum_j\frac{|\widehat{s_n}(j)|\,|Y_{mj}|}{|t_j|}\,.
\]
The two are not symmetric.  In the first, $T_{\rm circ}^{-1}$ acts on
$e_n$ \emph{before} $Y$, so the symbol is $t_n^{-1}$ and comes out of
the sum; in the second it acts on $s_n$, so the symbol is $t_j^{-1}$
and stays inside.  Both are
$\le C\rho_3^{|m-n|}$ by \eqref{eq:rmode-profile} and the geometric
convolution of Lemma~\ref{lem:classes}.  With \eqref{eq:amp-ratio} this
gives again $C\sqrt{d_md_n}\,\rho_3^{|m-n|}/\sqrt{mn}$.  At $m=n$ this
term is $O(d_m/m)$.  The impedance contribution sits here, one order
below the diagonal.
\emph{(c) Remainder--remainder.}  The identity part of $I+Y$ gives
$\sum_j|t_j|^{-1}|\widehat{s_n}(j)|\,|\widehat{s_m}(j)|
\le C\rho_3^{|m-n|}/\min(m,n)$, the factor $|t_j|^{-1}$ being
$O(1/\min(m,n))$ on the dominant modes and the small-$j$ terms
exponentially suppressed.  The $Y$ part is
$\sum_{i,j}|t_j|^{-1}|\widehat{s_n}(j)|\,|Y_{ij}|\,
|\widehat{s_m}(i)|$, and since $Y$ carries its weight on the row
index, $|Y_{ij}|\le C\rho_2^{|i-j|}(1+i)^{-1}$, two applications of
the geometric convolution of \cref{lem:classes} bound it by the same
$C\rho_3^{|m-n|}/\min(m,n)$ with one further enlargement of the base.
It is smaller by the extra factor $(1+i)^{-1}$.  Together these give
the smaller bound
$C\sqrt{d_md_n}\,\rho_3^{|m-n|}/\bigl(\sqrt{mn}\,\min(m,n)\bigr)$.
\emph{(d) Analytic tails.}  The Faber remainders of \eqref{eq:grunsky}
and all pairings between positive and negative frequencies (where
$T_{\rm circ}^{-1}$ is diagonal, $|Y_{-m,n}|\le C\rho_2^{m+n}$, and
$|\widehat{s_m}(-j)|\le C\rho_2^{m+j}$) contribute $O(\rho_2^{m+n})$.
The negative-frequency pair satisfies the estimates of (a)--(c)
verbatim, and the $2\times2$ blocks $E_{mn}$ collect the four sign
combinations.
Dividing by $\sqrt{d_md_n}$ and adding the weighted one-sided correction \eqref{eq:weighted-correction},
which perturbs entries within the same class,
every entry of $E$ obeys
$\norm{E_{mn}}_{\rm op}\le C\rho_1^{|m-n|}/\sqrt{mn}
\le C\rho_1^{|m-n|}/\min(m,n)$, which is the entrywise bound of
\eqref{eq:normal-form}.  The diagonal case $m=n$ gives the
required $O(1/m)$.  For the Schur test, the row
sum of $\rho_1^{|m-n|}/\min(m,n)$ over $n\ge m/2$ is $O(m^{-1})$ and the
sum over $n<m/2$ is $O(\rho_1^{m/2}\log m)$.  Columns are symmetric, and
in the mixed block the dominant columns sit near $M$, so the bound is
again $O(M^{-1})$.  Since $-e^{-i\pi/4}\gamma_2=-1/\sqrt{8\pi k}$ is real, no
additional phase enters the diagonal normalization, and
\eqref{eq:normal-form} follows.  The matrix $E$ carries a skew part,
fed by $\Imag T^{-1}$ and the impedance terms, which enter one order
below the hypersingular principal part.  The skew part obeys the same
entrywise bound, and the argument uses entrywise bounds alone.  For the disc this is visible in
Lemma~\ref{lem:gamma-ratio}, where
$\Imag\widetilde\kappa_m/\Real\widetilde\kappa_m=O(m^{-1})$.  This is
consistent with the passivity that makes $\Fh$ positive:
$\Imag\widetilde F\ge0$, strictly so for the disc by the energy bound
of \cref{sec:disc-prelim}, but of relative size $O(m^{-1})$ blockwise.
\end{proof}

\begin{remark}[Obstruction without centeredness]
\label{rem:noncentered-obstruction}
If the conformal center does not lie in $D$, then after the
normalization $c_0=0$ there is $v$ with $|v|\ge1$ and $\Psi(v)=0$;
$|v|>1$ when the center is outside $\overline D$, since $\Psi$ maps
$|w|>1$ bijectively onto that complement, and $|v|=1$ when it lies on
$\Gamma$, since $\Psi$ carries $|w|=1$ homeomorphically onto $\Gamma$.
These two cases exhaust the negation of \cref{def:centered}, and $v$ is
unique with $\Psi'(v)\ne0$, the extension of $\Psi$ to $|w|>\rho$ being
univalent.  Evaluating $P_n(\Psi(w))=w^n+\sum_{s\ge1}b_{ns}w^{-s}$ at
$w=v$---legitimate also at $|v|=1$, since $b_{ns}=O(\rho^{n+s})$ makes
the tail converge on $|w|>\rho$---gives
$P_n(0)=v^n+O(\rho^n)$, while
\eqref{eq:exact-Faber-primitive} at $z=0$, with
$\sum_{\ell\ge0}\beta_\ell v^{-\ell}=\Psi'(v)$ convergent for
$|v|\ge1$, gives
\begin{equation}\label{eq:outside-primitive-growth}
 (\mathcal JP_m)(0)
 =\frac{v^{m+1}\Psi'(v)}{m}
   \bigl(1+O(m^{-1})\bigr).
\end{equation}
Thus the $P_0$ coefficient of the zero-based primitive
$\mathcal IP_m$ grows like $|v|^m/m$ when $|v|>1$ and decays only like
$m^{-1}$ when $|v|=1$, rather than decaying geometrically.  The first
Vekua correction already carries this
coefficient into $\overline{\Psi(w)}$, so the reset estimate
\eqref{eq:reset-path-bound}, and hence the proof of
\eqref{eq:rmode-profile}, fails.  A theorem for a noncentered
obstacle would require a different Helmholtz continuation or an
additional cancellation argument.
\end{remark}

We next isolate the operator-theoretic step that converts the normal
form for $\Real\widetilde F$ into weighted bounds for the
symmetrization.  Passing from $\Real\widetilde F$ to
$|{\Real\widetilde F}|$ flips finitely many negative eigenvalues, and
the flip is small in the weight $D_c$.

\begin{lemma}[Weighted control of the negative part]\label{lem:flip}
Let $A_0=A_0^{*}$ be compact, let $D_c$ be as above, extended by the
identity on the orthogonal complement of
$\bigoplus_{m\ge M_0}\mathcal H_m$, and suppose there are $M_2$ and
$C_0$ such that for every $M\ge M_2$, with $D_{c,t}:=Q_MD_cQ_M$:
{\rm(i)} $Q_MA_0Q_M=D_{c,t}^{1/2}(I+G_M)D_{c,t}^{1/2}$ with
$\opnorm{G_M}\le\tfrac12$;
{\rm(ii)} $\opnorm{D_c^{-1/2}Q_MA_0(I-Q_M)}\le C_0/M$;
{\rm(iii)} $\opnorm{D_c^{-1/2}A_0}\le C_0$.
Then the negative spectral subspace $\mathcal N$ of $A_0$ is
finite-dimensional, $A_0\le-\mu_0<0$ on $\mathcal N$, and for all large
$M$,
\begin{equation}\label{eq:flip-bounds}
\begin{gathered}
 \opnorm{D_c^{-1/2}Q_M\bigl(|A_0|-A_0\bigr)Q_MD_c^{-1/2}}\le\frac C{M^2},\\
 \opnorm{D_c^{-1/2}Q_M\bigl(|A_0|-A_0\bigr)(I-Q_M)D_c^{-1/2}}\le\frac CM .
\end{gathered}
\end{equation}
\end{lemma}

\begin{proof}
By (i), $\ip{A_0v}v\ge\tfrac12\ip{D_cv}v>0$ for
$0\ne v\in\operatorname{ran}Q_{M_2}$, so
$\mathcal N\cap\operatorname{ran}Q_{M_2}=\{0\}$ and
$\dim\mathcal N\le\operatorname{corank}Q_{M_2}<\infty$.  Since $A_0$ is
compact and self-adjoint, its finitely many negative eigenvalues lie in
$[-\opnorm{A_0},-\mu_0]$ for some $\mu_0>0$.  Let $A_0v=-\mu v$ with
$\norm v=1$, $\mu\ge\mu_0$.  By (iii), $v=-\mu^{-1}A_0v$ lies in the
domain of $D_c^{-1/2}$ and
\begin{equation}\label{eq:flip-head}
 \norm{D_c^{-1/2}v}=\mu^{-1}\norm{D_c^{-1/2}A_0v}\le C_0/\mu_0 .
\end{equation}
On the tail, $Q_M(A_0+\mu)v=0$ gives
$Q_Mv=-(Q_MA_0Q_M+\mu)^{-1}Q_MA_0(I-Q_M)v$, and conjugating the
resolvent with (i),
\[
 D_{c,t}^{-1/2}\bigl(Q_MA_0Q_M+\mu\bigr)^{-1}D_{c,t}^{1/2}
 =\mu^{-1}\Bigl(I+\mu^{-1}(I+G_M)D_{c,t}\Bigr)^{-1},
\]
of norm at most $2/\mu_0$ once $\opnorm{D_{c,t}}=d_M\le\mu_0/3$.  With
(ii),
\begin{equation}\label{eq:flip-tail}
 \norm{D_c^{-1/2}Q_Mv}\le\frac2{\mu_0}\,\frac{C_0}M .
\end{equation}
Finally $|A_0|-A_0=2(A_0)_-=2\sum_i\mu_i\,v_i\otimes v_i$ over the
finitely many negative eigenpairs, so each term of the first operator
in \eqref{eq:flip-bounds} is a rank-one product of two factors
\eqref{eq:flip-tail}, and each term of the second combines one factor
\eqref{eq:flip-tail} with one factor \eqref{eq:flip-head}, $D_c$
commuting with $Q_M$.
\end{proof}

The block comparison requires a lower bound for $I+K$ on the heads
$\operatorname{ran}(I-Q_M)$, whose dimension grows with $M$.
Positivity and injectivity of the congruence
$\Fh=D_c^{1/2}(I+K)D_c^{1/2}$ do \emph{not} yield such a bound,
because $D_c^{1/2}$ has dense, nonclosed range.  The infimum over
growing heads could degenerate.  The needed estimate follows from the
quantitative content of the range identity.  The second Douglas factor of
\cref{prop:reduction} bounds $\ip{\Fh u}u$ from below by
$\norm{G^{*}u}_{H^{1/2}(\Gamma)}^2$, which in normalized coordinates
is $\norm{\mathcal Sx}$ for a \emph{single fixed} operator
$\mathcal S$.  Uniformity in $M$ is then a consequence rather than a
hypothesis: $\mathcal S$ and $\opnorm{\mathcal S^{-1}}$ do not depend
on $M$.  That $\mathcal S$ is invertible is a Fredholm argument in
which injectivity of $G$ enters once, at the end.

Equip $H^{1/2}(\Gamma)$ with the equivalent Sobolev--Fourier norm and
let
\[
 J:H^{1/2}(\Gamma)\longrightarrow\mathcal Y:=L^2(\mathbb S^1)
\]
be the resulting unitary map.  More explicitly, using the same
form pullback as in \eqref{eq:T-principal}, we set
\[
 Jf=\Lambda^{1/2}\widetilde f,
 \qquad \Lambda e_j=(1+|j|)e_j,
\]
and define the equivalent $H^{1/2}$ norm by $\norm f_{H^{1/2}}=
\norm{Jf}_{L^2}$.  Any analytic half-density in the form pullback is
included in $\widetilde f$.  Its banded order-zero matrix does not alter
the estimates.  Let
\[
 \mathcal R_J:H^{1/2}(\Gamma)\longrightarrow H^{-1/2}(\Gamma)
\]
be its Riesz isomorphism, characterized by
\[
 \langle\mathcal R_Jf,h\rangle_{-1/2,1/2}
 =(f,h)_{H^{1/2}}.
\]
Throughout this subsection $\mathcal U:=L^2(\mathbb S^1)$ is the mode
space.  Write $\Pi_n$ for the orthogonal
projection of $\mathcal Y$ onto the
$n$-th block $\operatorname{span}\{e_{\pm n}\}$ and
$\Pi_t:=\sum_{n\ge M_0}\Pi_n$ for the projection onto the tail.  These
act on the target and are not to be confused with the Faber
polynomials $P_n$.  On finite Fourier sums define
\begin{equation}\label{eq:S-core}
 \mathcal S_0x:=JG^*U^*D_c^{-1/2}x.
\end{equation}

\begin{lemma}[Coordinate form of the normalized map]
\label{lem:S-coordinate}
After increasing $M_0$, the tail-to-tail block of $\mathcal S_0$ has
the form
\begin{equation}\label{eq:S-coordinate}
 \Pi_n\mathcal S_0|_{\mathcal H_m}
 =\delta_{mn}\mathcal B_m+\mathcal R_{nm},
\end{equation}
where
\begin{equation}\label{eq:B-R-bounds}
 s_{\min}(\mathcal B_m)\ge b_->0,
 \qquad
 \opnorm{\mathcal R_{nm}}
 \le C\frac{\rho_2^{|m-n|}}{\sqrt{\min(m,n)}}.
\end{equation}
For every fixed output block $n<M_0$,
\begin{equation}\label{eq:S-fixed}
 \opnorm{\Pi_n\mathcal S_0|_{\mathcal H_m}}\le C_n\rho_2^m.
\end{equation}
Consequently $\Pi_t\mathcal S_0Q_{M_0}$ extends to a boundedly invertible
operator from $Q_{M_0}\mathcal U$ onto $\Pi_t\mathcal Y$.
\end{lemma}

\begin{proof}
From \cref{prop:reduction}, \eqref{eq:sandwich}, and
\eqref{eq:Bmode}, first define the normalized map in the
nonorthogonal $g$-coordinates,
\[
 \mathcal S_g:=JG^*\mathcal L D_c^{-1/2}.
\]
Its positive-frequency column is, by (F4) of \cref{fact:fm-inputs}
together with \eqref{eq:sandwich},
\begin{equation}\label{eq:S-column}
 \mathcal S_g\psi_m
 =\overline{\gamma_2}\,d_m^{-1/2}
 \Lambda^{1/2}T_{\mathrm{circ}}^{-*}(I+Y^*)
 \bigl(\mathcal A_m(me_m+s_m)\bigr),
\end{equation}
with no orthogonalization remainder.  The principal diagonal coefficient is
\[
 b_m=\overline{\gamma_2}\,
 \frac{\mathcal A_m m(1+m)^{1/2}}
 {\overline{t_m}\sqrt{d_m}}.
\]
Here \eqref{eq:amp-ratio} gives only an upper bound.  The lower bound is
the \emph{two-sided} comparison \eqref{eq:amp-stirling}, that is,
$|\mathcal A_m|/\sqrt{d_m}=(C_km)^{-1/2}\bigl(1+O(m^{-1})\bigr)$.  Together with $|t_m|\asymp m$, also
two-sided, this gives $0<b_-\le|b_m|\le b_+$.  The two signs therefore give uniformly
invertible $2\times2$ principal blocks $\mathcal B_m$.

Put $A=\Lambda^{1/2}T_{\mathrm{circ}}^{-*}(I+Y^*)$.  Its diagonal
part is $A_{nn}e_n$, with $|A_{nn}|=O(n^{-1/2})$, and the banded
calculus gives
\begin{equation}\label{eq:A-coordinates}
 |A_{nj}-\delta_{nj}A_{nn}|
 \le C\rho_2^{|n-j|}(1+n)^{-1/2}(1+j)^{-1}.
\end{equation}
For the $Y^*$ correction to $me_m$, use
\eqref{eq:amp-ratio} and \eqref{eq:A-coordinates}.  Its normalized
$n$-th coefficient is at most
\[
 C m^{-1/2}m\rho_2^{|n-m|}
 (1+n)^{-1/2}(1+m)^{-1}
 \le C\frac{\rho_2^{|n-m|}}{\sqrt{\min(m,n)}}.
\]
For the $s_m$ term, \eqref{eq:rmode-profile},
\eqref{eq:A-coordinates}, and the geometric convolution give
\[
 \sum_j|A_{nj}||\widehat s_m(j)|
 \le C(1+n)^{-1/2}\rho_2^{|n-m|}+C\rho_2^{m+n}.
\]
Multiplication by $|\mathcal A_m|/\sqrt{d_m}\le C/\sqrt m$ is
smaller than the asserted remainder.  Opposite-sign terms are
$O(\rho_2^{m+n})$.  We have therefore proved the diagonal-plus-small
estimate for $\mathcal S_g$.

The actual orthonormal-coordinate map is, by
\eqref{eq:weighted-correction},
\begin{equation}\label{eq:S-triangular-transfer}
 \mathcal S_0
 =JG^*\mathcal L\mathcal Ch D_c^{-1/2}
 =JG^*\mathcal L D_c^{-1/2}\mathcal Ch_c
 =\mathcal S_g\mathcal Ch_c.
\end{equation}
Since $\mathcal Ch_c-I$ satisfies the one-sided order $-1$ estimate
\eqref{eq:weighted-correction}, multiplication by $\mathcal Ch_c$ changes the
principal diagonal by $O(m^{-1})$ and adds a remainder smaller than
\eqref{eq:B-R-bounds}.  Hence
\eqref{eq:S-coordinate}--\eqref{eq:B-R-bounds} hold for
$\mathcal S_0$.  When the output index $n$ is fixed, every term in
$\mathcal S_g$ crosses a gap of size at least $m-O(1)$, while
$\mathcal Ch_c\psi_m$ only contains indices $j\ge m$.  Geometric
summation proves \eqref{eq:S-fixed}.

Let $\mathcal B=\bigoplus_{m\ge M_0}\mathcal B_m$ and
$\mathcal R=(\mathcal R_{nm})$.  The Schur test gives
$\opnorm{\mathcal R}\le C/\sqrt{M_0}$.  Increase $M_0$ so that this is
less than $b_-/2$.  Then
$\mathcal B+\mathcal R$ is invertible by the Neumann series.
\end{proof}

\begin{proposition}[Invertibility of the full normalized map]
\label{prop:S-invertible}
The operator $\mathcal S_0$ extends to a boundedly invertible map
$\mathcal S:\mathcal U\to\mathcal Y$.  In particular,
\begin{equation}\label{eq:S-lower}
 \norm{\mathcal Sx}\ge c_S\norm{x},
 \qquad c_S=\opnorm{\mathcal S^{-1}}^{-1}>0.
\end{equation}
\end{proposition}

\begin{proof}
Decompose domain and target into their fixed finite complements and
their tails.  The tail--tail block is invertible by
Lemma~\ref{lem:S-coordinate}.  The tail-to-finite block is bounded by
the square summability of \eqref{eq:S-fixed}.  Blocks with finite domain
are bounded because their finitely many columns belong to $\mathcal Y$.
Thus $\mathcal S_0$ has a bounded extension
\[
 \mathcal S=
 \begin{pmatrix}S_{hh}&S_{ht}\\S_{th}&S_{tt}\end{pmatrix}.
\]
Bounded Gaussian elimination reduces it to
\[
 \begin{pmatrix}
 S_{hh}-S_{ht}S_{tt}^{-1}S_{th}&0\\0&S_{tt}
 \end{pmatrix}.
\]
The finite Schur complement acts between spaces of the same dimension
$2M_0-1$, so it has index zero; $S_{tt}$ is invertible.  Hence
$\mathcal S$ is Fredholm with $\ind\mathcal S=0$.

It remains to show $\ker\mathcal S^*=\{0\}$.  If
$y\in\ker\mathcal S^*$, then for every finite Fourier sum $x$,
\begin{align*}
 0
 &=\ip{y}{JG^*U^*D_c^{-1/2}x}_{\mathcal Y}\\
 &=\ip{G\mathcal R_JJ^*y}{U^*D_c^{-1/2}x}_{\mathcal U}
 =\ip{UG\mathcal R_JJ^*y}{D_c^{-1/2}x}_{\mathcal U}.
\end{align*}
The vectors $D_c^{-1/2}x$ include all finite Fourier sums and are dense
in $\mathcal U$.  Thus $UG\mathcal R_JJ^*y=0$.  Equation
(F3) of \cref{fact:fm-inputs}, together with injectivity of $U$, $\mathcal R_J$,
and $J^*$, gives $y=0$.

Since $\ind\mathcal S=0$, the equality
$\dim\ker\mathcal S=\dim\ker\mathcal S^*$ gives
$\ker\mathcal S=\{0\}$.  Its range is closed because it is Fredholm and
dense because $\ker\mathcal S^*=\{0\}$.  Therefore it is surjective.
The bounded inverse theorem proves \eqref{eq:S-lower}.
\end{proof}

\begin{lemma}[Uniform head floor]
\label{lem:head-floor}
Let
\begin{equation}\label{eq:K-definition}
 I+K:=D_c^{-1/2}U\Fh U^*D_c^{-1/2}
\end{equation}
denote the bounded extension constructed below.  There is
$\kappa_0>0$, independent of $M$, such that
\begin{equation}\label{eq:head-floor}
 I+K_{hh}^{(M)}\ge\kappa_0I
 \quad\text{on }\ran(I-Q_M).
\end{equation}
One may take
\begin{equation}\label{eq:kappa}
 \kappa_0=\opnorm{C_2}^{-2}\opnorm{\mathcal S^{-1}}^{-2}.
\end{equation}
\end{lemma}

\begin{proof}
For $x\in\ran(I-Q_M)$ put
$u=U^*D_c^{-1/2}x$, which is defined because $x$ has finitely many
nonzero coordinates.  By
\cref{prop:reduction}, $G^*=C_2^*\Fh^{1/2}$, and hence
\begin{align*}
 \ip{(I+K)x}{x}
 &=\ip{\Fh u}{u}
 =\norm{\Fh^{1/2}u}^2\\
 &\ge\opnorm{C_2}^{-2}\norm{G^*u}_{H^{1/2}}^2
 =\opnorm{C_2}^{-2}\norm{\mathcal Sx}_{\mathcal Y}^2\\
 &\ge\opnorm{C_2}^{-2}\opnorm{\mathcal S^{-1}}^{-2}\norm{x}^2.
\end{align*}
This is \eqref{eq:head-floor}--\eqref{eq:kappa}, with a constant that
does not depend on the dimension of the head.
\end{proof}

Throughout write
$\gamma_m^{\rm pair}:=\tfrac12\bigl(\lambda_{2m}(\Fh)
+\lambda_{2m+1}(\Fh)\bigr)$ for the average of the $m$-th eigenvalue
pair, the eigenvalues of $\Fh$ sorted decreasingly, and let $m_\alpha$
be the largest $m$ with $\gamma_m^{\rm pair}\ge\alpha$.  For the disc
these agree with the $\gamma_m$ and the transition index of
\cref{sec:disc-prelim}.

\begin{theorem}[Capacity law for a connected analytic obstacle]
\label{thm:capacity-connected}
Let $D$ be bounded and simply connected with analytic boundary and
analytic impedance $\lambda\ge\lambda_0>0$, and assume that $D$ is
centered (\cref{def:centered}).  Centeredness enters only through
\cref{prop:trace-localization} and holds for every convex analytic
obstacle.  Then
\begin{equation}\label{eq:capacity-law}
 \gamma_m^{\rm pair}
 =C_k\Bigl(\frac{ek\operatorname{cap}(\overline D)}{2m}\Bigr)^{2m}
 \bigl(1+O(m^{-1})\bigr),
 \qquad m\to\infty,
\end{equation}
and consequently the transition index satisfies
\begin{equation}\label{eq:capacity-transition}
 m_\alpha=
 \frac{L_\alpha}
 {2\LambertW\bigl(L_\alpha/(ek\operatorname{cap}(\overline D))\bigr)}
 +O(1),
 \qquad \alpha\downarrow0 .
\end{equation}
\end{theorem}

\begin{proof}
\emph{Step 1: the symmetrization in the normal form.}
Taking real and imaginary parts in \eqref{eq:normal-form},
$A_0:=U(\Real\widetilde F)U^{*}=D_c^{1/2}(I+\Real E)D_c^{1/2}$ and
$U(\Imag\widetilde F)U^{*}=D_c^{1/2}(\Imag E)D_c^{1/2}$, with $\Real E$
and $\Imag E$ inheriting the entrywise bound of \eqref{eq:normal-form},
since $|(E^{*})_{mn}|=|\overline{E_{nm}}|$ obeys the same estimate.  We
verify the hypotheses of Lemma~\ref{lem:flip} for $A_0$.  Hypothesis
(i) holds with $G_M=\Real E_{tt}$, of norm $O(1/M)\le\tfrac12$ by the
Schur clause of Lemma~\ref{lem:faber-vekua}.  For (ii), on the block
range
$D_c^{-1/2}Q_MA_0(I-Q_M)
=\bigl[Q_M\Real E\,(I-Q_M)\bigr]D_c^{1/2}(I-Q_M)$, whose first factor
is $O(1/M)$ by the mixed-block clause (the transposed block obeys the
same Schur estimate), and whose second is bounded.  The coupling to the
finite complement is $O(\rho_1^M)$ by the final clause of
Lemma~\ref{lem:faber-vekua}, applied to each summand of
$2\Real\widetilde F=\widetilde F+\widetilde F^{*}$.  For (iii),
$D_c^{-1/2}A_0=(I+\Real E)D_c^{1/2}$ on the block range, bounded, plus
the same exponentially small complement coupling.
Lemma~\ref{lem:flip} yields \eqref{eq:flip-bounds} for $A_0$, and
conjugation by $U$ commutes with the absolute value, so
$U\Fh U^{*}=|A_0|+U(\Imag\widetilde F)U^{*}$.  Hence
$K:=D_c^{-1/2}\,U\Fh U^{*}D_c^{-1/2}-I$ is self-adjoint with
\begin{equation}\label{eq:K-clauses}
 \opnorm{Q_MKQ_M}\le\frac CM,
 \qquad
 \opnorm{Q_MK(I-Q_M)}\le\frac CM,
 \qquad M\ge M_0,
\end{equation}
and $I+K\ge\kappa_0>0$ on every head, uniformly in $M$, by
Lemma~\ref{lem:head-floor}.

\emph{Step 2: block-diagonal pinch.}
Fix $m\ge2M_0$ and put $M:=\lceil m/2\rceil$.  For $y=y_h+y_t$ with
$y_t=Q_My$, the coupling bound in \eqref{eq:K-clauses} gives
$|\langle Ky,y\rangle-\langle Ky_h,y_h\rangle-\langle Ky_t,y_t\rangle|
\le(C/M)(\norm{y_h}^2+\norm{y_t}^2)$, so
\[
 \Bigl(I+K_{hh}-\tfrac CM\Bigr)\oplus\Bigl(I+K_{tt}-\tfrac CM\Bigr)
 \le I+K\le
 \Bigl(I+K_{hh}+\tfrac CM\Bigr)\oplus\Bigl(I+K_{tt}+\tfrac CM\Bigr),
\]
and conjugation by $D_c^{1/2}$ pinches $U\Fh U^{*}$ between two
block-diagonal operators $A'\le U\Fh U^{*}\le A''$.  Every head
eigenvalue of $A'$ and $A''$ is at least
$\tfrac12\kappa_0d_{M-1}$, by Lemma~\ref{lem:head-floor} together with
$D_c\ge d_{M-1}$ on the head.  On the tail, the two-sided bound
\[
 \lambda_j(S^{*}D_{c,t}S)\in
 \bigl[\lambda_{\min}(S^{*}S)\,\lambda_j(D_{c,t}),\
 \lambda_{\max}(S^{*}S)\,\lambda_j(D_{c,t})\bigr],
 \qquad S=(I+K_{tt}\pm C/M)^{1/2},
\]
places the tail eigenvalues of $A'$ and $A''$ within factors
$1\mp2C/M$ of $\lambda_j(D_{c,t})$.  Since
$d_M/d_{M-1}=O(M^{-2})$, the head values exceed every tail value for
large $m$, so in both comparison operators the eigenvalue of global
index $(2M-1)+j$ is the $j$-th tail eigenvalue.  The head has dimension
$2M-1$, so the block-$m$ pair carries the global indices $2m$ and
$2m+1$ exactly.  Weyl monotonicity
$\lambda_j(A')\le\lambda_j(U\Fh U^{*})\le\lambda_j(A'')$ yields
\[
 \Bigl(1-\frac{4C}m\Bigr)d_m\le
 \lambda_{2m}(\Fh),\ \lambda_{2m+1}(\Fh)
 \le\Bigl(1+\frac{4C}m\Bigr)d_m,
\]
whose average is \eqref{eq:capacity-law}.  The adjacent ratio
$d_{m+1}/d_m=(kc)^2(4m^2)^{-1}(1+O(m^{-1}))$ then lets the proof of
Theorem~\ref{thm:disc-sharp}, with $a$ replaced by
$\operatorname{cap}(\overline D)$, give
\eqref{eq:capacity-transition} verbatim.
\end{proof}

Table~\ref{tab:capacity} reports the pair averages
$\gamma_m^{\rm pair}$ at
$m=4,6,8$ for six obstacles, computed with the
combined-field solver of \cref{sec:num-kite} at $k=\lambda=1$.  Three
shapes share capacity $1.30$ with the disc of the paper: two ellipses,
one of aspect ratio $10\!:\!3$, and a nonconvex four-lobed domain,
the image of the unit circle under
$\psi(w)=1.3\,(w+0.3w^{-1}+0.2w^{-3})$---univalent since
$\sum_{n\ge1}n|c_{-n}|/c=0.3+3\,(0.2)=0.9<1$---whose capacity is exactly $1.3$ while its
radius varies between $0.98$ and $1.95$.  By $m=8$ the equal-capacity
eigenvalues agree with the disc values to between $0.1\%$ and $3.2\%$,
the four-lobed shape to $0.1\%$.  The capacity-$1.55$ pair agrees
internally to $1\%$ and sits above the capacity-$1.30$ family by the
factor $16.4$, against $(1.55/1.30)^{16}=16.7$ predicted by the disc law.
The capacity law predicts that the impedance enters only through the
relative $O(m^{-1})$ correction.  The computations are consistent with
this prediction: replacing $\lambda=1$ by $\lambda=0.2$ or
$\lambda=5$ moves the $m=8$ pairs of the disc and of the four-lobed
shape by less than $17\%$, while the capacity change from $1.30$ to
$1.55$ moves them by the factor $16.4$.

\begin{table}[tbhp]
\centering
\caption{Pair-averaged eigenvalues of $\Fh$ for shapes of equal and
different logarithmic capacity ($k=\lambda=1$, $M=64$ directions,
combined-field solver).}
\label{tab:capacity}
\begin{tabular}{lcccc}
\toprule
shape & $\operatorname{cap}$ & $\gamma_4$ & $\gamma_6$ & $\gamma_8$\\
\midrule
disc $a=1.30$          & $1.30$ & $4.01\times10^{-3}$ & $1.20\times10^{-6}$ & $8.99\times10^{-11}$\\
ellipse $1.6\times1.0$ & $1.30$ & $4.03\times10^{-3}$ & $1.21\times10^{-6}$ & $9.04\times10^{-11}$\\
ellipse $2.0\times0.6$ & $1.30$ & $5.11\times10^{-3}$ & $1.24\times10^{-6}$ & $9.27\times10^{-11}$\\
four-lobed (conformal) & $1.30$ & $4.20\times10^{-3}$ & $1.22\times10^{-6}$ & $8.98\times10^{-11}$\\
disc $a=1.55$          & $1.55$ & $1.51\times10^{-2}$ & $9.61\times10^{-6}$ & $1.48\times10^{-9}$\\
ellipse $2.0\times1.1$ & $1.55$ & $1.91\times10^{-2}$ & $9.73\times10^{-6}$ & $1.49\times10^{-9}$\\
\bottomrule
\end{tabular}
\end{table}

The capacity law proved above is limited to a connected analytic boundary,
while the boundary-growth and reconstruction estimates are established only
for the disc.  Multiple components and sampling asymptotics for general
analytic obstacles are not treated here.

\begin{remark}[Scope of the capacity law]
Both the capacity and the exact constant in
\eqref{eq:capacity-law} come from the Faber--Vekua analysis.
Proposition~\ref{prop:reduction} is capacity-free.  It compares $\Fh$
with $G$ up to bounded factors, which move $\log\gamma_m^{\rm pair}$ by
$O(1)$ and so fix neither the geometry nor the prefactor.  In the proof
above it enters only through its second Douglas factor $C_2$, in
\cref{lem:S-coordinate} and \cref{lem:head-floor}.  Nor is the
logarithmic form itself available by a capacity-free route.  The
sharpest one proved here, Theorem~\ref{thm:general-twosided}, brackets
$\log\gamma_m^{\rm pair}$ between $2m\log\bigl(ekd/(2m)\bigr)$ and
$2m\log\bigl(ekR'/(2m)\bigr)$ up to $O(\log m)$, and the gap
$2m\log(R'/d)$ is of order $m$ for every admissible pair $d<R'$.  The
capacity enters only through the normal form of
Lemma~\ref{lem:faber-vekua}, as the leading Laurent coefficient $c$ of
$\Psi$ carried into the mode amplitudes $\mathcal A_m$ and hence into
$d_m$.  The proof is restricted to a single component.
For several components the normalized exterior map is replaced by
Walsh's map onto a lemniscatic domain, whose centers and component
weights enter the analogue of \eqref{eq:grunsky}; an extension would
require the corresponding Faber--Walsh normal form.  Boundaries with corners fall outside the normal
form---the annular extension of $\Psi$ and the banded classes
fail---though they retain the two-sided law of
Theorem~\ref{thm:general-twosided}.  The analysis here does not develop
an impedance analogue of the Dirichlet and Neumann monotonicity relations
in Albicker and Griesmaier~\cite{AG20}.
\end{remark}

\subsection{Reconstructions under the two-sided law}
In imaging terms, Theorem~\ref{thm:general-twosided} counts usable
modes.  At regularization $\alpha$, every obstacle supplies $N_\alpha$
reliable modes, and the theorem brackets this count additively,
$N(d)-O(1)\le N_\alpha\le N(R')+O(1)$ for each $R'>R$, and for a single
analytic
component Theorem~\ref{thm:capacity-connected} identifies the value:
$N_\alpha$ tracks $N(\operatorname{cap}(\overline D))$.
Table~\ref{tab:mode-budget} compares this prediction with the measured counts for
four of the obstacles of \cref{tab:capacity} at exact data and
$k=1$.  The three capacity-$1.30$ shapes pass \emph{identical} numbers of modes, to
within one (although their inradii range from $0.60$ to $1.30$) and track
$N(\operatorname{cap})$ to within one mode, while the capacity-$1.55$
shape passes two more modes at each level, against the one mode
predicted by $N(\operatorname{cap})$ and within the additive $O(1)$ of
\eqref{eq:general-count}.

\begin{table}[tbhp]
\centering
\caption{Measured number $N_\alpha$ of passed modes
(exact data, $k=1$) against the capacity prediction
$N(\operatorname{cap})$, for four of the obstacles of
Table~\ref{tab:capacity}.  Every count lies within one mode of $[N(d),N(R)]$.
Theorem~\ref{thm:general-twosided} itself gives the wider bracket
$[N(d),N(R')]$ for each $R'>R$.}
\label{tab:mode-budget}
\begin{tabular}{lccccccccc}
\toprule
& & & & \multicolumn{2}{c}{$\alpha=10^{-4}$}
& \multicolumn{2}{c}{$\alpha=10^{-7}$}
& \multicolumn{2}{c}{$\alpha=10^{-10}$}\\
shape & $\operatorname{cap}$ & $d$ & $R$
& $N_\alpha$ & $N(\operatorname{cap})$
& $N_\alpha$ & $N(\operatorname{cap})$
& $N_\alpha$ & $N(\operatorname{cap})$\\
\midrule
disc $a=1.30$          & $1.30$ & $1.30$ & $1.30$ & $9$  & $9.9$  & $13$ & $13.0$ & $15$ & $15.9$\\
ellipse $2.0\times0.6$ & $1.30$ & $0.60$ & $2.00$ & $9$  & $9.9$  & $13$ & $13.0$ & $15$ & $15.9$\\
four-lobed             & $1.30$ & $0.98$ & $1.95$ & $10$ & $9.9$  & $13$ & $13.0$ & $16$ & $15.9$\\
ellipse $2.0\times1.1$ & $1.55$ & $1.10$ & $2.00$ & $11$ & $10.8$ & $15$ & $14.1$ & $17$ & $17.1$\\
\bottomrule
\end{tabular}
\end{table}

\begin{figure}[tbhp]
\centering
\includegraphics{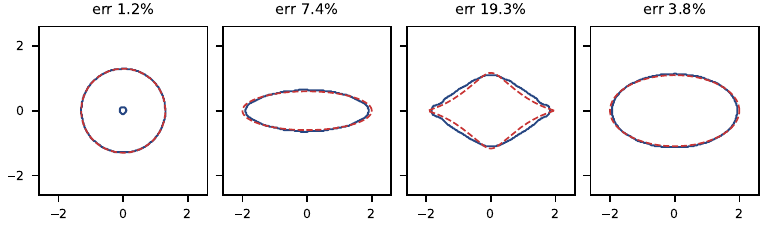}
\caption{Reconstructions of the four obstacles of
Table~\ref{tab:mode-budget} at $k=5$ from data with $2\%$ noise (first
test draw of each; solid blue: reconstructed boundary; dashed red: true
boundary), with the
calibrate-on-one-draw protocol of \cref{sec:num-kite}.}
\label{fig:general-recon}
\end{figure}

Figure~\ref{fig:general-recon} shows the corresponding reconstructions
at $k=5$ under $2\%$ data noise, with the protocol of
\cref{sec:num-kite}.  Over twelve test draws the misclassified areas are
$1.3\%\pm0.2\%$ for the disc, $7.7\%\pm0.3\%$ for the thin ellipse,
$19.1\%\pm0.3\%$ for the four-lobed star, and $4.0\%\pm0.2\%$ for the
larger ellipse.  The error ordering follows the geometry.  The
four-lobed star passes the same number of usable modes as the disc of
its capacity, but those modes must describe a more complicated shape.
Its reconstruction rounds the concave sides and tips, which lie beyond
the resolution scale of Remark~\ref{rem:resolution}.

\subsection{Numerical assessment of the capacity-informed
regularization level}
\label{sec:worked-example}

The prediction is tested on the ellipse with semiaxes $2.0$ and $1.1$
and on the kite, from $64\times64$ far-field matrices at $k=5$ with $2\%$
noise.  For the ellipse $\operatorname{cap}=(2.0+1.1)/2=1.55$ exactly, so
Theorem~\ref{thm:capacity-connected} supplies an a priori model for the
spectrum: the equal-capacity disc of radius $1.55$, whose
eigenvalues at $k=5$ are explicit Bessel expressions.  On noiseless
data the measured mode counts match that model exactly at every level
tested, $21$, $23$ and $25$ at $3.0\times10^{-2}$, $8.2\times10^{-3}$
and $10^{-3}$, although individual pairs agree only within a bounded
factor ($0.59$--$0.77$ at $m=12,\dots,18$).  Counts are protected by the
factorial gaps, prefactors are not.  At $2\%$ noise the rectified
per-mode floor (the median of twelve noisy spectra over the sorted
indices $j=45,\dots,60$, where the noise-free ladder is numerically
zero) is $8.2\times10^{-3}$, which leaves $23$ usable modes: the
monopole plus blocks $m=1,\dots,11$, so $23=1+2\cdot11$ with marginal
pair $\lambda_{22}=\lambda_{23}=2.73\times10^{-2}$, exhibiting the
$\lambda_{2m},\lambda_{2m+1}$ bookkeeping of
Theorem~\ref{thm:capacity-connected}.  The levels satisfying this
signal-retention criterion obey
$\alpha\le\lambda_{23}/19\approx1.4\times10^{-3}$, the condition that the
filter factor $q_{\alpha,j}=\lambda_j^2/(\lambda_j+\alpha)^2$ of the last
useful pair be at least $(19/20)^2=0.9025$.  Smaller $\alpha$ changes
the attenuation of these signal modes only modestly; this criterion does
not by itself control the contribution of noisy tail modes.  The
noise-proportional choice
$\alpha=\delta\opnorm F=3.0\times10^{-2}$ lies above this range.
Because $\delta$ is defined here through the Frobenius norm of the
discrete kernel, this is a heuristic noise scale, not an operator-norm
discrepancy estimate.
By the hard count $\lambda_j\ge\alpha$, it passes $21$ modes and places
the marginal useful pair just below the cutoff.  The conservative level
$\alpha^{*}=8.2\times10^{-2}$, ten times the floor, passes the same $21$.
At the noise-proportional level the Tikhonov filter does not remove that pair, but
attenuates it to weight $0.23$, and assigns near-unit weight,
$q_{\alpha,j}\ge0.95$, to only $10$ of the $23$ useful modes, against
$25$ modes at $\alpha=10^{-4}$.  At $\alpha^{*}$ the same pair is
attenuated further, to $0.06$.  The two thresholds are not the same
condition: $q_{\alpha,j}\ge0.95$ is strictly stronger than the retention
requirement $q_{\alpha,j}\ge(19/20)^2$, which every useful mode already
meets throughout this range.

Phantom calibration
selects $\alpha=10^{-4}$, within this range, with misclassified area
$4.0\%\pm0.2\%$ over twelve draws.  The noise-proportional level, granted its best
threshold in hindsight, gives $10.4\%$, a factor $2.6$ attributable
primarily to attenuation of the marginal pair, while a half-dynamic-range
threshold at the optimal $\alpha$ gives $106\%$.  A fixed fraction does
not locate the boundary on a slowly varying exterior field.  With the
fixed $19/20$ retention criterion, the capacity-informed level is the last
equal-capacity-disc eigenvalue above the measured floor divided by
$19$, giving $1.2\times10^{-3}$ and $4.3\%$ for the ellipse, within a
third of a percent of the calibrated value.  For the kite, where the
capacity has no closed form and is obtained from the equilibrium
problem on $\Gamma$, $\operatorname{cap}=1.3535$ prescribes
$1.0\times10^{-3}$ and yields $17.3\%$ against $23.9\%$ for the
noise-proportional choice.  Because the rule uses only a spectral level and
the error curve has a plateau, it tolerates error in the capacity.  A $6\%$
error moves $\alpha$ by a factor $2.3$ and the area error by $0.3$
points, against a $14$-point spread across
$\alpha\in[10^{-1},10^{-5}]$.  The kite further shows that robustness
across an inter-pair spectral gap depends on near-degeneracy.  Its $\pm m$
pairs split by
$12\%$ and $32\%$ at $j=18,19$ and $j=20,21$, leaving no empty band, so
the noise-proportional level falls inside a split pair and retains $20$ of the
$21$ useful modes.  Its residual $17.3\%$ is governed mainly by the
resolution in the data rather than by $\alpha$.  Exact data at the same
level give $16.6\%$ and
$\alpha\to0$ only $14.1\%$, against $3.7\%$ for the ellipse.  The
reconstruction is $8.8\%$ too large in area with $88\%$ of the excess
in the concave notch and errs by $9.8\%$ against the convex hull, the
notch depth $0.49$, measured from the convex hull, being comparable to
the resolvable scale
$|\Gamma|/(2m)=9.3/20\approx0.47$ at $21$ modes, $|\Gamma|$ the
perimeter.

The capacity-informed level also reduces spurious interior components.
At the noise-proportional level both reconstructions contain such components (one of $0.69\%$ of
the reconstructed area for the ellipse, three totalling $0.41\%$ for the
kite), which the capacity-informed level reduces to $0.12\%$ and removes
outright.  The spectral weights explain the difference.  That level carries $21$ and
$20$ modes at near-unit weight against $10$ in either noise-proportional
reconstruction, and the relative spread $(\max-\min)/\mathrm{mean}$ of
$W_\alpha$ over
the obstacle grows by some $70\%$, so noise no longer pushes isolated
cells across the threshold.
\section{Conclusions}

For an impedance disc, every positive-order H\"older source condition
fails at each noncentral interior sampling point
(\cref{prop:hoelder-fails}).  Relative to $\gamma_m^{1/2}$, the sampling
coefficients decay only geometrically, whereas an additional positive power
of $A$ requires factorial-scale decay.  Logarithmic source conditions provide
the appropriate replacement: the disc satisfies \cref{ass:log-source} for
every $\theta>0$ (\cref{prop:disc-log}).  These conditions yield upper
bounds, but do not determine the sharp bias rate or its dominant spectral
indices.

Inverting the sharpened eigenvalue asymptotics gives the transition index
$m_\alpha$ explicitly in terms of the Lambert $W$ function
(\cref{thm:disc-sharp}).  This index governs stabilization inside the disc,
iterated-logarithmic growth on the boundary, and divergence outside.  It also
determines the additive threshold in \cref{thm:disc-recon}.
Together with \cref{prop:disc-margin}, this gives upper bounds of
$O(1/m_\alpha)$ for the misclassified area and boundary offset.  Numerically,
$m_\alpha|D^{\tau_\alpha}_\alpha\triangle D|/|D|$ rises from $0.49$ at
$\alpha=10^{-10}$ and stays near $0.65$ for
$10^{-20}\ge\alpha\ge10^{-200}$.  Over the full range from $10^{-10}$ to
$10^{-200}$, the measured
symmetric difference decreases by only a factor of about seven, consistent
with the slow growth of $m_\alpha$.

For a general impedance obstacle, \cref{prop:reduction} compares the
spectral scale of $\Fh^{1/2}$ with the singular-value scale of the
data-to-pattern operator $G$ up to bounded factors.  The impedance remains
encoded in $G$ and the comparison constants.  Without analyticity,
\cref{thm:general-twosided} brackets the transition index by the inradius and
any radius exceeding the circumradius.  For an analytic boundary with
positive analytic impedance and conformal center inside the obstacle, the
Faber--Vekua normal form of \cref{lem:faber-vekua} identifies logarithmic
capacity as the effective radius, with the same leading constant as for the
disc.  The numerical results are consistent with this law: three obstacles
of capacity $1.30$ have mode counts differing by at most one at the tested
regularization levels.  Combining the capacity-based spectral model with the
measured noise floor also improves the reported reconstructions and reduces
spurious interior components.

Several extensions remain open.  The normal form does not cover boundaries
with corners or obstacles whose conformal center lies outside; the latter
failure is explained in \cref{rem:noncentered-obstruction}.  Multiple
boundary components are not treated; an extension would require a
Faber--Walsh normal form on a lemniscatic domain.  The Dirichlet and Neumann
inclusion analysis of Albicker and Griesmaier~\cite{AG20} is not extended here to impedance
boundaries.

\appendix
\renewcommand{\theHsection}{appendix.\Alph{section}}
\section{The conformal pull-back of the middle operator}
\label{app:pullback}

This appendix proves the normal form \eqref{eq:T-principal} used in the
proof of Lemma~\ref{lem:middle-pullback}, together with the banded matrix
estimates invoked there.  Throughout, $D$, $\Gamma$, $\Psi$, $\rho$, and
$\lambda$ are as in \cref{sec:general-capacity}, and all operators on
$\Gamma$ are transported to the parameter circle through
$\widetilde\varphi:=\varphi\circ\Psi|_{|w|=1}$.

\subsection{Banded matrix classes}

For $\sigma\in\mathbb R$ and $\rho_1\in(0,1)$ let
$\mathcal B^\sigma_{\rho_1}$ denote the class of block matrices
$A=(A_{mn})_{m,n\ge0}$ with
\[
 \opnorm{A_{mn}}\le C_A\,\rho_1^{|m-n|}\,(1+n)^{\sigma},
\]
and write $\mathcal B^\sigma:=\bigcup_{\rho_1<1}\mathcal B^\sigma_{\rho_1}$.
By the Schur test, $\mathcal B^{0}\subset\mathcal B(\ell^2)$.

\begin{lemma}[Composition and inverse localization]\label{lem:classes}
\begin{enumerate}
\item[\rm(i)] If $A\in\mathcal B^{\sigma_1}_{\rho_1}$ and
$B\in\mathcal B^{\sigma_2}_{\rho_1}$ with $\sigma_2\le0$, then
$AB\in\mathcal B^{\sigma_1+\sigma_2}_{\rho_2}$ for every
$\rho_2\in(\rho_1,1)$.  If $A$ carries the weight on the row
index, $\opnorm{A_{mn}}\le C\rho_1^{|m-n|}(1+m)^{-1}$, then so does
$AB$ for every $B\in\mathcal B^{0}_{\rho_1}$.
\item[\rm(ii)] If $X\in\mathcal B^{-1}_{\rho_1}$ and $I+X$ is
invertible on $\ell^2$, then $(I+X)^{-1}=I+Y$ with
$Y\in\mathcal B^{-1}_{\rho_2}$ for some $\rho_2\in(\rho_1,1)$.  If $X$
carries the weight on the row index, so does $Y$, that is
$\opnorm{Y_{mn}}\le C\rho_2^{|m-n|}(1+m)^{-1}$.
\item[\rm(iii)] $A\in\mathcal B^{\sigma}_{\rho_1}$ implies
$A^{*}\in\mathcal B^{\sigma}_{\rho_2}$ for every
$\rho_2\in(\rho_1,1)$.
\end{enumerate}
\end{lemma}

\begin{proof}
The entries of $AB$ are bounded by
$C_AC_B\sum_j\rho_1^{|m-j|+|j-n|}(1+j)^{\sigma_1}(1+n)^{\sigma_2}$ when
the weights sit as displayed.  Splitting the sum at $j=\min(m,n)/2$ and using, for the finitely many
dominant $j$, the bound
$(1+j)^{\sigma_1}\le C\,\rho_2^{-|m-j|}(1+m)^{\sigma_1}$, the geometric
convolution
$\sum_j\rho_1^{|m-j|+|j-n|}\le C\rho_2^{|m-n|}$ leaves the weight
$(1+m)^{\sigma_1}(1+n)^{\sigma_2}$, carried by the \emph{row} index
$m$ in its first factor.  One further transfer puts it on the column
index, as the class $\mathcal B^{\sigma_1+\sigma_2}$ requires.  Since
\[
 \frac{1+m}{1+n}\le1+|m-n|,
\]
one has
$(1+m)^{\sigma_1}\le(1+n)^{\sigma_1}\bigl(1+|m-n|\bigr)^{|\sigma_1|}$,
and the polynomial factor $(1+|m-n|)^{|\sigma_1|}$ is absorbed by a
final enlargement of the base.  This gives
$(1+n)^{\sigma_1+\sigma_2}$ and hence the claim.  A row weight
$(1+m)^{-1}$ multiplies every term of the $j$-sum by the same factor
and passes through unchanged, no transfer being needed.  The same
transfer makes each class closed under adjoints.  The entry $(A^{*})_{mn}$ is the
adjoint of $A_{nm}$ and therefore carries the weight on the row index,
so $A\in\mathcal B^{\sigma}_{\rho_1}$ implies
$A^{*}\in\mathcal B^{\sigma}_{\rho_2}$ for every
$\rho_2\in(\rho_1,1)$.

Part (ii) is Jaffard's theorem \cite{Jaf90}.  The inverse of a
boundedly invertible matrix with exponential off-diagonal decay
inherits that decay, so $(I+X)^{-1}\in\mathcal B^{0}_{\rho_2}$ for some
$\rho_2\in(\rho_1,1)$.  Since $X$ commutes with $(I+X)^{-1}$,
\[
 Y=-(I+X)^{-1}X=-X(I+X)^{-1}.
\]
The left-hand factorization is part~(i) with $\sigma_1=0$ and
$\sigma_2=-1$.  It supplies the gain of one order and the column weight
$Y\in\mathcal B^{-1}_{\rho_2}$.  The right-hand one has $X$ on the
left and $(I+X)^{-1}\in\mathcal B^{0}_{\rho_2}$ on the right, and
$\mathcal B^{0}_{\rho_1}\subset\mathcal B^{0}_{\rho_2}$, so the
row-weight clause of part~(i) applies at the base $\rho_2$ and carries
a row weight of $X$ over to $Y$.
\end{proof}

\begin{lemma}[Kernels to matrices]\label{lem:kernels}
Let $K(t,\tau)$ be a kernel on the parameter torus.
\begin{enumerate}
\item[\rm(i)] If $K$ extends holomorphically to
$|\Imag t|,|\Imag\tau|<\delta$, its matrix satisfies
$|K_{mn}|\le Ce^{-\delta(|m|+|n|)}$.  In particular it lies in
$\mathcal B^{\sigma}$ for every $\sigma$.
\item[\rm(ii)] If $K(t,\tau)=a(t,\tau)\,k_0(t-\tau)$ with $a$ analytic
on the torus and $k_0$ a $2\pi$-periodic convolution kernel whose
symbol obeys $|\widehat{k_0}(n)|\le C(1+|n|)^{\sigma}$, then
$K\in\mathcal B^\sigma$.
\end{enumerate}
\end{lemma}

\begin{proof}
(i) is Cauchy's estimate in each variable.  For (ii), expand
$a(t,\tau)=\sum_{p,q}\widehat a_{pq}e^{i(pt+q\tau)}$ with
$|\widehat a_{pq}|\le Ce^{-\delta(|p|+|q|)}$.  Then
$K_{mn}=\sum_{p+q=m-n}\widehat a_{pq}\,\widehat{k_0}(n+q)$, the
factor $e^{iq\tau}$ shifting the frequency that $k_0$ sees before
$e^{ipt}$ shifts the output; and
$(1+|n+q|)^{\sigma}\le(1+|n|)^{\sigma}(1+|q|)^{|\sigma|}$ lets the
exponential factors absorb the shift, so that
$|K_{mn}|\le C\rho_1^{|m-n|}(1+|n|)^{\sigma}$ with
$\rho_1=e^{-\delta/2}$.
\end{proof}

\begin{lemma}[A diagonal zero lowers the convolution order]
\label{lem:diagonal-zero}
Let $k_0$ be a periodic convolution kernel whose symbol $p_n$ satisfies
$|p_{n+1}-p_n|\le C(1+|n|)^{\sigma-1}$.  If $a$ is real analytic on the
torus and vanishes on the diagonal, $a(t,t)=0$, then
$a(t,\tau)k_0(t-\tau)\in\mathcal B^{\sigma-1}$.
\end{lemma}

\begin{proof}
Put $u=e^{it}$, $v=e^{i\tau}$.  Since $a$ vanishes on the diagonal,
division by $u-v$ is removable, so $a(t,\tau)=(e^{it}-e^{i\tau})
b(t,\tau)$ with $b$ real analytic on a slightly smaller complex
neighborhood of the torus.  Write
$b=\sum_{p,q}b_{pq}e^{i(pt+q\tau)}$, whose coefficients decay
exponentially in $|p|+|q|$.  Acting on $e_n$, the two terms carried by
$b_{pq}$ have the \emph{same} output frequency $n+p+q+1$, so they
combine, and their joint coefficient is the first difference
\[
 b_{pq}\bigl(p_{n+q}-p_{n+q+1}\bigr)
\]
in the Fourier convention of Lemma~\ref{lem:kernels}.  The hypothesis
on the differences, the exponential decay of $b_{pq}$, and the
geometric convolution of Lemma~\ref{lem:classes} give
$|K_{mn}|\le C\rho^{|m-n|}(1+|n|)^{\sigma-1}$ after enlarging $\rho<1$
once.
\end{proof}

\subsection{Kernel decompositions on an analytic curve}
\label{app:kernels}

Write $z(t):=\Psi(e^{it})$, $z'(t)=ie^{it}\Psi'(e^{it})$, and
\begin{equation}
 Q(t,\tau):=\frac{\Psi(e^{it})-\Psi(e^{i\tau})}{e^{it}-e^{i\tau}} .
\end{equation}
Since $\Gamma$ is analytic, $\Psi$ extends univalently to
$|w|>\rho'$ for some $\rho'<1$---$\Psi'$ is nonvanishing on $|w|=1$ and
injectivity on the circle propagates to a neighborhood---so $Q$ extends
holomorphically to a complex neighborhood of the torus and vanishes
nowhere.  Therefore $\log|Q|$ is analytic on the torus, and
$Q(t,t)=-ie^{-it}z'(t)$.

The Helmholtz kernel splits as
$\tfrac i4H_0^{(1)}(k|z|)=-\tfrac1{2\pi}J_0(k|z|)\log|z|+A(|z|^2)$ with
$A$ entire.  Hence, with $ds_\tau=|z'(\tau)|\dd\tau$, the transported
single-layer kernel is
\begin{equation}\label{eq:Skernel}
 \Phi_k\bigl(z(t),z(\tau)\bigr)
 =-\frac1{2\pi}\,J_0\bigl(k|z(t)-z(\tau)|\bigr)
 \Bigl(\log\bigl|2\sin\tfrac{t-\tau}2\bigr|+\log|Q(t,\tau)|\Bigr)
 +A_1(t,\tau),
\end{equation}
where $A_1$ and $\log|Q|$ are analytic on the torus, and
$J_0(k|z(t)-z(\tau)|)$ is analytic in $(t,\tau)$ because $J_0$ is even
and entire.  By Lemma~\ref{lem:kernels}, the transported single layer is
the circular logarithmic convolution, of symbol $(2\max(|n|,1))^{-1}$ up
to a rank-one term, times an analytic coefficient, plus analytically
smoothing terms.  It lies in $\mathcal B^{-1}$; its principal part is that
of the circular single layer multiplied by the analytic factor
$|\Psi'(e^{i\tau})|$.

The double-layer kernels of $K_k$ and $K_k'$ contain the factor
$\nu_\tau\cdot(z(t)-z(\tau))$, which vanishes to second order on the
diagonal.  On an analytic curve the Laplace parts are analytic kernels,
and the Helmholtz corrections are of the form covered by
Lemma~\ref{lem:kernels}(ii) with $\sigma\le-2$.  Hence
$K_k,K_k'\in\mathcal B^{0}$, and likewise the multiplications
$Z_k=ikM_\lambda$, with $\lambda$ analytic, lie in $\mathcal B^{0}$ by
Lemma~\ref{lem:kernels}(i) applied to $\lambda(t)\delta$-type
convolution, that is, their matrices are banded Toeplitz.

For the hypersingular part we use the integration-by-parts formula of
Maue~\cite{Maue49}, which converts the hypersingular form into two
single-layer forms carrying tangential derivatives; see also
\cite{CK19,Ned01}.  As a bilinear form on $H^{1/2}(\Gamma)$,
\begin{equation}\label{eq:maue}
 \ip{N_k\varphi}{\chi}_\Gamma
 =\int_\Gamma\!\!\int_\Gamma\Phi_k(x,y)\,
 \frac{d\varphi}{ds}(y)\,\overline{\frac{d\chi}{ds}}(x)\,\dd s_y\dd s_x
 -k^2\int_\Gamma\!\!\int_\Gamma\Phi_k(x,y)\,
 \nu(x)\!\cdot\!\nu(y)\,\varphi(y)\overline{\chi(x)}\,\dd s_y\dd s_x .
\end{equation}

\subsection{Proof of the normal form \texorpdfstring{\eqref{eq:T-principal}}{}}

Use the exact form pull-back $T_{\rm par}$ defined before
\eqref{eq:T-principal}.  The half-density conjugate
$|\Psi'|^{1/2}T|\Psi'|^{1/2}$, interpreted after pull-back, differs from it by commutators with the
analytic function $|\Psi'|^{1/2}$.  Each commutator kernel carries the
factor $|\Psi'(e^{it})|^{1/2}-|\Psi'(e^{i\tau})|^{1/2}$, which is
analytic on the torus and vanishes on the diagonal, so on each
convolution kernel of the decomposition of Appendix~\ref{app:kernels}
Lemma~\ref{lem:diagonal-zero} applies and lowers the order by one.  The
commutators therefore lie in $\mathcal B^{0}$, with
Lemmas~\ref{lem:classes} and~\ref{lem:kernels} handling the remaining
analytic products.
They are therefore included in $R_0$.  Thus $R_0$, and consequently
$X=T_{\rm circ}^{-1}R_0$ and $Y$, refer throughout to the exact operator
$T_{\rm par}$ used in \eqref{eq:assembly-form}.

In the first term of \eqref{eq:maue}, the substitution
$d/ds=|z'|^{-1}d/dt$ meets the arclength elements: with
$\dd s=|z'|\dd t$,
\[
 \int\!\!\int\Phi_k\,
 \frac{d\varphi}{ds}\,\overline{\frac{d\chi}{ds}}\,\dd s\,\dd s
 =\int\!\!\int\Phi_k\bigl(z(t),z(\tau)\bigr)\,
 \widetilde\varphi'(\tau)\,\overline{\widetilde\chi'(t)}\,\dd\tau\dd t,
\]
with \emph{no} conformal factor.  This is the cancellation used in
Lemma~\ref{lem:faber-vekua}.  Substituting \eqref{eq:Skernel}, the
leading kernel $-\tfrac1{2\pi}\log|2\sin\tfrac{t-\tau}2|$ produces
exactly the circular hypersingular form of radius $c$ after the two
derivatives---the radius enters only through the modes $g_{\pm m}$, not
through the form---while every remaining kernel is either analytic or of
the type of Lemma~\ref{lem:kernels}(ii) with $\sigma\le-1$.  After the
two derivatives these contribute classes $\mathcal B^{\sigma+2}$ with
$\sigma+2\le1$.  The order-one piece is
$a(t,\tau)L(t-\tau)$ with $L$ the circular logarithmic convolution and
$a:=J_0\bigl(k|z(t)-z(\tau)|\bigr)-1$, which is real analytic and
vanishes on the diagonal because $J_0(0)=1$.  The logarithmic symbol is
$O((1+|n|)^{-1})$ and its first difference is $O((1+|n|)^{-2})$, so
Lemma~\ref{lem:diagonal-zero} places $aL$ in $\mathcal B^{-2}$.  Two
derivatives then land it in $\mathcal B^{0}$.  The second term of \eqref{eq:maue} and the
remaining pieces $Z_kS_kZ_k^{*}$, $Z_kK_k$, $K_k'Z_k^{*}$ of $T$
transport to $\mathcal B^{0}$ (indeed $\mathcal B^{-1}$ for the first)
by Appendix~\ref{app:kernels} and Lemma~\ref{lem:classes}.  Collecting terms,
\[
 T_{\rm par}=T_{\rm circ}+R_0,
 \qquad R_0\in\mathcal B^{0},
\]
where $T_{\rm circ}$ is the middle operator of the circle of radius
$c$ with constant impedance $\lambda_0$.  The difference between the
transported impedance terms and the circle's is again in
$\mathcal B^{0}$.  This is \eqref{eq:T-principal}, with the banded
refinement used in the proof of Lemma~\ref{lem:middle-pullback}.  The symbol
of $T_{\rm circ}$ satisfies $|t_n|\ge c_T(1+|n|)$ for large $|n|$ by
the hypersingular principal part, and $t_n\ne0$ for every $n$ because
the circle's middle operator is invertible \cite{CCH13}, so
$X=T_{\rm circ}^{-1}R_0$ carries the row weight $(1+m)^{-1}$ of
Lemma~\ref{lem:classes}, as claimed.

The proof of Lemma~\ref{lem:middle-pullback} uses the
invertibility of $I+X$ on $L^2(\mathbb S^1)$, which
\eqref{eq:T-principal} does not give by substitution.  It gives the
statement one space higher.  Let $C_\Psi f=f\circ\Psi$ denote pull-back
to the circle.  From the defining form identity,
\[
 T_{\rm par}=M_{|\Psi'|}\,C_\Psi T C_\Psi^{-1}.
\]
The pull-back $C_\Psi$ is an isomorphism on $H^{\pm1/2}$, and
multiplication by the analytic nonvanishing function $|\Psi'|$ is an
isomorphism on $H^{-1/2}$.  Since
$T:H^{1/2}(\Gamma)\to H^{-1/2}(\Gamma)$ is an isomorphism, so is
$T_{\rm par}:H^{1/2}(\mathbb S^1)\to H^{-1/2}(\mathbb S^1)$.  Also
$T_{\rm circ}:H^{1/2}(\mathbb S^1)\to H^{-1/2}(\mathbb S^1)$ is
an isomorphism, its symbol obeying $|t_n|\asymp1+|n|$ and vanishing for
no $n$.  Hence \eqref{eq:T-principal}, read as
$T_{\rm circ}(I+X)=T_{\rm par}$, exhibits $I+X$ as an
isomorphism of $H^{1/2}(\mathbb S^1)$, in particular as an injection
there.  A Fredholm argument descends to $L^2(\mathbb S^1)$, identified
with $\ell^2$ through the Fourier basis.  A matrix of
$\mathcal B^{-1}$ differs in norm from its finite section of size $M$
by $O(1/M)$, by the Schur test applied to the two tail pieces, so $X$
is compact and $I+X$ is Fredholm of index zero.  If
$f\in L^2(\mathbb S^1)$ and $(I+X)f=0$, then $f=-Xf$, and the row
weight bounds the entries of $X$ by $C\rho^{|m-n|}(1+m)^{-1}$, so the
sequence $(1+m)(Xf)_m$ lies in $\ell^2$ by the Schur test.  Hence
$f\in H^{1}(\mathbb S^1)\subset H^{1/2}(\mathbb S^1)$ and $f=0$ by
the injectivity above.  Index zero and injectivity give the
invertibility of $I+X$ on $L^2(\mathbb S^1)$.  Only
\eqref{eq:T-principal}, the symbol bounds for $t_n$, the Schur test,
the composition rule of Lemma~\ref{lem:classes}(i), and the
isomorphism property of $T$ enter, so Lemma~\ref{lem:classes}(ii),
whose hypothesis is this invertibility, is not used.
\qed

\section*{Funding}
This research received no specific grant from any funding agency in the
public, commercial or not-for-profit sectors.

\bibliographystyle{amsplain-init}
\bibliography{references}

\begin{thebibliography}{27}
\providecommand{\natexlab}[1]{#1}
\providecommand{\url}[1]{\texttt{#1}}
\expandafter\ifx\csname urlstyle\endcsname\relax
  \providecommand{\doi}[1]{doi: #1}\else
  \providecommand{\doi}{doi: \begingroup \urlstyle{rm}\Url}\fi

\bibitem[Albicker and Griesmaier(2020)]{AG20}
A.~Albicker and R.~Griesmaier.
\newblock Monotonicity in inverse obstacle scattering on unbounded domains.
\newblock \emph{Inverse Problems}, 36:\penalty0 085014, 2020.
\newblock \doi{10.1088/1361-6420/ab98a3}.

\bibitem[Aleksandrov and Peller(2011)]{AP11}
A.~B. Aleksandrov and V.~V. Peller.
\newblock Estimates of operator moduli of continuity.
\newblock \emph{J. Funct. Anal.}, 261\penalty0 (10):\penalty0 2741--2796, 2011.
\newblock \doi{10.1016/j.jfa.2011.07.009}.

\bibitem[Arens(2004)]{Are04}
T.~Arens.
\newblock Why linear sampling works.
\newblock \emph{Inverse Problems}, 20:\penalty0 163--173, 2004.
\newblock \doi{10.1088/0266-5611/20/1/010}.

\bibitem[Arens and Lechleiter(2009)]{AL09}
T.~Arens and A.~Lechleiter.
\newblock The linear sampling method revisited.
\newblock \emph{J. Integral Equations Appl.}, 21:\penalty0 179--202, 2009.
\newblock \doi{10.1216/JIE-2009-21-2-179}.

\bibitem[Chamaillard et~al.(2014)Chamaillard, Chaulet, and Haddar]{CCH13}
M.~Chamaillard, N.~Chaulet, and H.~Haddar.
\newblock Analysis of the factorization method for a general class of boundary
  conditions.
\newblock \emph{J. Inverse Ill-Posed Probl.}, 22:\penalty0 643--670, 2014.
\newblock \doi{10.1515/jip-2013-0013}.

\bibitem[Colton and Kress(2019)]{CK19}
D.~Colton and R.~Kress.
\newblock \emph{Inverse Acoustic and Electromagnetic Scattering Theory},
  volume~93 of \emph{Applied Mathematical Sciences}.
\newblock Springer, Cham, 4 edition, 2019.
\newblock \doi{10.1007/978-3-030-30351-8}.

\bibitem[Corless et~al.(1996)Corless, Gonnet, Hare, Jeffrey, and
  Knuth]{CGHJK96}
R.~M. Corless, G.~H. Gonnet, D.~E.~G. Hare, D.~J. Jeffrey, and D.~E. Knuth.
\newblock On the {L}ambert {$W$} function.
\newblock \emph{Adv. Comput. Math.}, 5:\penalty0 329--359, 1996.
\newblock \doi{10.1007/BF02124750}.

\bibitem[Douglas(1966)]{Dou66}
R.~G. Douglas.
\newblock On majorization, factorization, and range inclusion of operators on
  {H}ilbert space.
\newblock \emph{Proc. Amer. Math. Soc.}, 17:\penalty0 413--415, 1966.
\newblock \doi{10.1090/S0002-9939-1966-0203464-1}.

\bibitem[Engl et~al.(1996)Engl, Hanke, and Neubauer]{EHN96}
H.~W. Engl, M.~Hanke, and A.~Neubauer.
\newblock \emph{Regularization of Inverse Problems}.
\newblock Kluwer, Dordrecht, 1996.
\newblock \doi{10.1007/978-94-009-1740-8}.

\bibitem[Gohberg and Krein(1969)]{GK69}
I.~C. Gohberg and M.~G. Krein.
\newblock \emph{Introduction to the Theory of Linear Nonselfadjoint Operators},
  volume~18 of \emph{Translations of Mathematical Monographs}.
\newblock American Mathematical Society, Providence, 1969.

\bibitem[Hofmann et~al.(2007)Hofmann, Kaltenbacher, P{\"o}schl, and
  Scherzer]{HKPS07}
B.~Hofmann, B.~Kaltenbacher, C.~P{\"o}schl, and O.~Scherzer.
\newblock A convergence rates result for {T}ikhonov regularization in {B}anach
  spaces with non-smooth operators.
\newblock \emph{Inverse Problems}, 23:\penalty0 987--1010, 2007.
\newblock \doi{10.1088/0266-5611/23/3/009}.

\bibitem[Hohage(1997)]{Hoh97}
T.~Hohage.
\newblock Logarithmic convergence rates of the iteratively regularized
  {G}auss--{N}ewton method for an inverse potential and an inverse scattering
  problem.
\newblock \emph{Inverse Problems}, 13:\penalty0 1279--1299, 1997.
\newblock \doi{10.1088/0266-5611/13/5/012}.

\bibitem[Hohage and Weidling(2015)]{HW15}
T.~Hohage and F.~Weidling.
\newblock Verification of a variational source condition for acoustic inverse
  medium scattering problems.
\newblock \emph{Inverse Problems}, 31:\penalty0 075006, 2015.
\newblock \doi{10.1088/0266-5611/31/7/075006}.

\bibitem[Hohage and Weidling(2017)]{HW17}
T.~Hohage and F.~Weidling.
\newblock Characterizations of variational source conditions, converse results,
  and maxisets of spectral regularization methods.
\newblock \emph{SIAM J. Numer. Anal.}, 55:\penalty0 598--620, 2017.
\newblock \doi{10.1137/16M1067445}.

\bibitem[Jaffard(1990)]{Jaf90}
S.~Jaffard.
\newblock Propri{\'e}t{\'e}s des matrices bien localis{\'e}es pr{\`e}s de leur
  diagonale et quelques applications.
\newblock \emph{Ann. Inst. H. Poincar{\'e} Anal. Non Lin{\'e}aire}, 7:\penalty0
  461--476, 1990.

\bibitem[Jung and Lim(2021)]{JL21}
Y.~Jung and M.~Lim.
\newblock Series expansions of the layer potential operators using the {F}aber
  polynomials and their applications to the transmission problem.
\newblock \emph{SIAM J. Math. Anal.}, 53:\penalty0 1630--1669, 2021.
\newblock \doi{10.1137/20M1348698}.

\bibitem[Kato(1973)]{Kato73}
T.~Kato.
\newblock Continuity of the map {$S\mapsto|S|$} for linear operators.
\newblock \emph{Proc. Japan Acad.}, 49:\penalty0 157--160, 1973.
\newblock \doi{10.3792/pja/1195519395}.

\bibitem[Kirsch(1998)]{Kir98}
A.~Kirsch.
\newblock Characterization of the shape of a scattering obstacle using the
  spectral data of the far field operator.
\newblock \emph{Inverse Problems}, 14:\penalty0 1489--1512, 1998.
\newblock \doi{10.1088/0266-5611/14/6/009}.

\bibitem[Kirsch and Grinberg(2008)]{KG08}
A.~Kirsch and N.~Grinberg.
\newblock \emph{The Factorization Method for Inverse Problems}.
\newblock Oxford University Press, Oxford, 2008.
\newblock \doi{10.1093/acprof:oso/9780199213535.001.0001}.

\bibitem[Ma and Hu(2022)]{MH22}
G.~Ma and G.~Hu.
\newblock Factorization method with one plane wave: from model-driven and
  data-driven perspectives.
\newblock \emph{Inverse Problems}, 38\penalty0 (1):\penalty0 015003, 2022.
\newblock \doi{10.1088/1361-6420/ac38b5}.

\bibitem[Math{\'e} and Pereverzev(2003)]{MP03}
P.~Math{\'e} and S.~V. Pereverzev.
\newblock Geometry of linear ill-posed problems in variable {H}ilbert scales.
\newblock \emph{Inverse Problems}, 19:\penalty0 789--803, 2003.
\newblock \doi{10.1088/0266-5611/19/3/319}.

\bibitem[Maue(1949)]{Maue49}
A.-W. Maue.
\newblock Zur {F}ormulierung eines allgemeinen {B}eugungsproblems durch eine
  {I}ntegralgleichung.
\newblock \emph{Z. Physik}, 126:\penalty0 601--618, 1949.
\newblock \doi{10.1007/BF01328780}.

\bibitem[Moiola et~al.(2011)Moiola, Hiptmair, and Perugia]{MHP11}
A.~Moiola, R.~Hiptmair, and I.~Perugia.
\newblock {V}ekua theory for the {H}elmholtz operator.
\newblock \emph{Z. Angew. Math. Phys.}, 62:\penalty0 779--807, 2011.
\newblock \doi{10.1007/s00033-011-0142-3}.

\bibitem[N{\'e}d{\'e}lec(2001)]{Ned01}
J.-C. N{\'e}d{\'e}lec.
\newblock \emph{Acoustic and Electromagnetic Equations: Integral
  Representations for Harmonic Problems}, volume 144 of \emph{Applied
  Mathematical Sciences}.
\newblock Springer, New York, 2001.
\newblock \doi{10.1007/978-1-4757-4393-7}.

\bibitem[Olver et~al.(2026)]{DLMF}
F.~W.~J. Olver et~al.
\newblock {NIST} digital library of mathematical functions, 2026.
\newblock \url{https://dlmf.nist.gov/}, accessed July 2026.

\bibitem[Ransford(1995)]{Ran95}
T.~Ransford.
\newblock \emph{Potential Theory in the Complex Plane}, volume~28 of
  \emph{London Mathematical Society Student Texts}.
\newblock Cambridge University Press, Cambridge, 1995.

\bibitem[Tautenhahn(1998)]{Tau98}
U.~Tautenhahn.
\newblock Optimality for ill-posed problems under general source conditions.
\newblock \emph{Numer. Funct. Anal. Optim.}, 19:\penalty0 377--398, 1998.
\newblock \doi{10.1080/01630569808816834}.

\end{thebibliography}

\end{document}